%% file: integral_7_Sept2014_2.tex
\documentclass[12pt] {article}
\usepackage{amsmath,amsthm}
\usepackage{amssymb,latexsym}
\usepackage{amscd}
\usepackage{epic,eepic}
\title{Valuations on manifolds and integral geometry.}
\date{}
\author{ Semyon Alesker \footnote{Partially supported by ISF grant 701/08.}
\\  { \normalsize Department of Mathematics, Tel Aviv University, Ramat Aviv}
 \\  { \normalsize 69978 Tel Aviv,
Israel }
\\ {\normalsize e-mail: semyon@post.tau.ac.il}}

\def\eps{\varepsilon}
\def\alp{\alpha}
\def\ome{\omega}
\def\Ome{\Omega}
\def\lam{\lambda}
\def\Lam{\Lambda}

\def\to{\rightarrow}
\def\qed { Q.E.D. }
\def\pt{\partial}
\def\grc{{}\!^ {\textbf{C}} Gr}
\def\grk{{}\!^ {\mathbb{K}} Gr}

\def\RR{\mathbb{R}}
\def\CC{\mathbb{C}}

\def\NN{\mathbb{N}}
\def\ZZ{\mathbb{Z}}
\def\HH{\mathbb{H}}

\def\PP{\mathbb{P}}
\def\KK{\mathbb{K}}

\def\One{{1\hskip-2.5pt{\rm l}}}
\swapnumbers
\newtheorem{theorem}{Theorem}[subsection]
\newtheorem{corollary}[theorem]{Corollary}
\newtheorem{lemma}[theorem]{Lemma}
\newtheorem{proposition}[theorem]{Proposition}
\newtheorem{claim}[theorem]{Claim}
\theoremstyle{definition}
\newtheorem{example}[theorem]{Example}
\newtheorem{definition}[theorem]{Definition}
\newtheorem{remark}[theorem]{Remark}

\theoremstyle{proposition-definition}
\newtheorem{proposition-definition}[theorem]{Proposition-Definition}

\numberwithin{equation}{subsection}

\def\cf{{\cal F}}

\def\ca{{\cal A}}  \def\cc{{\cal C}}
\def\cd{{\cal D}} \def\ce{{\cal E}} \def\cf{{\cal F}}
  
  \def\cl{{\cal L}}
\def\cm{{\cal M}} \def\cn{{\cal N}} \def\co{{\cal O}}
\def\cp{{\cal P}} 

\def\cs{{\cal S}} \def\ct{{\cal T}} \def\cu{{\cal U}}
\def\cv{{\cal V}}  
 
\def\inj{\hookrightarrow }

\def\vi{V^\infty}
\def\vmi{V^{-\infty}}

\input diagram
\begin{document}
\maketitle
\begin{abstract}
One constructs new operations of pull-back and push-forward on
valuations on manifolds with respect to submersions and immersions.
A general Radon type transform on valuations is introduced using
these operations and the product on valuations. It is shown that the
classical Radon transform on smooth functions, and the well known
Radon transform on constructible functions with respect to the Euler
characteristic are special cases of this new Radon transform. An
inversion formula for the Radon transform on valuations has been
proven in a specific case of real projective spaces. Relations of
these operations to yet another classical type of integral geometry,
Crofton and kinematic formulas, are indicated.
\end{abstract}

\tableofcontents \setcounter{section}{-1}

\section{Introduction.}\label{S:introduction}

\subsection{A general overview of the article.}\label{Ss:general-overview}

The theory of valuations on manifolds was introduced in
\cite{part1}, \cite{part2}, \cite{part3}, \cite{part4}; see also the
survey of these results \cite{alesker-survey}. This theory
generalizes in some direction the classical theory of valuations on
convex subsets of $\RR^n$ (see e.g. the surveys McMullen-Schneider
\cite{mcmullen-schneider} and McMullen \cite{mcmullen-survey}).

Let us describe a basic idea behind this new notion of a valuation
on a manifold in an informal way. Valuation is a finitely additive
complex valued functional defined on "nice" subsets of a smooth
manifold $X$ with some extra properties to be discussed.  By "nice"
subsets we usually mean compact submanifolds with corners, while in
the classical theory of valuations one considers the class of finite
unions of convex compact subsets of $X=\RR^n$.

The class of all finitely additive measures is too large and is
probably out of control. In order to have a more restricted class of
measures which still covers interesting examples from geometry and
analysis, one could try to look for an extra condition of analytic
nature which would tell how the measures of sets behave with respect
to limits. In the classical measure theory of Lebesgue this
condition is the countable additivity. As everybody knows, this
condition turned out to be extremely useful in numerous situations.

In the new theory of valuations on manifolds we impose a rather
different analytic condition (some version of continuity or
smoothness). While on a general manifold this condition is somewhat
technical (see \cite{part2}), it is basically motivated by the
classical theory of valuations on convex sets in combination with
some ideas from geometric measure theory: in the case of $X=\RR^n$
the main condition is the continuity of the valuation $\phi$ on the
class of {\itshape convex} compact sets with respect to the
Hausdorff metric, namely if $\{K_k\}$ is a sequence of convex
compact sets converging in the Hausdorff metric to another convex
compact set $K$ then $\phi(K_k)\to \phi(K)$. To the best of our
knowledge, this condition of continuity was first introduced and
systematically studied by Hadwiger, see his book
\cite{hadwiger-book}. Notice here also that results of geometric
measure theory imply that the right generalization of this
continuity condition to a more general class of non-convex sets is
as follows: if $\{A_k\}$ is a sequence of "nice" compact subsets of
a manifold $X$ such that the sequence of their normal cycles
$\{N(A_k)\}$ (see Section \ref{Ss:normal-cycles} below) converges in
the flat topology on currents to the normal cycle of a set $A$, then
$\phi(A_k)\to \phi(A)$. On convex sets this kind of convergence is
equivalent to the convergence in the Hausdorff metric.

\hfill

We would like to emphasize that valuations may not be countably
additive. Nevertheless smooth countably additive measures (i.e.
those which locally can be written as $f(x_1,\dots,x_n) dx_1\dots
dx_n$ with the function $f$ being $C^\infty$-smooth) are simplest
examples of valuations in the new theory. Usually valuations cannot
be defined on too broad class of sets, say on Borel sets. The most
general class of sets for this theory is not very clear for the
moment. Nevertheless if the manifold $X$ is real analytic then any
smooth valuation can be naturally evaluated on any subanalytic
relatively compact subset. It seems very likely that on a general
smooth manifold $X$ valuations can naturally be evaluated also on
compact subsets of positive reach (though a detailed proof has not
appeared anywhere).

Another most basic example of a valuation is the Euler
characteristic; definitely it is not countably additive, and it is
not defined on say Borel sets. The classical theory of valuations on
convex sets provides us with many more examples of valuations of
geometric importance. Let us consider the Euclidean plane $\RR^2$.
The length of the boundary of a set is an example of a smooth
valuation. In particular, it is defined on all compact submanifolds
with corners and is continuous with respect to the Hausdorff metric
on convex compact subsets.

More generally let us fix in $\RR^n$ convex compact domains
$A_1,\dots,A_k$, $1\leq k\leq n-1$,  with smooth positively curved
boundaries. Consider the mixed volume
\begin{eqnarray}\label{09oct1}
\phi(K):=V(\underset{n-k \mbox{
times}}{\underbrace{K,\dots,K}},A_1,\dots,A_{n-k}).
\end{eqnarray}
Classically it is defined only for convex sets $K$ (see e.g.
\cite{schneider-book}), but in fact it can be naturally extended to
all compact submanifolds with corners in $\RR^n$, and this extension
is a smooth valuation. For example if we take $A_1=\dots =A_k$ to be
the unit Euclidean ball then the valuation defined by \ref{09oct1}
is (up to a constant) the so called $(n-k)$-th intrinsic volume (in
another terminology it is also called the Lipschitz-Killing
curvature). In particular taking $n=2,\, k=1$, and $A_1$ being the
Euclidean ball, we recover the length of the boundary of a subset of
$\RR^2$. There are also some other natural sources of valuations,
for example coming from integral geometry.

Let us also mention that the space of translation invariant
(continuous) valuations on convex subsets of $\RR^n$ is a very
classical object; this space almost coincides with the space of
translation invariant valuations in the new setting of manifolds
(more precisely, the former space contains the latter as a dense
subspace). For $n\geq 2$ this space is infinite dimensional: it
contains not only the Lebesgue measure, but also the Euler
characteristic and all the mixed volumes. Nevertheless this space is
not out of control, and much is known about it \cite{hadwiger-book},
\cite{mcmullen-schneider}, \cite{mcmullen-survey},
\cite{alesker-gafa01}, \cite{alesker-fourier}.

\hfill


Let us denote by $\vi(X)$ the space of smooth valuations on a
manifold $X$ (see Section \ref{Ss:valuations-mfld} for a formal
definition). It is naturally a Fr\'echet space. It has rich
structures (see the survey \cite{alesker-survey}). In this article
the multiplicative structure on $\vi(X)$ will be of a particular
importance. The space $\vi(X)$ has a canonical product
$$\vi(X)\times \vi(X)\to \vi(X)$$
which is continuous, and $\vi(X)$ becomes a commutative associative
algebra; the Euler characteristic $\chi$ is its unit element. The
product on valuations was introduced first in \cite{alesker-gafa04}
for so called polynomial valuations on convex subsets of $\RR^n$,
then extended to any smooth valuations on $\RR^n$ in \cite{part1},
and eventually in \cite{part3} it was shown that it is independent
of the affine structure on $\RR^n$, and hence extends to any
manifold. A rather different construction of this product has been
recently given in \cite{alesker-bernig}.


The first main goal of this article is to introduce the operations
of pull-back and push-forward on valuations under smooth maps of
manifolds. This is done under rather restrictive assumptions on the
maps: they are assumed to be either submersions or immersions (more
general cases will be discussed elsewhere). Nevertheless these cases
are sufficient for some applications to integral geometry.

The second main goal of the article is to introduce a general Radon
type transform on valuations. The construction of it involves the
product on valuations and the operations of pull-back and
push-forward mentioned above. We show that the classical Radon
transform on {\itshape smooth} functions can be considered as a very
special case of this new Radon transform (see Section
\ref{Ss:gelfand}). Moreover this new Radon transform generalizes
(partly conjecturally) the well known and seemingly completely
different Radon transform on {\itshape constructible} functions with
respect to the Euler characteristic (see Section
\ref{Ss:constructible}).

Eventually we prove in Theorem \ref{T:radon-val-constr} an inversion
formula for this new Radon transform on smooth valuations in a very
concrete situation, generalizing the Khovanskii-Pukhlikov inversion
formula \cite{khovanskii-pukhlikov} for the Radon transform with
respect to the Euler characteristic on constructible functions on
the real projective space $\RR\PP^n$.

In Section \ref{Ss:chern} we discuss the relations of the product,
pull-back, and push-forward on valuations to yet another classical
and quite different type of integral geometry: kinematic and Crofton
type formulas.

Thus valuations on manifolds provide a general set up to unify
different directions of integral geometry.

\subsection{A more detailed overview of the main
results.}\label{Ss:detailed-overview} Let us denote by $\cp(X)$ the
family of all compact submanifolds with corners of a manifold $X$.
Smooth valuations are finitely additive functionals $\cp(X)\to \CC$
satisfying some additional conditions (see Section
\ref{Ss:valuations-mfld}).

We have to recall few more general properties of valuations on
manifolds (see e.g. the survey \cite{alesker-survey}). We denote by
$\vi_c(X)$ the subspace of $\vi(X)$ of compactly supported
valuations. By \cite{part4} we have the integration functional
$$\int_X\colon \vi_c(X)\to \CC$$
given by $\int_X\phi:=\phi(X)$. This is a linear functional which is
continuous in appropriate topology on $\vi_c(X)$.

We have a continuous bilinear map
$$\vi_c(X)\times \vi(X)\to \CC$$
defined by $(\phi,\psi)\mapsto \int_X\phi\cdot \psi$ where under the
integral we have the product on valuations mentioned above. By
\cite{part4} (see also \cite{bernig-quat} for a simpler proof) this
bilinear map is a perfect pairing, i.e. the induced map
\begin{eqnarray}\label{E:poincare}
\vi(X)\to (\vi_c(X))^*
\end{eqnarray}
is injective and has a dense image in the weak topology. Thus we
denote by $\vmi(X):=(\vi_c(X))^*$ and call it the space of {\itshape
generalized valuations}. Via the map (\ref{E:poincare}) we identify
the space $\vi(X)$ of smooth valuations with its image in the space
$\vmi(X)$ of generalized valuations:
$$\vi(X)\subset \vmi(X).$$


\hfill

Let us describe now the relations between valuations and
constructible functions. In the case of real analytic manifold $X$
this relation can be stated more elegantly. We will explain this
case after considering first the more general case of a smooth
manifold.

For a smooth manifold $X$ consider the natural map
$$\Xi_{\cp}\colon \cp(X)\to \vmi(X)$$
given by $P\overset{\Xi_{\cp}}{\mapsto}[\phi\mapsto \phi(P)]$. This
map is injective and the span of its image is dense in
$V^{-\infty}(X)$ in the weak topology. Thus for any smooth manifold
$X$ we have
\begin{eqnarray}\label{E:smooth-constr}
V^\infty(X)\subset V^{-\infty}(X)\supset \cp(X).
\end{eqnarray}
In the special case when $X$ is a real analytic manifold it is more
elegant to consider instead of $\cp(X)$ the space $\cf(X)$ of
constructible functions. By definition $\cf(X)$ consists of
functions $f\colon X\to \CC$ which are finite linear combinations
with complex coefficients of integer valued functions $g\colon X\to
\ZZ$ such that for any $m\in \ZZ$ the set $g^{-1}(m)$ is
subanalytic, and the family of sets $\{g^{-1}(m)\}_{m\in\ZZ}$ is
locally finite.\footnote{Essentially this class of constructible
functions was studied in detail in Ch. 9 of
\cite{kashiwara-schapira} where also the definition of a subanalytic
subset can be found.} $\cf(X)$ is an algebra over $\CC$ with the
pointwise product. Any smooth valuation can be naturally evaluated
on a relatively compact subanalytic set, even if it is not a compact
submanifold with corners (see \cite{part4}). Hence we have the map
(also denoted by $\Xi_{\cp}$)
$$\Xi_{\cp}\colon\cf(X)\to V^{-\infty}(X)$$ given by
$\sum_i\alp_i\One_{P_i}\overset{\Xi_{\cp}}{\mapsto}
[\phi\mapsto \sum_i\alp_i\phi(P_i)]$. By \cite{part4}, Section 8.1,
this map is injective and has a dense image in the weak topology.
Thus for a real analytic manifold $X$ we have imbeddings of two
linear dense subspaces analogous to (\ref{E:smooth-constr})
\begin{eqnarray}\label{E:subanalytic-constr}
V^\infty(X)\subset V^{-\infty}(X)\supset \cf(X).
\end{eqnarray}

The imbeddings (\ref{E:smooth-constr}), (\ref{E:subanalytic-constr})
are often very useful since various structures on valuations can be
restricted to $\cp(X)$ or, better to constructible functions
$\cf(X)$, where they can be interpreted in more familiar terms. For
example the above mentioned product on $\vi(X)$ has the following
interpretation. It was shown in \cite{alesker-bernig} that $\vmi(X)$
has a partially defined product, i.e. two generalized valuations can
be multiplied if they are in "generic position" to each other (the
precise technical conditions are formulated on the language of wave
front sets). For two smooth valuations this condition is satisfied
automatically, and the product coincides with the above mentioned
product on $\vi(X)$. On the other hand, the restriction of this
partial product to $\cp(X)$ is just the usual intersection of sets:
more precisely given $P_1,P_2\in \cp(X)$ which are transversal to
each other then the product of their images in $V^{-\infty}(X)$ is
well defined and is equal to the image of $P_1\cap P_2$; this was
proved in \cite{alesker-bernig}. It is expected that for real
analytic $X$ the restriction of the partial product on
$V^{-\infty}(X)$ to $\cf(X)$ coincides with the pointwise product on
constructible functions provided the functions are in "generic
position" to each other.
\begin{remark}
When talking on the map $\Xi_{\cp}\colon \cp(X)\to V^{-\infty}(X)$
for general smooth manifold $X$ we will often identify $P\in \cp(X)$
with its indicator function $\One_P$ and thus write
$\Xi_{\cp}(\One_P)$.
\end{remark}

\hfill

After this reminder let us discuss the main results of this article.
Let us discuss the pull-back and push-forward of valuations. In this
introduction the discussion will be not completely rigorous and
partly conjectural. The precise, but more technical, results can be
found in the main text of the article. We will try however to
indicate what is rigorous and what is not.

Let $f\colon M\to N$ be a smooth proper map of smooth manifolds.
Under appropriate assumptions one should be able to define a linear
map
$$f_*\colon \vi(M)\to \vmi(N)$$
called push-forward. The main idea is that $f_*$ should satisfy the
usual (in comparison to the classical measure theory) property
$$(f_*\phi)(P)=\phi(f^{-1}(P))$$
for any $P\in \cp(N)$. On the technical level this is not well
defined of course since the set $f^{-1}(P)$ may not belong to
$\cp(M)$. Nevertheless if $f$ is a {\itshape submersion} we get a
well defined map $f_*\colon \vi(M)\to\vi(N)$ (see Section
\ref{Ss:push-submersion}). For immersions the push-forward
$f_*\colon\vmi(M)\to \vmi(N)$ is defined rigorously in Section
\ref{Ss:push-gener-immer}. One expects that $(f\circ g)_*=f_*\circ
g_*$.

The pull-back map $f^*\colon \vi(N)\to \vmi(M)$ should be defined as
the dual of $f_*$. Rigorously this has been done for immersions in
Section \ref{Ss:pull-immersion}, and for submersions in Section
\ref{Ss:submersions}. One expects that $(f\circ g)^*=g^*\circ f^*$.
Actually in order to define $f^*$ one does not have to assume that
the map $f$ is proper.

Let us discuss these operations in some examples. First recall that
smooth measures form a subspace of smooth valuations. Then the
restriction of $f_*$ to smooth measures should coincide with the
classical push-forward of measures; we have proven this when $f$ is
either submersion or immersion.

Under appropriate assumptions $f_*$ and $f^*$ should be extended by
(sequential) continuity to larger subspaces of {\itshape generalized} valuations.
Then we could restrict them to constructible functions (which are,
perhaps, in "generic position" to the map $f$). Then on
constructible functions $f^*$ should coincide with the usual
pull-back of functions (this is proved rigorously for submersions in
Section \ref{Ss:submersions} and for closed imbeddings in
Proposition \ref{P:pull-back-constuct}).

The push-forward $f_*$ should coincide with the less trivial
operation of integration with respect to the Euler characteristic
along the fibers. It means the following: if $\One_P$ denotes the
indicator function of a set $P\in \cp(M)$ then
$$(f_*(\One_P))(y)=\chi(P\cap f^{-1}(y))$$
for any $y\in N$. If $P$ is in "generic position" with respect to
$f$ then $f_*(\One_P)$ is a constructible function. This fact seems
to be particularly technical to prove in full generality. We have
proven it under quite restrictive assumptions (see Section
\ref{Ss:pushforward-gen-submersions} and the proof of Proposition
\ref{P:radon-convex}). This fact is necessary for applications in
integral geometry in Section \ref{Ss:constructible}.

Notice that pull-back of the Euler characteristic should be again
the Euler characteristic whenever the pull-back is defined. We would
like to have one more remark on $f^*$ when $f\colon M\to N$ is a
submersion. In this case we have shown in Section
\ref{Ss:submersions} that $f^*\colon \vmi(N)\to \vmi(M)$ is a
continuous linear map. Let us restrict this map to smooth measures.
The classical measure theory tells us that on measures there is no
operation of pull-back. The new thing is that we can define the
pull-back of a measure to be a (generalized) valuation. With an
oversimplification, we have for a smooth measure $\mu$ on $N$:
$$(f^*\mu)(P)=\int_N\chi(P\cap f^{-1}(y))d\mu(y).$$

\hfill

Let us discuss now our applications of these constructions to
integral geometry. A double fibration is a diagram
$$X\overset{q_1}{\leftarrow}Z\overset{q_2}{\to} Y$$
of smooth submersive maps of smooth manifolds such that the map
$q_1\times q_2\colon Z\to X\times Y$ is a closed imbedding. We
assume in addition that $q_2$ is proper. Let us fix a smooth
valuation $\gamma\in \vi(Z)$. We define a map which we call the
Radon transform on valuations
$$R_\gamma\colon \vi(X)\to \vmi(Y)$$
by $R_\gamma(\phi)=q_{2*}(\gamma\cdot q_1^*\phi)$ where the
pull-back, push-forward, and the product are understood in the sense
of valuations. Theorem \ref{T:radon-val-1} says that $R_\gamma$ is a
well defined continuous linear operator. Here is a technical
difficulty: though $\phi$ is a smooth valuation, the pull-back
$q_1^*\phi$ is not smooth in general. Hence in order to define
$R_\gamma(\phi)$ in general it was necessary to extend (partially)
the push-forward map to generalized valuations.

Moreover Corollary \ref{COR:radon-val-2} says that under certain
assumption on the double fibration, $R_\gamma$ takes values in
smooth valuations $\vi(Y)$ for any $\gamma\in \vi(Z)$, and it is a
continuous linear map
$$\vi(X)\to\vi(Y).$$
Also under a similar assumption $R_\gamma$ extends (uniquely) by
continuity to a linear map on generalized valuations
$$\vmi(X)\to \vmi(Y).$$

Next we show that in the case when $\gamma$ is a smooth measure
considered as a smooth valuation, the above Radon transform
$R_\gamma$ can be considered as the classical Radon transform on
smooth functions. More precisely, $R_\gamma$ uniquely factorizes in
the following way:
$$\Vtriangle<1`3`-1;500>[\vi(X)`\vmi(Y)\hookleftarrow\cm^\infty(Y)`C^\infty(X); R_\gamma``]$$
where the surjection $\vi(X) \twoheadrightarrow C^\infty(X)$ is the
evaluation-on-points map, i.e.
$$\phi\mapsto[x\mapsto\phi(\{x\})],$$
and the map $C^\infty(X)\to\cm^\infty(Y)$ is the classical Radon
transform from the space $C^\infty(X)$ of smooth functions to the
space $\cm^\infty(Y)$ of smooth measures (see Section
\ref{Ss:gelfand}).

Let us consider another special case when $\gamma=\chi$ is the Euler
characteristic. Assume that $R_\chi$ extends by continuity to a
linear map $\vmi(X)\to\vmi(Y)$ (which is the case under appropriate
assumptions; see Corollary \ref{COR:radon-val-2}). Assume in
addition that $X,Y,Z,q_1,q_2$ are real analytic. Then we conjecture
that $R_\chi(\cf(X))\subset \cf(Y)$ (at least for "generic"
constructible functions) and the map $R_\chi\colon \cf(X)\to \cf(Y)$
coincides with the Radon transform with respect to the Euler
characteristic. Such kind of the Radon transform was considered by
several authors under various modifications of the notion of
constructibility: see \cite{khovanskii-pukhlikov},
\cite{khovanskii-pukhlikov2}, \cite{schapira2}; in the earlier paper
\cite{viro} the constructibility was considered in the sense of
{\itshape complex} analytic geometry. We have not proven this
conjecture in full generality, but some special cases of it were
proven and used in Section \ref{Ss:constructible}.

In Section \ref{Ss:chern} we describe relations of the product,
pull-back, and push-forward on valuations to a very different type
of integral geometry: kinematic and Crofton type formulas.

\hfill

Finally let us describe our last main result: an inversion formula
for the Radon transform on smooth valuations in a very concrete
situation. Let $\RR\PP^n$ denote the real projective space, i.e. the
space of lines through the origin in $\RR^{n+1}$. Let
$\RR\PP^{n\vee}$ denote the dual projective space, i.e. the space of
hyperplanes through the origin in $\RR^{n+1}$. Let $Z$ be the
incidence variety, namely
$$Z:=\{(l,H)\in \RR\PP^n\times \RR\PP^{n\vee}|\, l\subset H\}.$$
Let
$\RR\PP^n\overset{q_1}{\leftarrow}Z\overset{q_2}{\to}\RR\PP^{n\vee}$
be the double fibration where $q_1,q_2$ are the obvious projections.
We have the Radon transform with respect to the Euler characteristic
\begin{eqnarray}\label{E:march-radon}
R_\chi\colon \vi(\RR\PP^n)\to \vi(\RR\PP^{n\vee})
\end{eqnarray}
defined as above, i.e. (taking into account that $\chi$ is the unit
element) $R_\chi=q_{2*}q_1^*$. Then $R_\chi$ is invertible up to a
multiple of the Euler characteristic. More precisely consider the
dual Radon transform
$$R^t_\chi\colon \vi(\RR\PP^{n\vee})\to \vi(\RR\PP^n)$$
defined similarly by $R_\chi^t:=q_{1*}q_2^*$. Then we prove (Theorem
\ref{T:radon-val-constr}) that for any $\phi\in\vi(\RR\PP^n)$ one
has
\begin{eqnarray}\label{E:march-invert}
(-1)^{n-1}(R^t_\chi\circ
R_\chi)\phi=\phi+\frac{1}{2}((-1)^{n-1}-1)\left(\int_{\RR\PP^n}\phi\right)\cdot
\chi.
\end{eqnarray}
This inversion formula was motivated by the inversion formula for
the Radon transform with respect to the Euler characteristic on
"constructible" functions due to Khovanskii and Pukhlikov
\cite{khovanskii-pukhlikov}. \footnote{Actually their notion of
constructibility was more restrictive one than what we have
described above.} The same inversion formula was generalized by
Schapira \cite{schapira2} to functions constructible in the
subanalytic sense.
\begin{remark}
Our proof of the inversion formula (\ref{E:march-invert}) on smooth
valuations uses in fact the above mentioned result by Khovanskii and
Pukhlikov and generalizes it in the following sense. The Radon
transform $R_\chi$ can be extended by continuity to a linear map on
generalized valuations
$$R_\chi\colon \vmi(\RR\PP^n)\to \vmi(\RR\PP^{n\vee}).$$
Similarly $R^t_\chi$ extends by continuity to generalized
valuations.

If we restrict $R_\chi$ to the subspace of constructible functions
in the sense of Khovanskii-Pukhlikov (i.e. finite linear
combinations of indicator functions of convex compact polytopes in
$\RR\PP^n$) then we show that $R_\chi$ coincides with the Radon
transform with respect to the Euler characteristic. Since such
functions are also dense in $\vmi(\RR\PP^n)$, it suffices to prove
the inversion formula (\ref{E:march-invert}) for them, and this was
done by Khovanskii-Pukhlikov \cite{khovanskii-pukhlikov}. We
conjecture that the restriction of $R_\chi$ to the broader class of
subanalytic constructible functions also coincides with the Radon
transform with respect to the Euler characteristic on this class.
This would generalize then (a special case of) the result by
Schapira \cite{schapira2} in the same way as we have generalized the
result of Khovanskii-Pukhlikov.
\end{remark}

\hfill

{\bf Acknowledgements.} I thank J. Bernstein for numerous useful
discussions, and V. Milman for his attention to this work. I thank
also V. Palamodov and S. Schochet for their explanations on wave
front sets. I thank R. Schneider for some references.

\subsection{Notation list.}\label{Ss:notation}

$\bullet$ $V^\infty(X)$ - the space of smooth valuations on a
manifold $X$.

$\bullet$ $V^{-\infty}(X)$ - the space of generalized valuations on
a manifold $X$.

$\bullet$ $\cf(X)$ - the space of constructible functions on $X$.


$\bullet$ $\cp(X)$ - the family of compact submanifolds with corners
in $X$.

$\bullet$ $R_\gamma$ - the Radon transform on valuations with the
kernel $\gamma$.

$\bullet$ $\PP_+(V)$ - the oriented projectivization of a real
vector space (or a vector bundle) $V$, i.e.
$(V\backslash\{0\})/\RR_{>0}$.

$\bullet$ $\Ome^k$ - either the vector bundle or the space of smooth
$k$-forms.

$\bullet$ $|\ome_X|$ - the line bundle of densities on $X$.

$\bullet$ $\cm^\infty(X)$ - the space of smooth measures (densities)
on a manifold $X$.

$\bullet$ $C^\infty(X,\ce)$ - the space of smooth sections of a
vector bundle $\ce$ over $X$.

$\bullet$ $C^{-\infty}(X,\ce)$ - the space of generalized sections
of a vector bundle $\ce$.

$\bullet$ $WF(u)$ - the wave front of a generalized section $u$.

$\bullet$ $C^{-\infty}_\Lam(X,\ce)$ - the space of generalized
sections of a vector bundle $\ce$ with wave front contained in a
subset $\Lam\subset T^*X$.

$\bullet$ $[[Z]]$ - a current associated to a submanifold $Z$.

$\bullet$ $\cd_k(X)$ - the space of $k$-currents on $X$.

$\bullet$ $D$ - the Rumin differential operator on forms on a
contact manifold.

$\bullet$ $\PP_X$ - the spherical cotangent bundle of $X$, i.e.
$\PP_X:=\PP_+(T^*X)$.

$\bullet$ $N(P)$ - the normal cycle of a subset $P$.

$\bullet$ $\One_P$ - the indicator function of a subset $P\subset
X$.

$\bullet$ $T^*_ZX$ - the conormal bundle of a submanifold $Z\subset
X$.

\section{Background}\label{S:background}

\subsection{Manifolds with corners.}\label{Ss:corners}
\begin{definition} \label{D:corners}
A closed subset $P$ of an $n$-dimensional smooth manifold $X$ is
called a {\itshape submanifold with corners} of dimension $k$ if any
point $p\in P$ has an open neighborhood $U\ni p$ and a
$C^\infty$-diffeomorphism $\phi\colon U\stackrel{\sim}{\to} \RR^n$
such that for some $r\geq 0$
\begin{eqnarray*}
\phi(p)=0,\\
\phi(P\cap U)=\RR^r\times \RR^{k-r}_{\geq 0}\times 0_{\RR^{n-k}}
\end{eqnarray*}
where $0_{\RR^{n-k}}$ is the zero element of $\RR^{n-k}$. The set of
submanifolds with corners is denoted by $\mathcal{P}(X)$.
\end{definition}

The number $r$ is defined uniquely by the point $p$, but may depend
on it. This $r$ is called {\itshape type} of a point $p$.

\begin{example}
\begin{enumerate}
\item Any smooth submanifold of $X$ with or without boundary is a
submanifold with corners.
\item A convex compact $n$-dimensional polytope $P\subset \RR^n$ is a
submanifold with corners if and only if $P$ is simplicial, namely
every vertex has exactly $n$ adjacent edges. The vertices are
precisely the points of type $0$.
\end{enumerate}
\end{example}

A submanifold with corners $P\subset X$ has a natural finite
stratification by locally closed smooth submanifolds as follows. For
any integer $r$, $0\leq r\leq n$ let us denote by $S_r(P)$ the union
of all points of type $r$. Then the $S_r(P)$ are locally closed
smooth disjoint submanifolds of $X$ and
\begin{displaymath}
P=\bigsqcup_{r=0}^{\dim P} S_r(P).
\end{displaymath}

This stratification of $P$ will be called {\itshape canonical
stratification}.

\begin{definition}\label{D:transversality}
Let $P$ and $Q$ be two closed submanifolds with corners of $X$. We
say that $P$ and $Q$ {\itshape intersect transversally} if each
stratum of the canonical stratification of $P$ is transversal to
each stratum of the canonical stratification of $Q$.
\end{definition}

It is well known (se e.g. \cite{part2}, Lemma 2.1.10) that if two
closed submanifold with corners $P$ and $Q$ intersect transversally
in the sense of Definition \ref{D:transversality} then $P\cap Q$ is
also a closed submanifold with corners.

\subsection{Construction of oriented blow up along a
submanifold.}\label{Ss:blow-up} Let $X$ be a smooth manifold and
$Y\subset X$ be its closed submanifold (both $X$ and $Y$ are without
corners or boundary). We remind here the so called oriented blow up
of $X$ along $Y$ which is a manifold with boundary $\tilde X$ and a
smooth map $\alp\colon \tilde X\to X$.

Let $T_YX:=(TX|_Y)/TY$ be the normal bundle of $Y$. Let
$\PP_+(T_YX):=(T_YX\backslash\underline{0})/\RR_{>0}$ be the
oriented projectivization of it. Set theoretically, one defines
$$\tilde X:=\PP_+(T_YX)\bigsqcup (X\backslash Y).$$
The map $\alp\colon \tilde X\to X$ is equal to $Id$ on $X\backslash
Y$, and to the natural projection $\PP_+(T_YX)\to Y\subset X$ on
$\PP_+(T_YX)$.

For the details of defining the smooth structure of $\tilde X$ we
refer to \cite{melrose}, Ch. 5 (see also \cite{alesker-bernig} for a
special situation which suffices too for the purposes of this
article). Here we give an idea how to do that for
$X=\RR^n=\RR^k\times \RR^{n-k}$, $Y=\RR^k\times 0_{n-k}\subset
\RR^n$. Let us consider the imbedding
$$X\backslash Y\to X\times \PP_+(\RR^{n-k})=\RR^k\times
\RR^{n-k}\times \PP_+(\RR^{n-k})$$ given by $(x,y)\mapsto (x,y,l)$
where $x\in \RR^k,\, y\in \RR^{n-k}\backslash\underline{0}$, and $l$
is the only oriented line in $\RR^{n-k}$ passing through $0$ and $y$
and oriented from $0$ to $y$. Then define $\tilde X$ to be the
closure of the image of $X\backslash Y$ under this map. The map
$\alp\colon \tilde X\to X$ is the restriction to $\tilde X$ of the
natural projection $X\times \PP_+(\RR^{n-k})\to X$.

It is not hard to see that the case of a general manifold $X$ can be
deduced from this one by doing this construction in each coordinate
chart and then gluing. Moreover any diffeomorphism $f\colon X\to X$
such that $f(Y)=Y$ has a unique lifting to a diffeomorphism $\tilde
f\colon \tilde X\to \tilde X$ such that the diagram
\begin{eqnarray*}
\square[\tilde X`\tilde X`X`X;\tilde f`\alp`\alp`f]
\end{eqnarray*}
is commutative.

\subsection{The Rumin differential operator.}\label{Ss:Rumin}
Let $M$ be a contact manifold of dimension $2n-1$. Recall that this
means that $M$ is given a smooth distribution of codimension $1$
(i.e. a smooth field of hyperplanes in the tangent bundle) which is
completely non-integrable. More explicitly, locally there exists a
1-form $\alpha$ such that the field of hyperplanes is equal to $Ker
\alpha$ with the property that $\alpha \wedge d\alpha^{n-1} \neq 0$.
This form $\alpha$, which is called contact form, is unique up to
multiplication by a non-vanishing smooth function.

A form $\omega \in \Omega^*(M)$ is called {\it vertical} if it
vanishes on the contact distribution. Given a local contact form
$\alpha$, the form $\omega$ is vertical if and only if $\ome\wedge
\alp=0$, or equivalently $\omega=\alpha \wedge \phi$ for some $\phi
\in \Omega^*(M)$.

Given $\omega \in \Omega^{n-1}(M)$, there exists a {\it unique}
vertical form $\omega' \in \Omega^{n-1}$ such that
$d(\omega+\omega')$ is vertical (see \cite{rum94}). We define the
projection operator $Q:\Omega^{n-1}(M) \to \Omega^{n-1}(M)$ by
setting $Q\omega:=\omega+\omega'$. This operator is a first order
linear differential operator containing vertical forms in its
kernel. The {\it Rumin operator} is the second order linear
differential operator
\begin{displaymath}
D:=d \circ Q:\Omega^{n-1}(M) \to \Omega^n(M).
\end{displaymath}

The Rumin operator is the main ingredient in some differential
complex, called Rumin-de Rham complex, whose cohomology is
isomorphic to the de Rham cohomology \cite{rum94}. We shall not use
this isomorphism in the sequel.

\subsection{Generalized sections of vector bundles and wave
fronts.}\label{Ss:wave fronts} We recall the definition and the main
properties of the wave front of a generalized functions referring
for more details to \cite{hormander-pde1} or
\cite{guillemin-sternberg}, Ch. VI.

Let $X$ be a smooth manifold (always countable at infinity, in
particular paracompact). Let ${\mathcal E}\to X$ be a finite
dimensional vector bundle; for definiteness we assume that
${\mathcal E}$ is a complex bundle though for real bundles the
theory is exactly the same. We denote by $C^\infty(X,{\mathcal E})$
the space of $C^\infty$-sections of ${\mathcal E}$. It is a
Fr\'echet space with topology of uniform convergence on compact
subsets of $X$ of all partial derivatives. We denote by
$C^\infty_c(X,{\mathcal E})$ the space of $C^\infty$-sections of
${\mathcal E}$ with compact support. Naturally
$C^\infty_c(X,{\mathcal E})$ is a locally convex topological vector
space with topology of (strict) countable inductive limit of
Fr\'echet spaces. The inclusion
\begin{displaymath}
C^\infty_c(X,{\mathcal E})\hookrightarrow C^\infty(X,{\mathcal E})
\end{displaymath}
is continuous and has dense image.

Let us denote by $|\omega_X|$ the line bundle over $X$ of complex
densities; thus the fiber of $|\omega_X|$ over a point $x\in X$ is
equal to the one dimensional space of complex valued Lebesgue
measures on the tangent space $T_xX$.

We have a separately continuous bilinear map
\begin{displaymath}
 C^\infty(X,{\mathcal E})\times C^\infty_c(X,{\mathcal E}^*\otimes|\omega_X|)\to \mathbb{C}
\end{displaymath}
given by $(f,g)\mapsto \int_X \langle f,g\rangle$. This map is a
non-degenerate pairing. In other words the induced map
\begin{displaymath}
C^\infty(X,{\mathcal E})\to (C^\infty_c(X,{\mathcal
E}^*\otimes|\omega_X|))^*
\end{displaymath}
is continuous, injective and has a dense image when the target space
is equipped with the weak topology.

\begin{definition}
The space $(C^\infty_c(X,{\mathcal E}^*\otimes|\omega_X|))^*$ is
called the space of {\itshape generalized sections} of ${\mathcal
E}$. It is denoted by $C^{-\infty}(X,{\mathcal E})$.
\end{definition}
Thus $C^\infty(X,{\mathcal E})\subset C^{-\infty}(X,{\mathcal E})$.
Generalized sections of the trivial line bundle are called {\itshape
generalized functions}. Generalized sections of ${\mathcal
E}=|\omega_X|$ are called {\itshape generalized densities}.


The main technical tool in the following will be wave front sets of
a generalized section of a vector bundle. Let us remind this notion
following \cite{hormander-pde1}, Ch. VIII. For simplicity we will
discuss the case of the trivial line bundle; the discussion easily
extends to the general case. Moreover we assume that our manifold
$X$ is $\RR^n$; the considerations in the general case make use of
local coordinates, and independence of the constructions of the
local coordinates is proved in \cite{hormander-pde1}, Ch. VIII. Thus
let $u\in C^{-\infty}_c(\RR^n)$ be a compactly supported generalized
function. Then its Fourier transform $\hat u$ is defined and is a
continuous function on $\RR^n$.

Let us define the subset $\cs(u)\subset\RR^n$ consisting of points
$\eta$ with the property that there exists an open
$\RR_{>0}$-invariant neighborhood $V$ of $\eta$ such that for any
$N\in \NN$ there exists a constant $C_N$ such that
$$|\hat u(\xi)|\leq C_N(1+|\xi|)^{-N} \mbox{ for any } \xi\in V.$$
Denote by $\Sigma(u)$ the complement of $\cs(u)$. Clearly
$\Sigma(u)$ is a closed $\RR_{>0}$-invariant subset of $\RR^n$.

For a point $x\in \RR^n$ define
$$\Sigma_x(u):=\cap_f \Sigma (f\cdot u),\, f\in
C^\infty_c(\RR^n),\, f(x)\ne 0.$$ Finally one defines the wave front
of $u$ by
$$WF(u):=\{(x,\xi)\in\RR^n\times (\RR^n\backslash\{0\})|\, \xi\in
\Sigma_x(u)\}.$$


As we have mentioned, the definition of the wave front generalizes
to any manifolds and to generalized sections of smooth vector
bundles. Then $WF(u)$ is a subset of the cotangent bundle $T^*X$
with the zero section removed. Let us summarize the main properties
of wave fronts relevant for our applications.

\begin{proposition} (\cite{hormander-pde1} or \cite{guillemin-sternberg}, Ch.VI \S
3).\label{P:wave_front_properties} \\
Let $u\in C^{-\infty}(X,{\mathcal E})$.
\begin{itemize}
\item[(i)] The wave front $WF(u)$ is a closed
$\RR_{>0}$-invariant subset of $T^*X\backslash\underline{0}$, where
$\underline{0}$ denotes the zero-section.
\item[(ii)] $WF(u)=\emptyset$ if and only if $u$ is infinitely smooth.
\item[(iii)] \begin{displaymath} WF(u\boxtimes v)\subset\left( WF(u)\times
WF(v)\right)\cup \left( WF(u)\times \underline{0}\right)\cup \left(
\underline{0} \times WF(v)\right).
\end{displaymath}
\end{itemize}
\end{proposition}

Let us fix a closed $\RR_{>0}$-invariant subset $\Lam\subset
T^*X\backslash\underline{0}$. Let us denote by
$C^{-\infty}_\Lam(X,\ce)$ the set of generalized sections of $\ce$
whose wave front is contained in $\Lam$. This is a linear subspace
of $C^{-\infty}(X,\ce)$. Moreover $C^{-\infty}_\Lam(X,\ce)$ is
equipped with a locally convex linear topology. In the case of $\ce$
being the trivial line bundle and $X=\RR^n$ the topology is defined
below; the general case is defined by covering $X$ by a system of
local charts and trivializing $\ce$ in each chart. Thus let us
assume $X=\RR^n$ and $\ce$ is trivial. For any $N\in \NN$, $\phi\in
C^\infty_c(\RR^n)$ and any closed $\RR_{>0}$-invariant subset
$V\subset \RR^n$ such that
$$\Lam\cap (supp (\phi)\times V)=\emptyset$$
define the semi-norm on $C^{-\infty}_\Lam(X,\ce)$
$$||u||_{\phi,V,N}:=\sup_{\xi\in V}|\xi|^N|\widehat{\phi u}(\xi)|.$$
Then one equips $C^{-\infty}_\Lam(X,\ce)$ with the weakest locally
convex topology which is stronger than the weak topology on
$C^{-\infty}(X,\ce)$ and such that all semi-norms
$\{||\cdot||_{\phi,V,N}\}$ are continuous. Notice that if
$\Lam=T^*X\backslash\underline{0}$ then
$C^{-\infty}_\Lam(X,\ce)=C^{-\infty}(X,\ce)$ as linear topological
spaces.


\hfill

Let $f\colon Y\to X$ be a smooth map, ${\mathcal E}\to X$ be a
vector bundle. Then one has the obvious pull-back map on smooth
sections
\begin{displaymath}
f^*\colon C^\infty(X,{\mathcal E})\to C^\infty(Y,f^*{\mathcal E}).
\end{displaymath}
It turns out that $f^*$ can be extended in a natural way to some
generalized sections of ${\mathcal E}$ satisfying appropriate
assumptions. Now we are going to discuss these assumptions referring
for details to \cite{guillemin-sternberg}, Ch. VI \S3.

\def\cme{C^{-\infty}(X,\mathcal{E})}
\def\cmela{C^{-\infty}_\Lambda(X,\mathcal{E})}

Let us fix a closed conic subset $\Lambda \subset
T^*X\backslash\underline{0}$. Let us consider the linear subspace
$C^{-\infty}_\Lambda(X,\mathcal{E})\subset
C^{-\infty}(X,\mathcal{E})$ consisting of generalized sections of
$\mathcal{E}$ with the wave front contained in $\Lambda$. The space
$C^{-\infty}_\Lambda(X,\mathcal{E})$ is equipped with the locally
convex linear topology defined above.

\begin{definition}\label{D:transvesality}
Let $\Lambda\subset T^*X\backslash\underline{0}$ be a closed conic
subset. A smooth map $f\colon Y\to X$ is {\itshape transversal to
$\Lambda$} if for any $\xi\in T^*_{f(y)}X \cap \Lambda$ one has
$df^*_y(\xi)\ne 0$.
\end{definition}

Note that a submersion is transversal to each $\Lambda \subset
T^*X\backslash\underline{0}$.

\begin{proposition}\label{P:lifting}
If a smooth map $f\colon Y\to X$ is transversal to a closed conic
subset $\Lambda\subset T^*X\backslash\underline{0}$ then there is a
unique linear sequentially continuous map, also called pull-back,
\begin{displaymath}
f^*\colon \cmela\to C^{-\infty}_{f^*\Lambda}(Y,f^*\mathcal{E})
\end{displaymath}
where
\begin{displaymath}
f^*(\Lambda):=\{(y,\eta)\in T^*Y|\, \eta\in
df^*_y(\Lambda|_{f(y)})\}
\end{displaymath}
whose restriction to smooth sections of $\mathcal{E}$ is equal to
the pull-back on smooth sections discussed above.
\end{proposition}
Usually in the literature $f^*(\Lam)$ is denoted by $df^*(\Lambda)$.
We will use the shorter notation as above throughout the article.

We denote by $T_A^*M$ the conormal bundle of any smooth submanifold
$A\subset M$.

\begin{remark}\label{restriction-gen-func}
Let $Y$ be a closed submanifold of $X$, and $f\colon Y\to X$ be the
identity imbedding. Then Proposition \ref{P:lifting} says in
particular that $f^*u$ is defined provided
\begin{displaymath}
 WF(u)\cap T^*_YX=\emptyset.
\end{displaymath}
\end{remark}

From Proposition \ref{P:wave_front_properties} and Remark
\ref{restriction-gen-func} one can easily deduce the following
result on product of generalized sections (see
\cite{guillemin-sternberg}, Ch. VI \S 3, Proposition 3.10). Below we
denote by $s: \PP_X \to \PP_X$ the involution
$s((x,[\xi]))=(x,[-\xi])$, and for a subset $Z\subset \PP_X$ we
denote by $Z^s$ the image $s(Z)$.

\begin{proposition} \label{P:tens-prod-wf}
Let ${\mathcal E}_1,{\mathcal E}_2\to X$ be two vector bundles. Let
$\Lambda_1,\Lambda_2$ be closed conic subsets of $T^*X\backslash
\underline{0}$. Let us assume that
\begin{eqnarray}\label{E:tens-prod-wf}
\Lambda_1\cap \Lambda_2^s=\emptyset.
\end{eqnarray}
Let us define a new subset $\Lambda\subset T^*X\backslash
\underline{0}$ such that for any point $x\in X$
\begin{eqnarray}\label{E:wf-of-tens-prod}
\Lambda|_x= (\Lambda_1|_x+\Lambda_2|_x)\cup \Lambda_1|_x\cup
\Lambda_2|_x.
\end{eqnarray}
Then $\Lambda$ is also a closed conic subset, and moreover there
exists a unique bilinear jointly sequentially continuous map
\begin{displaymath}
C^{-\infty}_{\Lambda_1}(X,\mathcal{E}_1)\times
C^{-\infty}_{\Lambda_2}(X,\mathcal{E}_2)\to
C^{-\infty}_{\Lambda}(X,\mathcal{E}_1\otimes \mathcal{E}_2)
\end{displaymath}
whose restriction to smooth sections is the tensor product map.
\end{proposition}

We will need two further technical propositions. We denote by
$\cd(M)$ the space of currents on a manifold $M$ (if $M$ is oriented
$\cd(M)$ coincides with with the space of generalized differential
forms).

\begin{proposition} \label{prop_P:1}
Let $A$ be a smooth submanifold of $M$. Let $\tilde M_A$ denote the
oriented blow up of $M$ along $A$ (see Section \ref{Ss:blow-up}).
Let $f\colon \tilde M_A\to M$ be the natural map. Let $T\in
\mathcal{D}(M)$. Then $f^*T \in \mathcal{D}(\tilde M_A)$ is defined
provided
\begin{displaymath}
 T_A^*M\cap WF(T)=\emptyset.
\end{displaymath}
\end{proposition}

\proof Let $x\in \tilde M_A$. If $x\not\in f^{-1}(A)$ then clearly
$f^*T$ is well defined in a neighborhood of $x$. Thus let us assume
that $x\in f^{-1}(A)$. Let $\xi\in WF(T)|_{f(x)}$. We have to show
that $(df)^*\xi \in T^*_x \tilde M_A$ does not vanish on $T_x(f^{-1}(A))$, or,
equivalently, that the restriction of $\xi$ to $df(T_x f^{-1}(A))$
does not vanish.

But $df(T_xf^{-1}(A))=T_{f(x)}A$, hence
the assumption implies that the
restriction of $\xi$ to $T_{f(x)}A$ does not vanish.
\endproof

\begin{proposition} \label{prop_P:2}
Let $\tilde M_A \stackrel{f}{\longrightarrow} M
\stackrel{g}{\longrightarrow} R$ be smooth maps where $f$ is the
oriented blow up map as in Proposition \ref{prop_P:1}, and $g$ is a
submersion. Then $(gf)^*T$ is well defined provided the following
condition is satisfied:
\newline
for any $a\in A$ and any $\zeta \in WF(T)\cap T^*_{g(a)}R$ the
restriction of $\zeta$ to $dg(T_aA)$ does not vanish.
\end{proposition}
\proof Since $g$ is a submersion, $g^*T$ is defined, and by
Proposition \ref{P:lifting}
\begin{eqnarray}\label{E:3}
WF(g^*T)|_a\subset (dg_a)^*(WF(T)|_{g(a)}).
\end{eqnarray}
By Proposition \ref{prop_P:1} $f^*(g^*T)$ is defined provided
\begin{displaymath}
T^*_AM\cap WF(g^*T)=\emptyset.
\end{displaymath}

By \eqref{E:3} this condition is satisfied provided that for any
$a\in A$
\begin{displaymath}
(T^*_AM)|_a\cap (dg_a)^*(WF(T)|_{g(a)})=\emptyset.
\end{displaymath}

The last condition is equivalent to the assumption of the
proposition.
\endproof

\def\dens{|\omega_X|}


Now let us discuss the push-forward of generalized sections. Let
$f\colon Y\to X$ be a smooth {\itshape proper} map. Let
$\mathcal{E}\to X$ be a smooth vector bundle as above. One can
define the push-forward map
\begin{displaymath}
 f_*\colon C^{-\infty}(Y,f^*(\mathcal{E}\otimes \dens^*)\otimes
|\omega_Y|)\to C^{-\infty}(X,\mathcal{E})
\end{displaymath}
as the dual map to
\begin{displaymath}
 f^*\colon C^\infty_c(X,\mathcal{E}^*\otimes\dens)\to
C^\infty_c(Y,f^*(\mathcal{E}^*\otimes\dens)).
\end{displaymath}

Note that $f^*$ indeed takes compactly supported sections to
compactly supported ones due to the properness of $f$.

\begin{remark}\label{R:push-forward-dens}
Let us take $\mathcal{E}=\dens$. Then we get the push-forward map on
generalized densities:
\begin{displaymath}
 f_*\colon C^{-\infty}(Y,|\omega_Y|)\to C^{-\infty}(X,\dens).
\end{displaymath}
In the case when $f$ is a proper submersion, $f_*$ is integration
along the fibers.
\end{remark}

For a closed conic subset $W\subset T^*Y\backslash\underline{0}$ let
us define a new conic subset
\begin{displaymath}
f_*W:=\{(x,\eta)\in T^*X \setminus \underline{0} \,|\,\exists y\in
f^{-1}(x) \mbox{ such that } (y,df_y^*(\eta))\in W|_x\cup\{0\}\}.
\end{displaymath}
The standard notation for $f_*W$ in the literature is $df_*W$, but
we prefer to abbreviate it as above. Since $f$ is proper, $f_*W$ is
closed. One has the following result (see
\cite{guillemin-sternberg}, Ch. VI \S 3, Proposition 3.9).
\begin{proposition}\label{P:wf-push-forward}
Let $W\subset T^*Y\backslash\underline{0}$ be a closed conic subset.
Then
\begin{displaymath}
f_*\colon C^{-\infty}_W(Y,f^*(\mathcal{E}\otimes \dens^*)\otimes
|\omega_Y|)\to C_{f_*W}^{-\infty}(X,\mathcal{E})
\end{displaymath}
is a sequentially continuous linear map.
\end{proposition}

\hfill

In the geometric measure theory, generalized sections of the space
of forms $\Ome^{n-k}(X)$, $\dim X=n$, are called $k$-currents. An
oriented compact $k$-submanifold (possibly, with corners) $Z\subset
X$ defines a current $[[Z]]\in C^{-\infty}(X,\Ome^{n-k})$ given by
$$\phi\mapsto \int_Z\phi.$$ The operator on currents dual to the de
Rham differential $d$ is denoted by $\pt$. The space of $k$-currents
is also denoted by $\cd_k(X)$.

By Proposition \ref{P:tens-prod-wf} we have a partial product
$$\cap\colon C^{-\infty}(X,\Ome^{n-k})\times
C^{-\infty}(X,\Ome^{n-l})\dashrightarrow
C^{-\infty}(X,\Ome^{n-k-l})$$ extending the usual wedge product on
smooth differential forms. For this reason we will denote sometimes
the $\cap$-product by $\wedge$. If $Z_1,Z_2\subset X$ are oriented
compact submanifolds (with corners) and $Z_1,Z_2$ intersect
transversally then
$$[[Z_1]]\cap[[Z_2]]=[[Z_1\cap Z_2]].$$
Also one has partially defined operations of pull-back and
push-forward on currents as in Propositions \ref{P:lifting} and
\ref{P:wf-push-forward}.


\subsection{Normal cycles.}\label{Ss:normal-cycles} Let us recall
the definition of normal cycle of a compact submanifold with
corners. This is a generalization of a unit conormal bundle of a
smooth submanifold.

Let $X$ be an $n$-dimensional smooth manifold. For the simplicity of
the notation we assume that $X$ is oriented though this assumption
can be easily removed. Let $TX$ denote the tangent bundle of $X$,
and $T^*X$ the cotangent bundle. Let $\PP_+(T^*X)=:\PP_X$ denote the
oriented projectivization of the cotangent bundle $T^*X$, namely
$\PP_X\to X$ is a bundle over $X$ such that its fiber over a point
$x\in X$ is equal to the quotient space
$(T_x^*X\backslash\{0\})/\RR_{>0}$. Notice that if one fixes a
Riemannian metric on $X$ then $\PP_X$ is naturally identified with
the spherical cotangent bundle of $X$. For this reason we call
$\PP_X$ the spherical cotangent bundle even if no metric is fixed.

Let $P\in \cp(X)$ be a compact submanifold with corners. For any
point $x\in P$ let $T_xP\subset T_xX$ denote the tangent cone of $P$
at the point $x$. It is defined as follows:
$$T_xP:=\{\xi\in T_xX| \mbox{ there exists a } C^1-\mbox{map }
\gamma\colon [0,1]\to P \mbox{ such that }\gamma (0)=x \mbox { and }
\gamma'(0)=\xi\}.$$ One can easily see that $T_xP$ is a convex
polyhedral cone. If $P$ has no corners or boundary then $T_xP$ is
the usual tangent space of $P$ at $x$. Let $(T_xP)^\circ \subset
T^*_xX$ denote the dual cone. Recall that the dual cone $C^o$ of a
convex cone $C$ in a linear space $W$ is defined by
$$C^o:=\{y\in W^*|\, y(x)\geq 0 \mbox{ for any } x\in C\}.$$
By definition the normal cycle of $P$ is a subset $N(P)\subset
\PP_X$ defined by
$$N(P):=\cup_{x\in P}((T_xP)^o\backslash\{0\})/\RR_{>0}.$$
It is well known (and easy to see) that $N(P)$ is
$(n-1)$-dimensional submanifold (with singularities in general); the
orientation of $X$ induces an orientation of $N(P)$; it is a cycle,
i.e. $\pt N(P)=0$; and it is Legendrian with respect to the
canonical contact structure on $\PP_X$.

\begin{remark}
J. Fu \cite{fu-94} has defined normal cycle of a subanalytic subset
of a real analytic manifold. An essentially equivalent notion of
characteristic cycle of such a set was introduced independently
using different tools by Kashiwara (see \cite{kashiwara-schapira}).
In this article we are not going to use this more complicated
situation.
\end{remark}

\subsection{Valuations on manifolds.}\label{Ss:valuations-mfld} The general reference for
 valuations on manifolds is the series of articles \cite{part1},
\cite{part2},\cite{part3},\cite{part4}, as well as the survey
\cite{alesker-survey}.

For simplicity of the notation we assume that $X$ is an oriented
manifold, though this assumption can be easily removed, and all
results and constructions can be generalized. A choice of
orientation on $X$ induces in a standard way a choice of orientation
on the normal cycle $N(P)$ for any $P\in \cp(X)$ (we will not
explain it here).

\begin{definition}\label{D:smooth-val}
A smooth valuation is a map $\phi\colon \cp(X)\to \CC$ such that
there exist a smooth $n$-form $\mu$ on $X$ (remind $n=\dim X$), and
a smooth  $(n-1)$-form $\omega$ on $\PP_X$ such that
$$\phi(P)=\int_P\mu+\int_{N(P)}\omega$$
for any $P\in \cp(X)$.
\end{definition}
\begin{remark}\label{remark-val-basic}
1) Any $\phi$ of such form is a finitely additive functional on
$\cp(X)$.

2) In its present form the definition of smooth valuation might look
unmotivated. The original definition of a smooth valuation from
\cite{part2} was different but equivalent to this one, and behind it
stay some non-trivial characterization results from the theory of
valuations on convex sets and facts from geometric measure theory
(see \cite{part2}, \cite{part3}).

3) Various pairs of forms $(\mu,\ome)$ may define the same
valuations. The pairs $(\mu,\ome)$ defining the zero valuation have
been described by an explicit system of integral and differential
equations by Bernig and Br\"ocker \cite{bernig-broecker} as follows:
$(\mu,\omega)$ defines the zero valuation if and only if
\begin{eqnarray*}
\pi_*\ome=0,\mbox{ and }\\
D\ome+\pi^*\mu=0
\end{eqnarray*}
where $\pi\colon \PP_X\to X$ is the natural projection, $D$ is the
Rumin operator (see Section \ref{Ss:Rumin}), $\pi_*$ is the operator
of integration of a differential form along the fibers.
\end{remark}

\begin{example}
1) Any smooth measure on $X$ is a smooth valuation. Indeed take a
pair $(\mu,0)$.

2) The Euler characteristic $\chi$ is a smooth valuations. This is
less trivial to see. This follows for example from an old theorem of
Chern \cite{chern-45} who constructed explicitly a pair $(\mu,\ome)$
defining the Euler characteristic; his construction depends on a
choice of an auxiliary Riemannian metric on $X$. A different proof
was given in \cite{part2}.
\end{example}

Thus the space $V^\infty(X)$ of smooth valuations is a quotient
space of $\Ome^n(X)\oplus \Ome^{n-1}(\PP_X)$, and hence it carries
the natural linear Fr\'echet topology. Moreover $V^\infty(X)$ has a
canonical filtration by closed subspaces
$$V^\infty(X)=W_0\supset W_1\supset W_2\supset \dots\supset W_n$$
which can be described as follows. The last term $W_n$ is equal to
the subspace of smooth measures. Each other subspace $W_i$, $0\leq
i<n$, consist of smooth valuations which can be given by a pair
$(\mu,\ome)\in \Ome^n(X)\oplus \Ome^{n-1}(\PP_X)$ where $\ome$
satisfies the following condition: for any point $z\in \PP_X$ and
any $(n-1)$-dimensional subspace $L\subset T_z\PP_X$ such that the
intersection of $L$ with the tangent space at $z$ to the fiber
$\pi^{-1}(\pi(z))$ is at least $(n-i)$-dimensional, the restriction
of $\ome$ to $L$ vanishes. In particular $W_0=V^\infty(X)$.

Next we have the evaluation on points map $V^\infty(X)\to
C^\infty(X)$, i.e. $\phi\mapsto [x\mapsto \phi(\{x\})]$. This linear
map is onto, its kernel is equal to $W_1$. Thus the Fr\'echet space
of smooth functions is canonically isomorphic to the Fr\'echet space
$V^\infty(X)/W_1$. The following result was proved in \cite{part3},
\cite{part4}.
\begin{theorem}
There exists a canonical product $V^\infty(X)\times V^\infty(X)\to
V^\infty(X)$ which is

(1) continuous;

(2) commutative and associative;

(3) the filtration $W_\bullet$ is compatible with it:
$$W_i\cdot W_j\subset W_{i+j}$$
where we define $W_{j} =0$ for $j>n$;

(4) the Euler characteristic $\chi$ is the unit in the algebra
$V^\infty(X)$;

(5) this product commutes with restrictions to locally closed
submanifolds; in particular the evaluation on points map
$V^\infty(X)\to C^\infty(X)$ is a homomorphism of algebras.
\end{theorem}

Let us describe a construction of the product on valuations
following \cite{alesker-bernig}. This is different from the
construction given in \cite{part3}, and it will be used in this
article. First let us consider the natural map
$$\pi_{X\times X}\colon \PP_{X\times X}\to X\times X.$$
Let us denote by $\PP$ the preimage under this map of the diagonal
in $X\times X$. Elements of $\PP$ can be written in homogeneous
coordinates $(x,[\xi:\eta])$ where $x\in X$, $\xi,\eta\in T^*_xX$.
Let $\cm$ denote the smooth submanifold of $\PP$
$$\cm:=\{(x,[\xi:0])\}\cup \{(x,[0:\eta])\}\cup \{(x,[\xi:-\xi])\}$$
where clearly the union is disjoint. Let $\bar\PP$ denote the
oriented blow up of $\PP$ along $\cm$ (see Section
\ref{Ss:blow-up}). We have the maps $\Phi\colon \PP\backslash \cm\to
\PP_X\times_X \PP_X$ given by
$\Phi((x,[\xi:\eta]))=(x,[\xi],[\eta])$. Also we have a smooth map
$p\colon \PP\backslash \cm\to \PP_X$ given by
$p(x,[\xi:\eta])=(x,[\xi+\eta])$. It is not hard to see that there
exist unique smooth maps
\begin{eqnarray}\label{E:bar-phi}
\bar\Phi\colon \bar\PP\to \PP_X\times_X\PP_X,\\
\label{E:bar-p} \bar p\colon \bar\PP\to \PP_X
\end{eqnarray}
extending $\Phi$ and $p$ respectively. Moreover $\bar \Phi$ and
$\bar p$ are proper. Let
\begin{eqnarray}\label{E:q-maps}
q_1,q_2\colon \PP_X\times_X\PP_X\to \PP_X
\end{eqnarray}
be the obvious projections onto the first and second factors
respectively. Let $s\colon \PP_X\to \PP_X$ be the involution given
by $s(x,[\xi])=(x,[-\xi])$.

The following result was proved in \cite{alesker-bernig}.
\begin{theorem} \label{mthm_prod}
Let $X$ be an $n$-dimensional oriented manifold; let $\phi_i \in
V^\infty(X)$ be represented by $(\omega_i,\mu_i) \in
\Omega^{n-1}(\PP_X) \times \Omega^n(X)$, $i=1,2$. Then the product
$\phi_1 \cdot \phi_2$ is represented by
\begin{align*}
\omega & = \bar p_*\bar\Phi^*(q_1^* \omega_1 \wedge q_2^* D\omega_2) + \omega_1 \wedge \pi^*\pi_* \omega_2 \in \Omega^{n-1}(\PP_X), \nonumber\\
\mu & = \pi_* (\omega_1 \wedge s^* (D\omega_2+\pi^* \mu_2))+ \mu_1
\wedge \pi_* \omega_2 \in \Omega^n(X).
\end{align*}
\end{theorem}

An important property of the product on valuations is a version of
the Poincar\'e duality. To describe it let us denote first by
$V_c^\infty(X)$ the subspace of $V^\infty(X)$ of compactly supported
valuations. Then $V_c^\infty(X)$ admits a continuous integration
functional
\begin{eqnarray}\label{intfunc}
\int\colon V^\infty_c(X)\to \CC
\end{eqnarray}
given by $\phi\mapsto \phi(X)$. Consider the bilinear map
$$V^\infty(X)\times V^\infty_c(X)\to \CC$$
defined by $(\phi, \psi)\mapsto \int\phi\cdot \psi$. The following
theorem was proved in \cite{part4}, and a simpler proof was given in
\cite{bernig-quat}.
\begin{theorem}\label{selfduality}
This bilinear form is a perfect pairing. In other words, the induced
map
$$V^\infty(X)\to (V^\infty_c(X))^*$$
is injective and has a dense image (with respect to the weak
topology on $(V^\infty_c(X))^*$).
\end{theorem}
Define $V^{-\infty}(X):=(V^\infty_c(X))^*$. Elements of this space
are called {\itshape generalized valuations}. Thus by the last
theorem $V^\infty(X)$ is identified with a dense subspace of
$V^{-\infty}(X)$: $$\Xi_\infty\colon V^\infty(X)\inj
V^{-\infty}(X).$$

Let us denote by $\Theta\colon \Ome^n_c(X)\oplus
\Ome^{n-1}_c(\PP_X)\to V_c^\infty(X)$ (where the subscript $c$ stays
for the compact support) the map given by
$$(\Theta(\mu,\ome))(P)=\int_P\mu+\int_{N(P)}\ome.$$
This map is onto. The dual map
$$\Theta^*\colon V^{-\infty}(X)\to (\Ome^n_c(X))^*\oplus
(\Ome^{n-1}_c(\PP_X))^*=C^{-\infty}(X)\oplus
C^{-\infty}(\PP_X,\Omega^{n})$$ induces an injective closed
imbedding. The result of Bernig and Br\"ocker (see Remark
\ref{remark-val-basic} (3)) implies that if a generalized valuation
$\psi$ corresponds under this imbedding to a pair of currents
$(C,T)$ then
\begin{eqnarray}
T\mbox{ is a cycle, i.e. }\pt T=0;\label{current1}\\
T\mbox{ is Legendrian;}\label{current2}\\
\pi_*T=\pt C\label{current3}.
\end{eqnarray}
Conversely any pair of currents $(C,T)$ satisfying
(\ref{current1})-(\ref{current3}) corresponds to a unique
generalized valuation. (Recall that $T$ is called Legendrian if for
a contact form $\alp$ one has $T\wedge\alp=0$.) The described above
product on smooth valuations has been extended in
\cite{alesker-bernig} to a partially defined product on generalized
valuations. More precisely the following result was proved in
\cite{alesker-bernig}. (Below for a subset $\Lam$ in the cotangent
bundle of a manifold $M$ we denote by $\Lam^s$ the image of $\Lam$
under the antipodal involution $s$ acting on cotangent vectors to
$M$.)
\begin{theorem} \label{thm_prod_gen}
Let $\Lambda_1,\Lambda_2\subset T^*X\backslash\underline{0},\,
\Gamma_1,\Gamma_2\subset T^*\PP_X\backslash\underline{0}$ be closed
conic subsets such that $d\pi_*(\Gamma_i)\subset \Lambda_i,\,
i=1,2$. Let us assume that these subsets satisfy the following
conditions:

\begin{enumerate}
\item[(a)] $\Lambda_1\cap \Lambda_2^s=\emptyset$. \label{condition_down}
\item[(b)] $\Gamma_1\cap (\pi^*(\Lambda_2))^s=\emptyset$.
\item[(c)] $\Gamma_2\cap (\pi^*(\Lambda_1))^s=\emptyset$.
\item[(d)] If $(x,[\xi_i],u_i,0) \in
\Gamma_i$ for $i=1,2$, then $u_1 \ne -u_2$.
\item[(e)] If $(x,[\xi])\in \PP_X$ and
\begin{align*}
(u,\eta_1) & \in \Gamma_1|_{(x,[\xi])}\\
(-u,\eta_2) & \in \Gamma_2|_{(x,[-\xi])}
\end{align*}
then
\begin{displaymath}
d\theta^*(0,\eta_1,\eta_2) \neq (0,l,-l)\in
T^*_{(x,[\xi],[\xi])}(\PP_X \times_X \PP_X),
\end{displaymath}
where $\theta\colon\PP_X \times_X \PP_X \to \PP_X \times_X \PP_X$ is
defined by $\theta(x,[\xi_1],[\xi_2])=(x,[\xi_1],[-\xi_2])$.
\end{enumerate}
Then there is a unique jointly sequentially continuous bilinear map, called a
partial product,
\begin{displaymath}
V^{-\infty}_{\Lambda_1,\Gamma_1}(X)\times
V^{-\infty}_{\Lambda_2,\Gamma_2}(X)\to V^{-\infty}(X)
\end{displaymath}
such that for any $\psi_1,\psi_2 \in V^\infty(X)$ the product
$\Xi_\infty(\psi_1) \cdot \Xi_\infty(\psi_2)$ exists (i.e. the
conditions (a)-(e) are satisfied), and
\begin{equation} \label{eq_extension_smooth}
\Xi_\infty(\psi_1) \cdot \Xi_\infty(\psi_2)=\Xi_\infty(\psi_1 \cdot
\psi_2).
\end{equation}
\end{theorem}

In fact this partial product of generalized valuations can be
described as follows. Consider again the maps
$$\PP_X\overset{\bar p}{\leftarrow} \bar \PP\overset{\bar\Phi}{\to} \PP_X\times_X \PP_X$$
as in (\ref{E:bar-phi})-(\ref{E:bar-p}). Assume that generalized
valuations $\psi_i\in V^\infty_{\Lam_i,\Gamma_i}$, $i=1,2$, are
represented by pairs of currents $(C_i,T_i)$. Then the product
$\psi_1\cdot \psi_2$ is represented by the pair $(C,T)$ where
\begin{eqnarray}\label{E:may0001}
T:=(-1)^n\bar p_*\bar\Phi^*(q_1^*T_1\cap q_2^*T_2)+\pi^*C_1\cap T_2
+T_1\cap \pi^*C_2\in C^{-\infty}(\PP_X,\Ome^n),\\\label{E:may0002}
C:=C_1\cap C_2\in C^{-\infty}(X).
\end{eqnarray}

\hfill

The filtration $\{W_\bullet\}$ on smooth valuations $V^\infty(X)$
discussed above extends naturally to $V^{-\infty}(X)$ as follows
(see \cite{part4}, Section 7.3, for the details). Let us denote by
$W_i(V^{-\infty}(X))$ the closure of $W_i$ in the weak topology on
$V^{-\infty}(X)$. Then by \cite{part4}, Corollary 7.3.3, one has
$$W_i(V^{-\infty}(X))\cap V^\infty(X)=W_i.$$
By the abuse of notation we will denote $W_i(V^{-\infty}(X))$ by
$W_i$.

\hfill

Let us explain now the explicit relation of product on valuations to
constructible functions assuming for simplicity of exposition that
$X$ is a real analytic manifold. As we have explained in Section
\ref{Ss:detailed-overview} of the introduction we have a canonical
imbedding with dense image of the space of constructible functions
$$\Xi_{\cp}\colon \cf(X)\inj V^{-\infty}(X).$$
Let us assume that compact submanifolds with corners $P_1,P_2$ are
transversal to each other. Then it was shown in
\cite{alesker-bernig} that the product of generalized valuations
$\Xi_{\cp}(\One_{P_1})$ and $\Xi_{\cp}(\One_{P_2})$ is well defined,
i.e. the assumptions of Theorem \ref{thm_prod_gen} are satisfied,
and
\begin{eqnarray}\label{E:prod-constr}
\Xi_{\cp}(\One_{P_1})\cdot
\Xi_{\cp}(\One_{P_2})=\Xi_{\cp}(\One_{P_1\cap P_2}).
\end{eqnarray}

The filtration $\{W_\bullet(V^{-\infty}(X))\}$ induces on $\cf(X)$
the filtration by codimension of support (see \cite{part4},
Proposition 8.2.2).

\hfill

The following technical lemma was proved in \cite{alesker-bernig},
Lemma 8.2, it will be useful later.
\begin{lemma}\label{C:approxim}
Let $\xi\in V^{-\infty}(\RR^n)$ be a generalized valuation given by
a pair of currents $(C,T)\in \cd_n(\RR^n)\times
\cd_{n-1}(\PP_{\RR^n})$. Let
$$\Lambda\subset T^*(\RR^n)\backslash\underline{0},\,
\Gamma\subset T^*(\PP_{\RR^n})\backslash\underline{0}$$ be closed
conic subsets such that $\pi_*(\Gamma)\subset \Lambda$ and
$$WF(C)\subset \Lambda,\, WF(T)\subset \Gamma.$$
Then there exists a sequence of smooth valuations $\{\xi_j\}\subset
V^\infty(\RR^n)$ corresponding to currents $(C_j,T_j)$ such that
\begin{eqnarray*}
C_j\to C \mbox{ in } C^{-\infty}_{\Lambda}(\RR^n),\\
T_j\to T \mbox{ in } C^{-\infty}_{\Gamma}(\PP_{\RR^n},\Ome^n).
\end{eqnarray*}
In particular $\xi_j\to\xi$ in $V^{-\infty}(X)$.
\end{lemma}

\section{Exterior product of generalized
valuations.}\label{S:exter-prod} \setcounter{subsection}{1} In this
section we define the exterior product of two generalized
valuations. The construction generalizes the construction of
exterior product of two smooth valuations from \cite{alesker-bernig}
though we do not prove this fact here. Let $X_1,X_2$ be smooth
manifolds (without boundary) of dimensions $n_1,n_2$ respectively.
As usual we assume for simplicity of the notation that they are
oriented. Set $X:=X_1\times X_2$. Let $\cm_1,\cm_2\subset\PP_{X}$ be
the submanifolds defined by
\begin{eqnarray*}
\cm_1=\{(x_1,x_2;[\xi_1:0])\},\\
\cm_2=\{(x_1,x_2;[0:\xi_2])\}.
\end{eqnarray*}
Let $F\colon \hat\PP\to\PP_{X}$ be the oriented blow-up map along
$\cm_1\bigsqcup\cm_2$. Let
$$\Phi\colon \hat\PP\to \PP_{X_1}\times \PP_{X_2}$$
be the only smooth map extending the map
$\PP\backslash(\cm_1\bigsqcup \cm_2)\to \PP_{X_1}\times \PP_{X_2}$
given by
$$(x_1,x_2;[\xi_1:\xi_2])\mapsto
\left((x_1,[\xi_1]);(x_2,[\xi_2])\right).$$

Let us also consider the diagrams
\begin{eqnarray*}
\PP_{X_1}\overset{p_1}{\leftarrow}\PP_{X_1}\times
X_2\overset{i_1}{\inj}\PP_X,\\
\PP_{X_2}\overset{p_2}{\leftarrow}X_1\times\PP_{X_2}\overset{i_2}{\inj}\PP_X
\end{eqnarray*}
where $p_1,p_2$ are the obvious projections, and $i_1$ is given by
$i_1((x_1,[\xi_1]),x_2)=(x_1,x_2;[\xi_1:0])$, and $i_2$ is defined
similarly.

Let us denote by $\tilde p_i\colon X_1\times X_2\to X_i$, $i=1,2$,
the obvious projections. For $i=1,2$ let us define the spaces
$$\cl_i:=\{(C_i,T_i)\in \cd_{n_i}(X_i)\oplus
\cd_{n_i-1}(\PP_{X_i})|\, T_i \mbox{ is Legendrian }, \pt T_i=0,
\pi_{X_{i}*}T_i=\pt C_i\}.$$

Let $\psi_i\in V^{-\infty}(X_i),i=1,2$. Then $\psi_i$ is represented
by by a unique pair $(C_i,T_i)\in \cl_i$. In fact this gives a topological isomorphism
between $V^{-\infty}(X_i)$ and $\cl_i$ (see Section 8 of \cite{alesker-bernig}).
 Let us consider
$$(C,T)\in \cd_{n_1+n_2}(X)\oplus
\cd_{n_1+n_2-1}(\PP_{X})$$ defined by
\begin{eqnarray}\label{D:ext-pr1}
C:=C_1\boxtimes C_2,\\\label{D:ext-pr2} T:=F_*\Phi^*(T_1\boxtimes
T_2)+(\tilde p_1\circ \pi_X)^*C_1\cdot (i_{2*}p_2^*T_2) +
(i_{1*}p_1^*T_1)\cdot (\tilde p_2\circ \pi_X)^*C_2
\end{eqnarray}
where $\pi_X\colon \PP_X\to X$ is the natural map.

\begin{claim}\label{Cl:well-defined}
The currents $C,T$ are well defined. Moreover the map
$\left((C_1,T_1), (C_2,T_2)\right)\mapsto (C,T)$ defined by the
formulas (\ref{D:ext-pr1}), (\ref{D:ext-pr2}) defines a bilinear
jointly sequentially continuous map
$$\left(\cd_{n_1}(X_1)\oplus \cd_{n_1-1}(\PP_{X_1})\right)\times
\left(\cd_{n_2}(X_2)\oplus \cd_{n_2-1}(\PP_{X_2})\right)\to
\cd_{n_1+n_2}(X)\oplus \cd_{n_1+n_2-1}(\PP_{X}).$$
\end{claim}
{\bf Proof.} The claim is obviously satisfied for $C$. To prove both
parts of the  claim for $T$, observe first of all that they are
satisfied for the first summand in the definition of  $T$ (i.e. in
the right hand side of (\ref{D:ext-pr2})) since $\Phi$ is a
submersion. The second and the third summands are symmetric, so let
us consider only the third summand.

We have
\[WF(i_{1*}p_1^*T_1)\subset
i_{1*}p_1^*(T^*\PP_{X_1}\backslash\underline{0}).\] In order to
describe the latter set, we will assume that $X_1,X_2$ are vector
spaces. Also let us fix a point $A\in \PP_X$ which we write in the
form
$$A=(x_1,x_2;[\xi_1:\xi_2]).$$
The fiber of $T^*\PP_X$ over $A$ we can write as
\begin{eqnarray*}
\{(x_1^*,x_2^*;\zeta)|\, x_1^*\in T^*_{x_1}X_1,\, x_2^*\in
T^*_{x_2}X_2,\, \zeta\in T^*_{[\xi_1:\xi_2]}\PP_+(X^*)\}.
\end{eqnarray*}
The fiber of $i_{1*}p_1^*(T^*\PP_{X_1}\backslash\underline{0})$ over
$A$ is non-empty if and only if $\xi_2=0_{X_2^*}$, and in this case
this fiber is equal to
\begin{eqnarray}\label{wf-ip}
\left( X_1^*\oplus \{0_{X_2^*}\}\oplus
T^*_{[\xi_1:0_{X_2^*}]}\PP_+(X^*)\right)\backslash\{0\}.
\end{eqnarray}
On the other hand the fiber over $A$ of $WF(( \tilde
p_2\circ \pi_X)^*C_2)$ is contained in the fiber of $(\tilde
p_2\circ \pi_X)^*WF(C_2)$ which is contained obviously in
\begin{eqnarray}\label{wf-phi2}
\{0_{X_1^*}\}\oplus X_2^*\oplus \{0\}
\end{eqnarray}
where the last zero belongs to $T^*_{[\xi_1:0_{X_2^*}]}\PP_+(X^*)$.
But clearly the sets (\ref{wf-ip}) and (\ref{wf-phi2}) satisfy the
condition of disjointness. The claim follows. \qed

\begin{claim}
Assume that $(C_i,T_i)\in \cl_i, i=1,2$. Let $C$, $T$ be defined by
the formulas (\ref{D:ext-pr1}), (\ref{D:ext-pr2}) respectively. Then

(i) $T$ is Legendrian;

(ii) $\pt T=0$.

(iii) $\pi_*T=\pt C$.
\end{claim}
{\bf Proof.} Let us make first few observations.

(a) For compact submanifolds with corners
$P_i\subset X_i,\, i=1,2,$ let us consider the corresponding pairs of currents $(C_i,T_i):=([[P_i]],[[N(P_i)]])$.
For such pairs one has
\begin{eqnarray*}
C=[[P_1\times P_2]],\\
T=[[N(P_1\times P_2)]].
\end{eqnarray*}
Clearly they satisfy the conditions (i)-(iii) of the claim, namely belongs to the subspace
$\cl$.

(b) Now in general the map
$$((C_1,T_1),(C_2,T_2))\mapsto (C,T)$$
is jointly sequentially continuous.

(c) The subspace $\cl\subset  \cd_{n_1+n_2}(X)\oplus\cd_{n_1+n_2-1}(\PP_X)$
is closed in the weak topology.

(d) There is a natural isomorphism of topological vector spaces $V^{-\infty}(X_i)\simeq \cl_i,\, V^{-\infty}(X)\simeq \cl$
(see Section 8 of \cite{alesker-bernig}).

These observations and Corollary 1.3 in \cite{alesker-sequential-closure} imply that any jointly (in fact, even separately) sequentially continuous map
$\cl_1\times \cl_2\to \cd_{n_1+n_2}(X)\oplus\cd_{n_1+n_2-1}(\PP_X)$ which maps pairs of indicator functions
of compact submanifolds with corners to $\cl$, has image contained in $\cl$.
\qed




\begin{claim}
The exterior product of constructible functions considered as
generalized valuations coincides with the usual product of
constructible functions.
\end{claim}
{\bf Proof.} This immediately follows from the construction. \qed







\section{Pull-back and push-forward of valuations.}\label{S:pull-push}

\subsection{Pull-back of smooth valuations under
immersions.}\label{Ss:pull-immersion} Let
$$f\colon X^n\to Y^m$$ be an immersion of smooth manifolds, $n<m$.
Let us define the pull-back map on smooth valuations
$$f^*\colon V^\infty(Y)\to V^\infty(X)$$
as follows. By the sheaf property of smooth valuations (see
\cite{part2}, Theorem 2.4.10) it is enough to define pull-back
locally on $X$. Then we may assume that $f$ is a closed imbedding.
Then define
\begin{eqnarray}\label{def-pull-immer}
(f^*\phi)(P):=\phi(f(P))
\end{eqnarray}
for any $P\in \cp(X)$. The following claim is obvious.
\begin{claim}
If $f\colon X\to Y$ is a smooth immersion then
$$f^*\colon V^\infty(Y)\to V^\infty(X)$$
is a continuous linear operator.
\end{claim}

Actually we will need yet another description of the pull-back under
immersions. Let us assume that $X\subset Y$ is a closed submanifold.
Let us consider the obvious maps
\begin{eqnarray}\label{def-q-theta}
X\overset{q}{\leftarrow} \PP_+(T^*_XY)\overset{\theta}{\inj}\PP_Y.
\end{eqnarray}

\def\xy{X\times_Y\PP_Y}
\def\txy{\widetilde{X\times_Y\PP_Y}}

Observe first that $\xy$ is a smooth submanifold of $\PP_Y$.
Moreover the oriented projectivization $\PP_+(T_X^*Y)$ of the
conormal bundle of $X$ in $Y$ is a closed submanifold in $\xy$. Let
us denote by $\txy$ the oriented blow up of $\xy$ along
$\PP_+(T_X^*Y)$. Let $\alp\colon \txy\to \PP_Y$ be the composition
of the blow up map $\txy\to\xy$ with the imbedding $\xy\to \PP_Y$.

Let   $$(df)^*\colon X\times_YT^*Y\to T^*X$$ be the map adjoint to
the differential $df$. It is easy to see that it induces a smooth
map to the oriented projectivization of $T^*X$
$$\beta\colon \txy\to \PP_X.$$
Thus consider the diagram
\begin{eqnarray}\label{def-alp-beta}
\PP_X\overset{\beta}{\leftarrow} \txy\overset{\alp}{\to} \PP_Y.
\end{eqnarray}

\begin{proposition}\label{P:smooth-pull-imbed}
Let $\xi\in V^\infty(Y)$ be a smooth valuation given by
$\xi(P)=\int_P\mu+\int_{N(P)}\ome$ where $\mu\in \Ome^m(Y),\,
\ome\in \Omega^{m-1}(Y)$. Then $f^*\xi$ is given by
$(f^*\xi)(P)=\int_P\mu'+\int_{N(P)}\ome'$ where
$$\mu'=q_*\theta^*\ome,\, \ome'=\beta_*\alp^*\ome.$$
\end{proposition}
{\bf Proof.} It follows immediately from the definition of $f^*$.
\qed

One can easily deduce from Proposition \ref{P:smooth-pull-imbed} the
following result that pull-back preserves the filtration
$\{W_\bullet\}$ on valuations discussed in Section
\ref{Ss:valuations-mfld}.
\begin{proposition}\label{P:pull-filtr-immer}
If $f\colon X\to Y$ is a smooth immersion then for any $i$ one has
on smooth valuations
$$f^*(W_i)\subset W_i.$$
\end{proposition}

\subsection{Push-forward of smooth valuations under
submersions.}\label{Ss:push-submersion} Let
$$f\colon X^n\to Y^m$$
be a smooth proper submersion. We define an operator of push-forward
$$f_*\colon V^\infty(X)\to V^\infty(Y)$$
by the formula
\begin{eqnarray}\label{def-pull-back-submer}
(f_*\phi)(P)=\phi(f^{-1}(P))
\end{eqnarray}
for any $P\in \cp(Y)$. The following claim follows from Proposition
\ref{P:push-smooth-formula} below.
\begin{claim}
If $f$ is a proper submersion then $f_*\colon V^\infty(X)\to
V^\infty(Y)$ is a continuous linear operator.
\end{claim}
\begin{remark}\label{Rem:nonproper-submersion}
If $f\colon X\to Y$ is a submersion which is not necessarily proper,
then similarly the formula (\ref{def-pull-back-submer}) defines a
continuous linear operator
$$f_*\colon V_c^\infty(X)\to V_c^\infty(Y)$$
when both spaces are equipped with topologies of inductive limits of
Fr\'echet spaces. This also follows from Proposition
\ref{P:push-smooth-formula} below.
\end{remark}

Let us describe the operation of push-forward on the language of
differential forms. Let $\phi\in V^{\infty}_c(X)$. By \cite{part4},
Lemma 2.1.1, there exist $(\nu,\omega)\in \Ome^n_c(X)\times
\Omega^{n-1}_c(\PP_X)$ such that $\phi(P)=\int_P\nu
+\int_{N(P)}\ome$ for any $P\in \cp(X)$. (In the case when $f$ is
proper, it is not necessary to consider valuations with compact
support.)

Let us consider the following diagram
$$\PP_X\overset{(df)^*}{\hookleftarrow}X\times_Y\PP_Y\overset{p}{\to}\PP_Y$$
where $p$ is the natural projection, and $(df)^*$ is the adjoint of
the differential $df$. Clearly $(df)^*$ is a closed imbedding since
$f$ is a submersion.
\begin{proposition}\label{P:push-smooth-formula}
The valuation $f_*\phi$ is represented by the pair
$(f_*\nu,p_*(df^*)^*\omega)$.
\end{proposition}
{\bf Proof}. It is straightforward and left to a reader. \qed

One can easily see from the definitions that $f_*$ respects the
filtration in the following way.
\begin{proposition}\label{P:push-filtr-smooth}
Let $f\colon X\to Y$ is a proper submersion. Then for any $i$ one
has on smooth valuations
$$f_*(W_i)\subset W_{i-\dim X+\dim Y}$$
where one defines $W_j=W_0$ for $j<0$.
\end{proposition}

\subsection{Pull-back of generalized valuations under
submersions.}\label{Ss:submersions} Let
$$f\colon X^n\to Y^m$$
be a submersion of smooth manifolds. Let us define a pull-back
operator
$$f^*\colon V^{-\infty}(Y)\to V^{-\infty}(X)$$
as the adjoint operator to $f_*\colon V_c^\infty(X)\to
V_c^\infty(Y)$ defined in Section \ref{Ss:push-submersion}.

\begin{remark}\label{R:pull-non-smooth}
Under submersions, pull-back of smooth valuations may be not smooth.
For example let us consider the projection $f\colon \RR^2\to \RR$
onto the first coordinate. Let $\mu$ be a Lebesgue measure on $\RR$
(which is a smooth valuation of course). Then $f^*\mu$ is not
smooth. In fact $f^*\mu$ can be given the following description
which is not completely rigorous at this point. For any convex
compact set $K\subset \RR^2$ one has
$$(f^*\mu)(K)=\mu(f(K))$$
is the Lebesgue measure of the projection of $K$.
\end{remark}
\begin{example}
The example from the previous remark can be generalized as follows.
Again, we do not make here the statements completely rigorous in
order to clarify the intuition. Let $f\colon \RR^n\to \RR^k$ be an
orthogonal projection. Let $\phi$ be a smooth valuation on $\RR^k$,
for example volume $vol_k$ or intrinsic volume $V_i$. Then for any
convex compact set $K\subset \RR^n$ one has
$$(f^*\phi)(K)=\phi(f(K)).$$
Moreover if $P\in \cp(\RR^n)$ is a compact submanifold with corners
which is "in generic position" to the projection $f$, then
$$(f^*vol_k)(P)=\int_{y\in \RR^k}\chi(P\cap f^{-1}(y))\,dvol_k(y).$$
\end{example}

Let us describe the pull-back on valuations in terms of currents. As
in Section \ref{Ss:push-submersion}, let us consider the following
diagram
$$\PP_X\overset{df^*}{\hookleftarrow}X\times_Y\PP_Y\overset{p}{\to}\PP_Y$$
where $p$ is the natural projection, and $df^*$ is the adjoint of
the differential $df$.  We have the following result.
\begin{proposition}\label{P:pull-submer-gener-expr}
Let $f\colon X^n\to Y^m$ be a smooth submersion. Let $\xi\in
V^{-\infty}(Y)$ corresponds to a pair of currents $(C,T)\in
\cd_m(Y)\times \cd_{m-1}(\PP_Y)$. Then $f^*\xi\in V^{-\infty}(X)$
corresponds to the pair of currents $(C',T')\in \cd_n(X)\times
\cd_{n-1}(\PP_X)$ where
\begin{eqnarray}
C'=f^*C,\\
T'=(df^*)_*p^*T.
\end{eqnarray}
\end{proposition}
{\bf Proof.} Let $\phi\in V_c^\infty(X)$ be an arbitrary smooth
valuation with compact support. By \cite{part4}, Lemma 2.1.1, there
exist $(\nu,\omega)\in \Ome^n_c(X)\times \Omega^{n-1}_c(\PP_X)$ such
that $\phi(P)=\int_P\nu +\int_{N(P)}\ome$ for any $P\in \cp(X)$. By
Proposition \ref{P:push-smooth-formula}
\begin{eqnarray*}
<\xi,f_*\phi>=<C,f_*\nu>+<T,p_*(df^*)^*\ome>=\\
<f^*C,\nu>+<(df^*)_*p^*T,\ome>.
\end{eqnarray*}
The proposition is proved. \qed

\begin{proposition}\label{P:pull-submer-gener-constr}
Let $f\colon X\to Y$ be a submersion. Let $P\in\cp(Y)$. Then
$f^*(\Xi_\cp(\One_P))=\Xi_\cp(\One_{f^{-1}(P)})$. In other words the
pull-back of $\One_P$ considered as a generalized valuation
coincides with the usual pull-back on (constructible) functions.
\end{proposition}
{\bf Proof.} This follows from Proposition
\ref{P:pull-submer-gener-expr}. \qed

By dualizing Proposition \ref{P:push-filtr-smooth} one easily sees
that pull-back is compatible with filtrations, namely we have the
following result.
\begin{proposition}\label{P:pull-back-filtr-gen}
Let $f\colon X\to Y$ be a submersion. Then on generalized valuations
one has
$$f^*(W_i(V^{-\infty}(Y)))\subset W_i (V^{-\infty}(X)).$$
\end{proposition}

\subsection{Push-forward of generalized valuations under
immersions.}\label{Ss:push-gener-immer} Let
$$f\colon X\to Y$$ be a smooth immersion. Let us define
the push-forward operator
$$f_*\colon V_c^{-\infty}(X)\to V_c^{-\infty}(Y)$$
as the adjoint operator to pull-back $f^*\colon V^\infty(Y)\to
V^\infty(X)$ defined in Section \ref{Ss:pull-immersion}. Similarly
for a proper immersion $f$ we define $f_*\colon V^{-\infty}(X)\to
V^{-\infty}(Y)$ as dual to $f^*\colon V^\infty_c(Y)\to
V^\infty_c(X)$.
\begin{remark}\label{R:p-f-gener}
(1) The push-forward under immersion of a smooth valuation is
usually not smooth. For example consider the imbedding $f\colon
\{0\}\inj \RR$. The $f_*(1)$ is the delta-measure at 0.

(2) The space of generalized densities (which is the completion of
the space of smooth measures in the weak topology) is a subspace of
generalized valuations by \cite{part4}, Proposition 7.3.5. Then the
push-forward under immersions (and, in fact, under more general
maps) on this subspace coincides with the usual push-forward of
generalized densities.

(3) If $f\colon X\to Y$ is a closed imbedding and $P\in\cp(X)$, then
$f_*(\Xi_\cp(\One_P))=\Xi_\cp(\One_{P})$. Thus push-forward under a
closed imbedding of constructible functions is just extension by
zero.
\end{remark}

\hfill

Let us describe now the push-forward $f_*$ using the language of
currents. We will assume for simplicity that $f\colon X\to Y$ is a
closed imbedding, and $\dim X<\dim Y$ (of course, if $\dim X=\dim Y$
the description is evident). By our definition, $f_*\colon
\vmi(X)\to \vmi(Y)$ is dual to $f^*\colon V^\infty_c(Y)\to
V^\infty_c(X)$.

Consider the diagram as in (\ref{def-alp-beta})
\begin{eqnarray}\label{dia-1}
\PP_X\overset{\beta}{\leftarrow} \txy\overset{\alp}{\to} \PP_Y
\end{eqnarray}
where $\alp=\delta\circ \eps$ with
\begin{eqnarray}\label{dia-2}
\txy\overset{\eps}{\to}X\times_Y\PP_Y\overset{\delta}{\to}\PP_Y
\end{eqnarray}
with $\eps$ being the oriented blow up map, $\delta$ being the
natural imbedding, and $\beta$ is induced by the dual of the
differential $df^*\colon X\times_YT^*Y\to T^*X$.

\begin{proposition}\label{P:push-forward-gener-imbed-currents}
Assume that $\phi\in \vmi(X)$ is given by a pair of currents
$(C,T)$. Then $f_*\phi$ is given by the pair $(0,\alp_*\beta^*T)$.
\end{proposition}
{\bf Proof.} Notice first of all that the current $\alp_*\beta^*T$
is well defined since $\beta$ is a submersion. The operator
$(C,T)\mapsto \alp_*\beta^*T$ is a continuous linear operator
$$ \vmi(X)\to \cd_{\dim Y-1}(\PP_Y).$$ Hence it is enough to prove
the proposition for a dense subspace of generalized valuations. Thus
it is sufficient to prove it for $\phi =\Xi_\cp(\One_P)$ for any
$P\in \cp(X)$. By Remark \ref{R:p-f-gener}(3)
$f_*(\Xi_\cp(\One_P))=\Xi_\cp(\One_P)$. But the normal cycle of $P$
as a subset of $Y$ is equal to $\alp_*\beta^*$ of the normal cycle
of $P$ considered as a subset of $X$. \qed

Dualizing Proposition \ref{P:pull-filtr-immer} one easily obtains
the following proposition.
\begin{proposition}
Let $f\colon X\to Y$ be a closed imbedding. Then for any $i$ one has
on generalized valuations
$$f_*(W_i(V^{-\infty}(X)))\subset W_{i-\dim X+\dim Y}(V^{-\infty}(Y)).$$
\end{proposition}

\subsection{Pull-back of generalized valuations under
immersions.}\label{Ss:pull-gener-immer} This case is somewhat more
involved than the previous ones. Let
$$f\colon X^n\to Y^m$$ be an immersion. We want to construct a
{\itshape partially defined} map $f^*\colon
V^{-\infty}(Y)\dashrightarrow V^{-\infty}(X)$ which coincides on
smooth valuations with the map from Section \ref{Ss:pull-immersion}.

Generalized valuations form a sheaf by \cite{part4}, Proposition
7.2.2. Thus it is sufficient to define $f^*$ locally in a way
compatible with restrictions to open subsets. Hence we may assume
that $f\colon X\to Y$ is a closed imbedding.


Consider the diagram as in (\ref{def-alp-beta})
\begin{eqnarray}\label{dia-1}
\PP_X\overset{\beta}{\leftarrow} \txy\overset{\alp}{\to} \PP_Y
\end{eqnarray}
where $\alp=\delta\circ \eps$ with
\begin{eqnarray}\label{dia-2}
\txy\overset{\eps}{\to}X\times_Y\PP_Y\overset{\delta}{\to}\PP_Y
\end{eqnarray}
with $\eps$ being the oriented blow up map, $\delta$ being the
natural imbedding.


Let $\xi\in V^{-\infty}(Y)$ be a generalized valuation given by a
pair of currents $(C,T)\in \cd_m(Y)\times \cd_{m-1}(\PP_Y)$. Then we
define $f^*\xi$ to correspond to the pair
\begin{eqnarray}\label{D:pull-gener-immer}
(C',T'):=(f^*C, \beta_*\alp^*T)
\end{eqnarray}
whenever this makes sense. Now we are going to formulate the precise
sufficient conditions that the above pair is well defined, and to show
that $T'$ is Legendrian, $\pt T'=0$, and $\pi_{X*} T'=\pt C'$, and
prove various compatibilities.

\hfill


\def\vlg{V^{-\infty}_{\Lam,\Gamma}}

Let us introduce more notation. Let
\begin{eqnarray*}
\Lambda\subset T^*Y\backslash\underline{0},\\
\Gamma\subset T^*\PP_Y\backslash\underline{0}
\end{eqnarray*}
be closed conic subsets. Let us denote by
$V^{-\infty}_{\Lam,\Gamma}(Y)$ the linear subspace of generalized
valuations on $Y$ corresponding to currents $(C,T)$ with
$WF(C)\subset \Lam,\, WF(T)\subset \Gamma$. Thus
$$\vlg(Y)\subset C^{-\infty}_\Lam(Y)\times
C^{-\infty}_\Gamma(\PP_Y,\Ome^{\dim Y})$$ is a closed subspace. Let
us equip $\vlg(Y)$ with the induced topology from the latter space
(see Section \ref{Ss:wave fronts}).


\begin{definition}\label{D:def-sets-transversal}
(1) We say that an imbedding $f\colon X\to Y$ is transversal to
$(\Lam,\Gamma)$ if the following conditions of transversality are
satisfied:

($TR_1$) $\Lam\cap T^*_XY=\emptyset$;

($TR_2$) $\Gamma\cap T^*_{X\times_Y\PP_Y}(\PP_Y)=\emptyset$;

($TR_3$) $\Gamma\cap T^*_{\PP_+(T^*_XY)}(\PP_Y)=\emptyset$.

(2) We say that an immersion $f\colon X\to Y$ is transversal to
$(\Lam,\Gamma)$ if for every point $x\in X$ there exist an open
neighborhood $U$ of $x$ and an open neighborhood $V$ of $f(x)$ such
that $f|_U\colon U\to V$ is a closed imbedding, and $f|_U$ is
transversal to $(\Lam,\Gamma)$ in the sense of part (1) of the
definition.
\end{definition}

\begin{definition}\label{D:map-val-transversal}
Let $f\colon X\to Y$ be an immersion. Let $\xi\in V^{-\infty}(Y)$ be
a generalized valuation given by a pair $(C,T)$. We say that $f$ is
{\itshape transversal} to $\xi$ if $f$ is transversal to
$(WF(C),WF(T))$ in the sense of Definition
\ref{D:def-sets-transversal}(2).
\end{definition}

\begin{remark}\label{R:transversality-map-val}
(1) It is easy to see that if $f$ is a smooth imbedding then parts
(1) and (2) of Definition \ref{D:def-sets-transversal} are
equivalent.

(2) If $\xi$ is a smooth valuation then any immersion $f\colon X\to
Y$ is transversal to $\xi$ for trivial reasons.
\end{remark}

\begin{claim}
Let $f\colon X\to Y$ be an immersion transversal to a generalized
valuation $\xi\in V^{-\infty}(Y)$. Then the currents
(\ref{D:pull-gener-immer}) are well defined, and $f^*\colon
V^{-\infty}_{\Lam,\Gamma}(Y)\to V^{-\infty}(X)$ is a sequentially continuous
linear map where $(\Lam,\Gamma):=(WF(C),WF(T))$.
\end{claim}
{\bf Proof.} This easily follows from Propositions \ref{P:lifting},
\ref{P:wf-push-forward}. \qed

\begin{proposition}\label{P:pull-immer-smooth-equiv}
If $f\colon X\to Y$ is an immersion, and $\xi\in V^\infty(Y)$ is a
smooth valuation then the pull-back of $\xi$ in the sense of
generalized valuations and in the sense of smooth valuations
(Section \ref{Ss:pull-immersion}) coincide.
\end{proposition}
{\bf Proof.} Recall that we have diagrams of maps
\begin{eqnarray*}
X\overset{q}{\leftarrow} \PP_+(T^*_XY)\overset{\theta}{\inj}\PP_Y,\\
\PP_X\overset{\beta}{\leftarrow} \txy\overset{\alp}{\to} \PP_Y.
\end{eqnarray*}

Let $\xi\in V^\infty(Y)$ be given by
$$\xi(P)=\int_P\mu +\int_{N(P)}\ome.$$
Then by \cite{bernig-quat} (see also Theorem \ref{mthm_prod} of this
article) $\xi$ corresponds to a pair of currents $(C,T)$ where
$$C=\pi_{Y*}\ome,\, T=s_Y^*(D_Y\ome +\pi_Y^*\mu),$$
and $D_Y$ is the Rumin operator on $\PP_Y$. Then the pull-back
$f^*\xi$ in the sense of generalized valuations is given by the pair
of currents $(C',T')$ where
$$C'=f^*(\pi_{Y*}\ome),\,
T'=\beta_*\alp^*s_Y^*(D_Y\ome+\pi_Y^*\mu)$$ where $\alp,\beta$ are
as in (\ref{def-alp-beta}).

On the other hand let us compute $f^*\xi$ in the smooth sense. By
Proposition \ref{P:smooth-pull-imbed} $f^*\xi$ is given by
$(f^*\xi)(P)=\int_P\mu'+\int_{N(P)}\ome'$ where
$$\mu'=q_*\theta^*\ome,\, \ome'=\beta_*\alp^*\ome.$$
Thus we have to show that
\begin{eqnarray*}
C'=\pi_{X*}(\beta_*\alp^*\ome),\\
T'=s_X^*(D_X(\beta_*\alp^*\ome)+\pi^*_{X}q_*\theta^*\ome).
\end{eqnarray*}
More explicitly
\begin{eqnarray*}
f^*(\pi_{Y*}\ome)=\pi_{X*}(\beta_*\alp^*\ome),\\
\beta_*\alp^*s_Y^*(D_Y\ome+\pi_Y^*\mu)=s_X^*(D_X(\beta_*\alp^*\ome)+\pi^*_{X}q_*\theta^*\ome).
\end{eqnarray*}
The first equality is straightforward. To prove the second one let us
observe first of all that
$$\alp^*s_Y^*\pi_Y^*\mu=\alp^*\pi_Y^*\mu=0.$$
Thus it remains to show that for any $\ome\in \Ome^{m-1}(Y)$ one has
\begin{eqnarray*}
\beta_*\alp^*s_Y^*D_Y\ome=s_X^*(D_X(\beta_*\alp^*\ome)+\pi^*_{X}q_*\theta^*\ome).
\end{eqnarray*}
Since $\alp$ commutes with $s_Y$, and $\beta$ intertwines $s_Y$ and
$s_X$, the last equality is equivalent to
\begin{eqnarray}\label{E:important-equality}
\beta_*\alp^*D_Y\ome=D_X(\beta_*\alp^*\ome)+\pi^*_{X}q_*\theta^*\ome.
\end{eqnarray}

We will need the following well known lemma which follows easily
from the Stokes formula.
\begin{lemma}\label{L:stokes}
Let $f\colon A\to B$ be a smooth proper fibration where $B$ is a
manifold without boundary, and $A$ is a manifold possibly with
boundary. Assume moreover that the restriction
$$f|_{\pt A}\colon \pt A\to B$$
is a fibration too. Then $$f_*(d\eta)=d(f_*\eta)+(-1)^{\dim
A+\deg\eta+1}(f|_{\pt A})_*\eta.$$
\end{lemma}
Let us apply this lemma to the form $\alp^*\ome$ and the map
$$\beta\colon \txy\to \PP_X.$$
We get
\begin{eqnarray}\label{E:push-push}
\beta_*\alp^*d\ome=d\beta_*\alp^*\ome+(-1)^{n+1}\left(\beta|_{\pt
(\txy)}\right)_*\alp^*\ome.
\end{eqnarray}

We will need two more lemmas whose proofs are postponed till the end
of the proof of Proposition \ref{P:pull-immer-smooth-equiv}.
\begin{lemma}\label{L:help1}
Let $\ome\in \Ome^\bullet(\PP_Y)$. If $d\ome$ is vertical then
$d\beta_*\alp^*\ome$ is also vertical.
\end{lemma}
\begin{lemma}\label{L:help2}
For any form $\ome\in \Ome^{m-1}(\PP_Y)$ one has
$$(-1)^{n+1}\left(\beta|_{\pt
(\txy)}\right)_*\alp^*\ome=\pi_X^*q_*\theta^*\ome.$$
\end{lemma}

Let us finish the proof of Proposition
\ref{P:pull-immer-smooth-equiv}. The equality (\ref{E:push-push})
and Lemma \ref{L:help2} imply the equality
\begin{eqnarray}\label{E:useful-eq}
\beta_*\alp^*d\ome=d\beta_*\alp^*\ome+\pi_X^*q_*\theta^*\ome.
\end{eqnarray}
Observe first of all that if $\gamma\in \Ome^\bullet(\PP_Y)$ is a
vertical form then $\theta^*\gamma=0$ since $\PP_+(T^*_XY)=N(X)$ is
a Legendrian cycle in $\PP_Y$. Hence, by the construction of the
Rumin operator, one may assume in equality (\ref{E:useful-eq}) that
$d\ome$ is vertical, namely $d\ome =D_Y\ome$. But then by Lemma
\ref{L:help1} we have $d\beta_*\alp^*\ome=D_X\beta_*\alp^*\ome$.
Thus we eventually obtain
\begin{eqnarray*}
\beta_*\alp^*D_Y\ome=D_X\beta_*\alp^*\ome+\pi_X^*q_*\theta^*\ome.
\end{eqnarray*}
This is exactly (\ref{E:important-equality}). Proposition
\ref{P:pull-immer-smooth-equiv} is proved. \qed

\hfill

{\bf Proof of Lemma \ref{L:help1}.} Let $\lambda_X$ be a contact
form on $X$, and $\lambda_Y$ be a contact form on $Y$. Then we claim
that
\begin{eqnarray}\label{E:contact-forms}
\beta^*(\lam_X)=h\cdot\alp^*(\lam_Y)
\end{eqnarray} where $h$ is a smooth non-vanishing function on $\txy$. Let us prove it
in local coordinates. Let us assume that $X=\RR^n$ is imbedded into
$Y=\RR^m$ using the first $n$ coordinates. Let us denote the
coordinates in $Y$ by $y_1,\dots,y_m$, and in the dual space by
$z_1,\dots,z_m$. Furthermore we will identify $\PP_Y$ with the unit
sphere bundle in $\RR^m\times S^{m-1}\subset\RR^m\times \RR^{m*}$,
and $\PP_X$ with the unit sphere bundle $\RR^n\times S^{n-1}$.
Contact form on a contact manifold is defined up to a product by a
non-vanishing function. Then we may choose
\begin{eqnarray*}
\lam_Y=\sum_{i=1}^{m}z_idy_i,\, \lam_X=\sum_{i=1}^{n}z_idy_i.
\end{eqnarray*}
It is clear that the restriction of $\lam_Y$ to
$X\times_YT^*Y=\RR^n\times \RR^{m*}$ is equal to $\sum_{i=1}^n
z_idy_i$. But obviously this is equal to the pull-back of $\lam_X$
under the natural map $X\times_YT^*Y\to T^*X$. This proves
(\ref{E:contact-forms}).

Next by the assumption of the lemma we have $d\ome\wedge\lam_Y=0$.
We want to show that
$$d(\beta_*\alp^*\ome)\wedge\lam_X=0.$$
We have for any test form $\phi\in \Ome^\bullet_c(\PP_X)$
\begin{eqnarray*}
\int_{\PP_X}\phi\wedge\lam_X\wedge
d(\beta_*\alp^*\ome)\overset{(\ref{E:push-push})}{=}
\int_{\PP_X}\phi\wedge\lam_X\wedge\left(\beta_*\alp^*d\ome+(-1)^{n}\left(\beta|_{\pt
(\txy)}\right)_*\alp^*\ome\right)=\\
\int_{\txy}\beta^*(\phi\wedge\lam_X)\wedge\alp^*(d\ome)+
(-1)^{n}\int_{\pt(\txy)}\left(\beta|_{\txy}\right)^*(\phi\wedge\lam_X)\wedge\alp^*\ome\overset{(\ref{E:contact-forms})}{=}\\
\int_{\txy}h\cdot \beta^*\phi\wedge\alp^*(\lam_Y\wedge
d\ome)+(-1)^{n}\int_{\pt(\txy)}h\cdot
\beta^*\phi\wedge\alp^*(\lam_Y\wedge\ome)=\\
(-1)^{n}\int_{\pt(\txy)}h\cdot
\beta^*\phi\wedge\alp^*(\lam_Y\wedge\ome)
\end{eqnarray*}
where the last equality is due to the fact that the expression under
the first integral vanishes by the assumption $\lam_Y\wedge
d\ome=0$. But the expression under the last integral also vanishes
since the restriction of $\alp^*\lam_Y$ to $\pt(\txy)$ vanishes:
indeed $\pt(\txy)$ is the primage of $\PP_+(T^*_XY)=N(X)$ under the
blow up map $\txy\to \xy$, and the restriction of $\lam_Y$ to $N(X)$
vanishes since $N(X)$ is Legendrian. Hence Lemma \ref{L:help1} is
proved. \qed

\hfill

{\bf Proof of Lemma \ref{L:help2}.} Consider the following
commutative diagram
\begin{eqnarray*}
\square<-1`1`1`-1;700`400>[\PP_X`\pt(\txy)`X`N(X);\beta`\pi_X`\alp_1`q]
\end{eqnarray*}
where $\alp_1\colon \pt(\txy)\to N(X)$ is the restriction of the
blow up map. It is easy to see that this is a Cartesian square.
Denote $\tau:=\theta^*\ome=\ome|_{N(X)}$. Thus we have to show that
$$(-1)^{n+1}\left(\beta|_{\pt(X\times_Y\PP_Y)}\right)_*\alp_1^*\tau=\pi_X^*q_*\tau.$$
Taking into account the orientations, the last formula holds always
for all Cartesian squares of proper fibrations (the base change
formula). \qed

\begin{proposition}
Let $f\colon X\to Y$ be an immersion transversal to $\xi\in
V^{-\infty}(Y)$. Let $\xi$ be given by a pair of currents $(C,T)\in
\cd_m(Y)\times \cd_{m-1}(\PP_Y)$. Define $(C',T')\in \cd_n(X)\times
\cd_{n-1}(\PP_X)$ by the formula (\ref{D:pull-gener-immer}). Then
$T'$ is Legendrian, $\pt T'=0$, and $\pi_{X*}(T')=\pt C'$.
\end{proposition}
{\bf Proof.} Let us denote
$$\Lam:=WF(C),\, \Gamma:=WF(T).$$
The formula (\ref{D:pull-gener-immer}) defines a sequentially continuous linear
operator
$$f^*\colon C^{-\infty}_{\Lam}(Y)\oplus
C^{-\infty}_{\Gamma}(\PP_Y,\Ome^m)\to C^{-\infty}(X)\oplus
C^{-\infty}(\PP_X,\Ome^n)$$ which we also denote by $f^*$. By Lemma
\ref{C:approxim} locally on $Y$ there exists a sequence of smooth
valuations $\{\xi_j\}\subset V^\infty(Y)$ such that
\begin{eqnarray*}
\xi_j\to \xi \mbox{ in }V^{-\infty}_{\Lam,\Gamma}(Y).
\end{eqnarray*}
But then from the properties of wave fronts we have
$$f^*\xi_j\to f^*\xi\mbox{ in }V^{-\infty}(X).$$
Let $\xi_j$ correspond to a pair of smooth currents $(C_j,T_j)$. By
Proposition \ref{P:pull-immer-smooth-equiv} the conclusion of our
proposition is satisfied for smooth currents $f^*((C_j,T_j))$
corresponding to $f^*\xi_j$. Hence by sequential continuity the proposition is
satisfied also by $(C',T')$. \qed

\hfill


\begin{remark}\label{R:map-wave-fronts}
It is important to notice that if $\Lam\subset
T^*Y\backslash\underline{0},\, \Gamma\subset
T^*(\PP_Y)\backslash\underline{0}$ are closed conic subsets, and if
an imbedding $f\colon X\to Y$ is transversal to $(\Lam,\Gamma)$ then
the pull-back
$$f^*\colon V^{-\infty}_{\Lam,\Gamma}(Y)\to
V^{-\infty}_{f^*\Lam,\beta_*\alp^*\Gamma}(X)$$ is a sequentially continuous
linear operator. This is obvious from the properties of wave fronts.
\end{remark}

\begin{proposition}\label{P:imbed-compos}
Let $g\colon X\to Y,\, f\colon Y\to Z$ be imbeddings. Let $\xi\in
\vlg(Z)$. Assume that $f$ and $f\circ g$ are transversal to
$(\Lam,\Gamma)$, and $g$ is transversal to
$(f^*\Lam,\beta_*\alp^*\Gamma)$. Then the maps
$$(f\circ g)^*,g^*\circ f^*\colon V^{-\infty}_{\Lam,\Gamma}(Z)\to
V^{-\infty}(X)$$ coincide.
\end{proposition}
{\bf Proof.} Since the statement is local, we may assume that $X,Y,Z$ are vector spaces.
The assumptions imply that the two linear operators
\begin{eqnarray*}
(f\circ g)^*,g^*\circ f^*\colon \vlg(Z)\to V^{-\infty}(X)
\end{eqnarray*}
are sequentially continuous. By Lemma \ref{C:approxim} smooth valuations are
sequentially dense in the source space. Hence it is enough to show that these two
operators coincide on smooth valuations. But this is obvious. \qed


\begin{proposition}\label{P:pull-back-constuct}
Let $f\colon X\to Y$ be a closed imbedding. Let $P\in \cp(Y)$ be a
compact submanifold with corners such that $X$ is transversal to any
strata of the canonical subdivision of $P$. Then $f$ is transversal
to the generalized valuation $\Xi_\cp(\One_P)$ and
$$f^*(\Xi_\cp(\One_P))=\Xi_\cp(\One_{P\cap X}).$$
\end{proposition}
{\bf Proof.} Let us denote $\Lam:=WF(\One_P),\,
\Gamma:=WF([[N(P)]])$. First we have to check the conditions of
transversality ($TR_1$)-($TR_3$)  from Definition
\ref{D:def-sets-transversal}.

The condition ($TR_1$) easily follows from the transversality of $X$
to every stratum of $P$ and the observation that $\Lam$ is contained
in the union of conormal bundles to all strata of $P$.

The condition ($TR_3$) immediately follows from the observation that
\begin{eqnarray}\label{E:febr0}
N(P)\cap \PP_+(T^*_XY)=\emptyset
\end{eqnarray} which again is a
simple consequence of transversality of $X$ to all the strata of
$P$.

Let us check $(TR_2$). This condition is equivalent to
\begin{eqnarray}\label{E:febr1}
WF([[N(P)]])\cap \pi_Y^*(T^*_XY)=\emptyset
\end{eqnarray}
where $\pi_Y\colon \PP_Y\to Y$ is the natural projection. Let us
denote $P_0:=P\cap X$ which is a submanifold with corners of $X$
(again due to the transversality of $X$ and strata of $P$). Let us
fix a point $x_0\in P_0$. Again as a consequence of transversality
of $X$ and strata of $P$, in a neighborhood of $x_0$ we can write
$Y=X\times Z$, $P=P_0\times Z$ where $X$ and $Z$ may and will be
assumed vector spaces (for the convenience of the notation),
$P_0\subset X$. Then one has
\begin{eqnarray}\label{E:febr2}
N(P)=\{(x,z;[\xi:0_{Z^*}])|\, (x,[\xi])\in N(P_0),\, z\in Z\}\subset
\PP_Y.
\end{eqnarray}
This description easily implies (\ref{E:febr1}). Thus we have shown
that $f^*(\Xi_{\cp}(\One_P))$ is well defined.

In order to prove the second part of the proposition, we notice that
the equality $P=P_0\times Z$ implies the equality
\begin{eqnarray}\label{E:febr3}
f^*(\One_P)=\One_{P\cap X}
\end{eqnarray}
(where $\One_P\in C^{-\infty}(Y)$, $\One_{P\cap X}\in
C^{-\infty}(X)$). Thus it remains to show that
$$\beta_*\alp^*([[N(P)]])=[[N(P\cap X)]]$$
where as previously we have
$$\PP_X\overset{\beta}{\leftarrow}\widetilde{X\times_Y\PP_Y}\overset{\alp}{\to}\PP_Y.$$
Recall that the map $\alp$ is equal to the composition
$$\widetilde{X\times_Y\PP_Y}\overset{\eps}{\to}
X\times_Y\PP_Y\overset{\delta}{\to}\PP_Y$$ where $\delta$ is the
natural imbedding, and $\eps$ is the oriented blow up map along
$\PP_+(T^*_XY)$.

From (\ref{E:febr2}) it is easy to see that
\begin{eqnarray}\label{E:febr4}
\delta^*([[N(P)]])=[[N(P)\cap
(X\times_Y\PP_Y)]]=[[\left\{(x;[\xi:0_{Z^*}])|\, (x,[\xi])\in
N(P_0)\right\}]].
\end{eqnarray}

Then (\ref{E:febr4}) and (\ref{E:febr0}) imply readily that
\begin{eqnarray}\label{E:febr5}
\beta_*\eps^*(\delta^*([[N(P)]])=[[N(P_0)]].
\end{eqnarray}
The equalities (\ref{E:febr3}) and (\ref{E:febr5}) imply the
proposition. \qed

\hfill

The following result will be needed later for integral geometric
applications. It is another special case of the heuristic fact that
$(f\circ g)^*=g^* \circ f^*$.
\begin{proposition}\label{P:compos-immer-submer}
Let $X\overset{g}{\to}Y\overset{f}{\to}Z$ be smooth maps such that
$g$ is an immersion, $f$ is a submersion, and the composition
$f\circ g\colon X\to Z$ is an immersion. Let $\xi\in V^{-\infty}(Z)$
be a generalized valuation. Assume that $f\circ g$ is transversal to
$\xi$. Then $g$ is transversal to $f^*\xi$ and
$$g^*(f^*(\xi))=(f\circ g)^*\xi.$$
\end{proposition}
{\bf Proof.} The proof is a straightforward tedious computation
which we are going to do. First by localization we may assume that
$g$ and $f\circ g$ are imbeddings. Localizing further we may assume
that $X$ is a vector space, and for some vector spaces $S$ and $L$
\begin{eqnarray*}
Z=X\times S,\, Y=X\times S\times L,\\
g\colon X\to X\times S\times L \mbox{ is given by } g(x)=(x,0_S,0_L),\\
f\colon X\times S\times L\to X\times S \mbox{ is the projection.}
\end{eqnarray*}
Thus $(f\circ g)(x)=(x,0_S)$. Let $\xi$ be given by a pair of
currents $(C,T)$. Let us denote
\begin{eqnarray*}
\Lambda:=WF(C)\subset T^*(X\times S)\backslash \underline{0},\\
\Gamma:=WF(T)\subset T^*(\PP_{X\times S})\backslash\underline{0}.
\end{eqnarray*}

Let us recall that we have the maps
$$\PP_Y\overset{df^*}{\hookleftarrow}\PP_Z\times L\overset{p}{\to} \PP_Z.$$
Recall that $f^*\xi$ is given by the pair
$(C',T'):=(f^*C,(df^*)_*p^*T)$.

\underline{Step 1.} Let us check that $g$ and $f^*\xi$ satisfy the
first condition of transversality $(TR_1)$. More precisely let us
check that $g^*(f^*C)$ is well defined, namely
$$WF(f^*C)\cap T^*_{X}(X\times S\times L)=\emptyset.$$

Identifying $T^*(X\times S\times L)=(X\times S\times L)\times
(X^*\times S^*\times L^*)$, let us  observe that $$T^*_X(X\times
S\times L)=\{(x,0_S,0_L;0_{X^*},s^*,l^*)\}.$$ Next we have
$$WF(f^*C)\subset f^*\Lambda=\{(x,s,l;x^*,s^*,0_{L^*})|\,
(x,s;x^*,s^*)\in \Lambda\}.$$ Thus
\begin{eqnarray*}
f^*\Lambda\cap T^*_X(X\times S\times
L)=\{(x,0_S,0_L;0_{X^*},s^*,0_{L^*})|\, (x,0_S;0_{X^*},s^*)\in
\Lambda\}.
\end{eqnarray*}
But
\begin{eqnarray*}
\{(x,0_S;0_{X^*},s^*)\in \Lambda\}=\Lambda\cap T^*_X(X\times S)
=\emptyset
\end{eqnarray*}
where the last equality follows from the assumption that $f\circ g$
is transversal to $\xi$. Thus $(TR_1)$ is proved.

\hfill

\underline{Step 2.} Let us check that $g$ and $f^*\xi$ satisfy the
third condition of transversality $(TR_3)$. Namely let us check that
$WF(T')$ does not intersect the conormal bundle of the submanifold
$\PP_+(T^*_X(X\times S\times L))\subset\PP_{X\times S\times L}$.

We have $WF(T')=WF((df^*)_*p^*T)\subset (df^*)_*p^*\Gamma$. Let us
write down the latter space more explicitly.

Let us fix an element $(z_0,[z_0^*])$ of the space $\PP_Z=Z\times
\PP_+(Z^*)$. The fiber of $T^*\PP_Z$ over this point is identified
with $Z^*\times T^*_{[z_0^*]}(\PP_+(Z^*))$. Let us fix also a point
$l_0\in L$. The fiber of $T^*(\PP_Z\times L)$ over the point
$(z_0,[z_0^*],l_0)$ is equal to $Z^*\times
T^*_{[z_0^*]}(\PP_+(Z^*))\times L^*$, elements of which will be
denoted by $(z^*,\xi^*,l^*)$. In this notation we easily see that
the fiber of $p^*\Gamma$ over the point $(z_0,[z_0^*],l_0)$ is equal
to
\begin{eqnarray}\label{descrip-1}
(p^*\Gamma)|_{(z_0,[z_0^*],l_0)}=\{(z^*,\xi^*,0_{L^*})\in Z^*\times
T^*_{[z_0^*]}(\PP_+(Z^*))\times L^*|\, (z^*,\xi^*)\in
\Gamma|_{(z_0,[z_0^*])}\}.
\end{eqnarray}

Let us describe the fiber of $\tilde\Gamma:=(df^*)_*p^*\Gamma$ over
the point $(z_0,[z_0^*:0_{L^*}],l_0)$ (clearly the fiber of $\tilde
\Gamma$ may be non-empty only over points of such form). Let us
denote for this purpose
$$\nu\colon T^*_{[z_0^*:0_{L^*}]}(\PP_+(Z^*\times L^*))\to
T^*_{[z_0^*]}(\PP_+(Z^*))$$ the natural map, namely the dual map to
the differential of the imbedding $\PP_+(Z^*)\inj \PP_+(Z^*\times
L^*)$. Then using the description (\ref{descrip-1}) one easily sees
that
\begin{eqnarray*}\label{dscrip-2}
\tilde\Gamma|_{(z_0,[z_0^*:0_{L^*}],l_0)}=\\
\{(z^*,\eta^*,0_{L^*})\in Z^*\times
T^*_{[z_0^*:0_{L^*}]}(\PP_+(Z^*\times L^*))\times L^* |\,
(z^*,\nu(\eta^*))\in \{0\}\cup
\Gamma|_{(z_0,[z_0^*])}\}\backslash\{0\}.
\end{eqnarray*}

We want to show that the intersection of $\tilde \Gamma$ with the
conormal bundle of $\PP_+(T^*_X(X\times S\times L))$ is empty.
Observe that
\begin{eqnarray*}
\PP_+(T^*_X(X\times S\times L))=\{(x,0_S,0_L;[0_{X^*}:s^*:l^*])\}.
\end{eqnarray*}

Thus we will consider the case
\begin{eqnarray*}
z_0=(x_0,0_S)\in X\times S=Z,\\
l_0=0_L\in L,\\
\left[z_0^*\right]= [0_{X^*}:s_0^*]\in \PP_+(X^*\times S^*).
\end{eqnarray*}
The fiber of the conormal bundle of $\PP_+(T^*_X(X\times S\times
L))\subset \PP_Y$ over the point $(z_0,0_L,[z_0^*:0_{L^*}])$ (which
we identify with $(x_0,0_S,0_L;[0_{X^*}:s_0^*:0_{L^*}])$) is equal
to
\begin{eqnarray}\label{conormal-fiber}
\left\{(0_{X^*},s^*,l^*;\xi^*)|\, \mu(\xi^*)=0\right\}
\end{eqnarray}
where $\mu\colon T^*_{[0_{X^*}:s_0^*:0_{L^*}]}(\PP_+(X^*\times
S^*\times L^*))\to T^*_{[s_0^*:0_{L^*}]}(\PP_+(S^*\times L^*))$ is
the dual of the differential of the imbedding $\PP_+(S^*\times
L^*)\subset \PP_+(X^*\times S^*\times L^*)$. The intersection of the
set (\ref{conormal-fiber}) with the fiber of $\tilde \Gamma$ is
equal to
\begin{eqnarray}\label{intersection-1}
\{(0_{X^*},s^*,0_{L^*};\xi^*)|\,(0_{X^*},s^*, \nu(\xi^*))\in
\{0\}\cup \Gamma \mbox{ and } \mu(\xi^*)=0\}\backslash\{0\}.
\end{eqnarray}
It is easy to see that if $\mu(\xi^*)=0$ and $\xi^*\ne 0$, then
$\nu(\xi^*)\ne 0$. Hence if the set (\ref{intersection-1}) is
non-empty and a vector $(0_{X^*},s^*,0_{L^*};\xi^*)$ belongs to it,
then $(0_{X^*},s^*,\nu(\xi^*))\in\Gamma$. But then the latter vector
belongs to the intersection of the fibers over
$(x_0,0_S;[0_{X^*}:s_0])$ of $\Gamma$ and of the conormal bundle of
$\PP_+(T^*_X(X\times S))\subset \PP_{X\times S}$ (using
$\mu(\xi^*)=0$). But the last intersection is empty since the map
$f\circ g$ is assumed to be transversal to $\xi$ and thus satisfies
the third condition of transversality ($TR_3$), namely the
intersection of $\Gamma$ with the conormal bundle of
$\PP_+(T^*_X(X\times S))\subset \PP_{X\times S}$ is empty. Hence
Step 2 is completed.

\hfill

\underline{Step 3.} Let us check that $g$ and $f^*\xi$ satisfy
$(TR_2)$. Namely let us check that $WF(T')$ does not intersect the
conormal bundle of the submanifold $X\times_{Y}\PP_{Y}\subset
\PP_{Y}$. Recall that $Y=X\times S\times L$.

Let us fix an arbitrary point
$B_0:=(x_0,0_S,0_L;[x_0^*:s_0^*:l_0^*])\in X\times_Y\PP_Y$. We
identify the  fiber of $T^*\PP_Y$ over this point with $X^*\times
S^*\times L^*\times T^*_{[x_0^*:s_0^*:l_0^*]}(\PP_+(X^*\times
S^*\times L^*))$.

The fiber of the conormal bundle of $X\times_{Y}\PP_{Y}$ over the
point $B_0$ is equal to
\begin{eqnarray}\label{set-2.1}
\{(0_{X^*},s^*,l^*;0)\}\subset (X^*\times S^*\times L^*)\times
T^*_{[x_0^*:s_0^*:l_0^*]}(\PP_+(X^*\times S^*\times L^*))
\end{eqnarray}
where the last zero in the first set belongs to
$T^*_{[x_0^*:s_0^*:l_0^*]}(\PP_+(X^*\times S^*\times L^*))$. As in
Step 2, the fiber of $\tilde \Gamma=(df^*)_*p^*\Gamma$ over $B_0$ is
equal to
\begin{eqnarray}\label{set-2.2}
\left\{(x^*,s^*,0_{L^*};\xi^*)|\, (x^*,s^*;\nu(\xi^*))\in \{0\}\cup
\Gamma\right\}.
\end{eqnarray}
Hence the intersection of the sets (\ref{set-2.1}) and
(\ref{set-2.2}) is equal to
\begin{eqnarray*}
\{(0_{X^*},s^*,0_{L^*};0)| \, (0_{X^*},s^*;0)\in \{0\}\cup
\Gamma\}=\{(0_{X^*},s^*,0_{L^*};0)| \, (0_{X^*},s^*;0)\in \Gamma\}.
\end{eqnarray*}
But the last set is empty since the set $\{(0_{X^*},s^*;0)\}\cap
\Gamma$ is empty. Indeed by the assumption $f\circ g$ is transversal
to $\xi$, and hence the second condition of transversality ($TR_2$)
is satisfied: the intersection of $\Gamma$ with the conormal bundle
of $X\times _{X\times S}\PP_{X\times S}\subset \PP_{X\times S}$ is
empty. Step 3 is proved.

\hfill

\underline{Step 4.} It remains to show that $(f\circ
g)^*(\xi)=g^*(f^*(\xi))$.

First let us prove it in the special case $X=Z$, i.e. $S=\{0\}$. In
this case we have
$$Z\overset{g}{\to}X\times L\overset{f}{\to}Z,$$
and $f\circ g=Id_Z$. For any closed conic subsets $\Lam\subset
T^*Z\backslash \underline{0},\, \Gamma \subset T^*\PP_Z\backslash
\underline{0}$, we denote, as in Step 2, by
$\tilde\Gamma:=(df^*)_*p^*\Gamma$. Then the maps
$$f^*\colon V^{-\infty}_{\Lam,\Gamma}(Z)\to
V^{-\infty}_{f^*\Lam,\tilde\Gamma}(X\times L),\, g^*\colon
V^{-\infty}_{f^*\Lam,\tilde\Gamma}(X\times L)\to V^{-\infty}(Z)$$
are sequentially continuous operators. In particular taking
$\Lam=T^*Z\backslash\underline{0},
\Gamma=T^*\PP_Z\backslash\underline{0}$ we deduce that
$$g^*\circ f^*\colon V^{-\infty}(Z)\to V^{-\infty}(Z)$$
is a sequentially continuous operator.

To show that $g^*\circ f^*=(f\circ g)^*(=Id)$ it is enough to check
that these operators coincide on characteristic functions of compact
submanifolds with corners; indeed if two sequentially continuous linear operators, defined on $V^{-\infty}(Z)$ as a source space,
coincide on constructible functions then they coincide by Corollary 1.3 in \cite{alesker-sequential-closure}.
Thus let $P\subset Z$ be a compact submanifold with corners.
Then by Proposition \ref{P:pull-submer-gener-constr}
$$f^*(\Xi_{\cp}(\One_P))=\Xi_{\cp}(\One_{P\times L}).$$
But observe that $Z$ is transversal to $P\times L$. Hence by
Proposition \ref{P:pull-back-constuct} we have
$g^*(\Xi_{\cp}(\One_{P\times L}))=\Xi_{\cp}(\One_{(P\times L)\cap
Z})=\Xi_\cp(\One_Z)$. Hence $g^*\circ f^*=Id$.

Now let us consider the general case. Let $\xi$ correspond to a pair
of currents $(C,T)$. As in Steps 1-3 let us denote
$\Lambda:=WF(C),\, \Gamma:=WF(T)$. From the computations in Steps
1-3 and from the behavior of wave fronts under pull-back and
push-forward (see Propositions \ref{P:lifting} and
\ref{P:wf-push-forward}) it follows that $(f\circ g)^*$ and
$g^*\circ f^*$ are sequentially continuous linear operators
$$\vlg(Z)\to V^{-\infty}(X).$$ By
Lemma \ref{C:approxim} it is enough to show that these operators
coincide on the subspace of smooth valuations. Thus we may assume
that $\xi$ is a smooth valuation.

Let us decompose $g$ as $X\overset{h}{\to}Z\overset{k}{\to} Y$ where
$h(x)=(x,0_S)$, $k(z)=(z,0_L)$. Both maps $h,k$ are imbeddings.
\begin{claim}\label{claim-on-prop}
These maps $h,k$ and the valuation $f^*\xi\in V^{-\infty}(Y)$
satisfy the assumptions of Proposition \ref{P:imbed-compos}.
\end{claim}
Let us postpone the proof of this claim and finish the proof of our
proposition. By Proposition \ref{P:imbed-compos},
$g^*(f^*\xi)=h^*(k^*(f^*\xi))$. By what we have just shown
$k^*(f^*(\xi))=\xi$. Hence $g^*(f^*(\xi))=h^*\xi$. On the other hand
$f\circ g=h$, hence $(f\circ g)^*\xi=h^*\xi$. Thus to finish the
proof of the proposition it remains to prove Claim
\ref{claim-on-prop}.

\hfill

{\bf Proof of Claim \ref{claim-on-prop}.} Let $f^*\xi\in
V^{-\infty}(Y)$ be given by a pair of currents $(\tilde C,\tilde
T)$. Since $\xi$ is smooth, $\tilde C=f^*C$ is also smooth. Thus
$WF(\tilde C)=\emptyset$. The wave front of $\tilde T$ is contained
in the conormal bundle of $Y\times_Z\PP_Z\subset \PP_Y$,  thus let
us denote
$$M:=T^*_{Y\times_Z\PP_Z}(\PP_Y)\backslash\underline{0}\subset
T^*\PP_Y\backslash\underline{0}.$$ We have
\begin{eqnarray}\label{descr-claim1}
Y\times_Z\PP_Z=\{(z,l;[z^*:0_{L^*}])\}\subset \PP_Y=\PP_{Z\times L}.
\end{eqnarray}

\underline{Step 1.} First we will check that $k\colon Z\to Y$ is
transversal to $(\emptyset,M)$, i.e. we have to verify the three
conditions of transversality ($TR_1$)-($TR_3$) (see Definition
\ref{D:def-sets-transversal}).

The first condition ($TR_1$) is satisfied for trivial reasons. Let
us check ($TR_2$), namely
\begin{eqnarray}\label{tr2-new}
M\cap T^*_{Z\times_Y\PP_Y}(\PP_Y)=\emptyset.
\end{eqnarray}
We have
\begin{eqnarray*}
Z\times_Y\PP_Y=Z\times \{0_{L}\}\times \PP_+(Z^*\times L^*).
\end{eqnarray*}

Fix points $z_0\in Z,\, [z_0^*]\in \PP_Z$. Denote
$A_0:=(z_0,0_L;[z_0^*:0_{L^*}])\in Y\times_Z\PP_Z$. The fiber of
$M\cup\underline{0}$ over $A_0$ is equal to
\begin{eqnarray}\label{E:fiber-M}
M|_{A_0}\cup \{0\} =\{0_{Z^*}\}\times \{0_{L^*}\}\times
T^*_{\PP_+(Z^*\times 0_{L^*})}(\PP_+(Y^*))\big|_{[z_0^*:0_{L^*}]}.
\end{eqnarray}

The fiber of $T^*_{Z\times_Y\PP_Y}(\PP_Y)$ over $A_0$ is equal to
\begin{eqnarray}\label{fiber-tzy}
T^*_{Z\times_Y\PP_Y}(\PP_Y)\big|_{A_0}=\{0_{Z^*}\}\times L^*\times
\{0\}
\end{eqnarray}
where the last zero on the right hand side belongs to
$T^*_{[z_0^*:0_{L^*}]}(\PP_+(Y^*))$. Next from (\ref{E:fiber-M}) and
(\ref{fiber-tzy}) we see that
$$M\cap T^*_{Z\times_Y\PP_Y}(\PP_Y)=\emptyset.$$ Thus the second condition
of transversality ($TR_2$) is checked.

Let us check the third condition ($TR_3$), namely
$$M\cap T^*_{\PP_+(T^*_ZY)}(\PP_Y)=\emptyset.$$

We have
\begin{eqnarray*}
\PP_+(T^*_ZY)=\{(z,0_L;[0_{Z^*}:l^*])\}\subset \PP_Y.
\end{eqnarray*}
It is clear that $\PP_+(T^*_ZY)\cap (Y\times_Z\PP_Z)=\emptyset$.
Hence trivially
$$T^*_{\PP_+(T^*_ZY)}(\PP_Y)\cap M=\emptyset,$$
thus the condition ($TR_3$) is satisfied.
\def\talp{\alp_1}
\def\tbeta{\beta_1}
\def\teps{\eps_1}
\def\tdelta{\delta_1}

\underline{Step 2.} Let us check that $h\colon X\to Z$ is
transversal to $k^*(f^*\xi)$. It suffices to check that $h\colon
X\to Z$ is transversal to $(\emptyset,\beta_{1*}\alp_1^*M)$ where
$$\PP_Z\overset{\beta_1}{\leftarrow}
\widetilde{Z\times_Y\PP_Y}\overset{\alp_1}{\to}\PP_Y$$ is the
diagram analogous to (\ref{dia-1}) but with $Z\overset{k}{\to}Y$
instead of $X\to Y$. We will prove a stronger statement:
$\beta_{1*}\alp_1^*M=\emptyset$.

Similarly to (\ref{dia-2}), let us decompose $\talp$ as
\begin{eqnarray*}
\widetilde{Z\times_Y\PP_Y}\overset{\teps}{\to}Z\times_Y\PP_Y\overset{\tdelta}{\to}
\PP_Y
\end{eqnarray*}
where $\teps$ is the oriented blow up map along $\PP_+(T^*_ZY)$, and
$\tdelta$ is the imbedding map.

Let us consider subset $J\subset \PP_Y$ defined by
$$J:=(Y\times_Z\PP_Z)\cap(Z\times_Y\PP_Y)=\{(z,0_L;[z^*:0_{L^*}])\}$$
where the intersection is obviously transversal. It is clear that
\begin{eqnarray*}
\delta_1^*M=T^*_J(Z\times_Y\PP_Y)\backslash\underline{0}.
\end{eqnarray*}

Next the oriented blow up map
$$\eps_1\colon \widetilde{Z\times_Y\PP_Y}\to Z\times_Y\PP_Y$$
is along the submanifold $\PP_+(T^*_ZY)$ which obviously does not
intersect $J$. It follows that
$$\eps_1^*\delta_1^*M=T^*_{\eps_1^{-1}(J)}(\widetilde{Z\times_Y\PP_Y})\backslash\underline{0},$$
and $\eps_1$ is a diffeomorphism of a neighborhood $\co$ of
$\eps_1^{-1}(J)$ onto its image. Hence while computing
$\beta_{1*}\eps^*_1\delta_1^*M$ we may consider the restriction of
$\beta_1$ to this neighborhood $\co$. This neighborhood $\co$ can be
taken to be
$$\co:=\{(z,0_L;[z^*:l^*])|\, z^*\ne 0\}.$$
Next $\beta_1\colon \co\to \PP_Z$ is given by
$$\beta_1((z,0_L;[z^*:l^*]))=(z;[z^*]).$$ Also
$$\eps_1^{-1}(J)=\{(z,0_L;[z^*:0_{L^*}])\}.$$
Then clearly the restriction of $\beta_1$ to $\eps_1^{-1}(J)$ is a
diffeomorphism $\eps^{-1}(J)\tilde\to \PP_Z$. This and the fact that
$\beta_1\colon \co\to \PP_Z$ is a fibration imply that
$$\beta_{1*}(T^*_{\eps_1^{-1}(J)}\co\backslash\underline{0})=\emptyset.$$
Hence
$\beta_{1*}\alp_{1*}M=\beta_{1*}\eps_1^*\delta_1^*M=\emptyset$. Thus
Step 2 is completed.

\underline{Step 3.} It remains to show that $g=k\circ h$ is
transversal to $(\emptyset, M)$. The condition ($TR_1$) is satisfied
trivially.

Let us check ($TR_2$). Let us suppose in the contrary that there
exists a point $A_0\in X\times_Y\PP_Y$ and $\xi\in M|_{A_0}$ which
belongs also to the fiber over $A_0$ of the conormal bundle
$T^*_{X\times_Y\PP_Y}(\PP_Y)$. The point $A_0$ has the form
$$A_0=(x_0,0_S,0_L;[x_0^*:s_0^*:0_{L^*}]).$$ Denote also
$$Q_0:=[x_0^*:s_0^*:0_{L^*}]\in \PP_+(X^*\times S^*\times
0_{L^*})\subset \PP_+(Y^*).$$

Using  (\ref{E:fiber-M}) one sees that $\xi$ must have the form
$$\xi=(0_{X^*},0_{S^*},0_{L^*};\zeta) \mbox{ with } \zeta\in
\left(T^*_{\PP_+(X^*\times S^*\times
0_{L^*})}(\PP_+(Y^*))\right)\big|_{Q_0}.$$ But the assumption that
$\xi\in \left(T^*_{X\times_Y\PP_Y}(\PP_Y)\right)\big|_{A_0}$ implies
that $\zeta=0$. Hence $\xi=0$ which is a contradiction. Thus
($TR_2$) is proved.

Let us prove ($TR_3$). Let us suppose in the contrary that there
exists a point $B_0\in \PP_+(T^*_XY)$ and $\eta\in M|_{B_0}$ such
that $\eta$ also belongs to the fiber over $B_0$ of the conormal
bundle $T^*_{\PP_+(T^*_XY)}(\PP_Y)$. Notice also that
$$\PP_+(T^*_XY)=\{(x,0_S,0_L;[0_{X^*}:s^*:l^*])\}.$$

Necessarily $B_0$ has the form
$$B_0=(\tilde x_0,0_S,0_L;[0_{X^*}:\tilde s_0^*:0_{L^*}]).$$
Let us denote $R_0:=[0_{X^*}:\tilde s_0^*:0_{L^*}]\in \PP_+(Y^*)$.
In this notation $\eta$ must have the form
$$\eta=(0_{X^*},0_{S^*},0_{L^*};\tau)$$
with $\tau\in \left(T^*_{\PP_+(X^*\times S^*\times
0_{L^*})}(\PP_+(Y^*))\right)\big|_{R_0}\cap
\left(T^*_{\PP_+(0_{X^*}\times S^*\times
L^*)}(\PP_+(Y^*))\right)\big|_{R_0}.$ But the last intersection is
equal to $\{0\}$. This is a contradiction, thus ($TR_3$) is proved.

Hence Claim \ref{claim-on-prop} is proved. Hence Proposition
\ref{P:compos-immer-submer} is proved completely too. \qed

\subsection{Push-forward of generalized valuations under
submersions.}\label{Ss:pushforward-gen-submersions} Let $f\colon
X\to Y$ be a proper submersion. In Section \ref{Ss:push-submersion}
we have defined a push-forward map on smooth valuations
$$f_*\colon V^\infty(X)\to V^\infty(Y)$$
which is a continuous operator. In this section we extend this map
to a {\itshape partially} defined map on generalized valuations.
This extension will be important for applications to integral
geometry in Section \ref{S:IG} below.

Ideologically, the push-forward map $f_*$ is dual to the pull-back
map $f^*$. Proposition \ref{P:pull-submer-gener-expr} implies that
for the submersion $f\colon X\to Y$ we have the continuous pull-back
map
$$f^*\colon V^\infty(Y)\to V^{-\infty}_{\emptyset,M}(X)$$
where $M\subset T^*\PP_X\backslash\underline{0}$ is the conormal
bundle of $X\times_Y\PP_Y\subset \PP_X$ with the zero section
removed. (In this case $f^*$ is topologically continuous and not just sequentially continuous.
This is due to the fact that the source space $V^\infty(Y)$ is metrizable.)
If in addition $f$ is proper then $f^*$ maps compactly
supported valuations to compactly supported:
$$f^*\colon V^\infty_c(Y)\to V^{-\infty}_{\emptyset,M,c}(X).$$

Next in \cite{alesker-bernig} there was defined a partial product on
generalized valuations on a manifold $X$ which will be needed now.
Let $\xi_1,\xi_2\in V^{-\infty}(X)$ be generalized valuations.
Assume that the wave fronts $(\Lam_i,\Gamma_i)$, $i=1,2$, of the
currents corresponding to $\xi_i$ satisfy some explicitly written
condition of disjointness (see Theorem \ref{thm_prod_gen}, or
(\ref{pro-a})-(\ref{pro-e}) below in this section). Then one can
define the product $\xi_1\cdot\xi_2\in V^{-\infty}(X)$. Moreover
this product defines a bilinear jointly sequentially continuous map
$$V^{-\infty}_{\Lam_1,\Gamma_1}(X)\times
V^{-\infty}_{\Lam_2,\Gamma_2}(X)\to V^{-\infty}(X)$$ which coincides
with the usual product on smooth valuations.

Now in order to define the push-forward map $f_*$ for a proper
submersion $f\colon X\to Y$, let us fix an arbitrary closed conic
subsets
$$\Lam\subset T^*X\backslash\underline{0},\, \Gamma\subset
T^*\PP_X\backslash\underline{0}$$ such that $(\Lam,\Gamma)$ and
$(\emptyset,M)$ satisfy the above mentioned condition of
disjointness. Let us define eventually
$$f_*\colon V^{-\infty}_{\Lam,\Gamma}(X)\to V^{-\infty}(Y)$$
as follows
$$<f_*\xi,\phi>:=\int_X \xi\cdot f^*\phi$$
for any $\xi\in V^{-\infty}_{\Lam,\Gamma}(X)$ and any test valuation
$\phi\in V^\infty_c(Y)$. The above discussion implies that we get a sequentially
continuous linear operator $f_*$.

By the discussion in Section \ref{Ss:submersions} it is clear that
the restriction of this $f_*$ to smooth valuations coincides with
the push-forward defined in Section \ref{Ss:push-submersion}. Thus
there is no abuse of notation.

\hfill






Let us remind again the condition of disjointness of wave fronts.
Assume that $\Lam_i\subset T^*Y\backslash \underline{0}$,
$\Gamma_i\subset T^*\PP_Y\backslash \underline{0}$, $i=1,2$, be two
closed conic subsets. One can define a jointly sequentially continuous
bilinear map
$$\vmi_{\Lam_1,\Gamma_1}(Y)\times \vmi_{\Lam_2,\Gamma_2}(Y)\to
\vmi(Y)$$ extending the product on smooth valuations, provided the
following conditions are satisfied:
\begin{eqnarray}\label{pro-a}
\Lam_1\cap s(\Lam_2)=\emptyset,\\\label{pro-b} \Gamma_1\cap
s(\pi^*_Y\Lam_2)=\emptyset,\\\label{pro-c} \Gamma_2\cap
s(\pi^*_Y\Lam_1)=\emptyset,\\\label{pro-d} \mbox{if }
(y,[\xi_i],u_i,0)\in \Gamma_i, i=1,2, \mbox{ then } u_1\ne -u_2,\\
\mbox{if } (y,[\xi])\in\PP_Y \mbox{ and } (u,\eta_1)\in
\Gamma_1|_{(y,[\xi])},\, (-u,\eta_2)\in \Gamma_2|_{(y,[-\xi])},
\mbox{ then }\nonumber\\\label{pro-e} d\theta^*(0,\eta_1,\eta_2)\ne
(0,l,-l)\in T^*_{(y,[\xi],[\xi])}(\PP_Y\times_Y\PP_Y) \\\mbox{ where
} \theta\colon \PP_Y\times_Y\PP_Y\to \PP_Y\times_Y\PP_Y \mbox{ is
defined by }
\theta(y,[\zeta_1],[\zeta_2])=(y,[\zeta_1],[-\zeta_2])\nonumber.
\end{eqnarray}

Now let us define the push-forward map $$f_*\colon
V^{-\infty}_{\Lam,\Gamma}(X)\to V^{-\infty}(Y)$$ for a proper
submersion $f\colon X\to Y$ and closed conic subsets $\Lam\subset
T^*X\backslash\underline{0}$, $\Gamma\subset T^*\PP_X\backslash
\underline{0}$ such that
\begin{eqnarray}\label{D:push-forw-gener-val-submer}
\mbox{the pairs } (\Lam,\Gamma) \mbox{ and }
(\emptyset,T^*_{X\times_Y\PP_Y}(\PP_X)) \mbox{ satisfy
(\ref{pro-a})-(\ref{pro-e})},
\end{eqnarray}
as follows:
$$<f_*\phi,\psi>=\int_Y\phi\cdot f^*\psi$$
for any $\phi\in V^{-\infty}_{\Lam,\Gamma}(X),\, \psi\in
V^\infty_c(Y)$. The previous discussion implies readily that this
map $f_*\colon V^{-\infty}_{\Lam,\Gamma}(X)\to V^{-\infty}(Y)$ is a
sequentially continuous linear operator.

Under the same assumption on $(\Lam,\Gamma)$ one defines in the same
way the push-forward on compactly supported valuations
$$f_*\colon V^{-\infty}_{\Lam,\Gamma,c}(X)\to V^{-\infty}_c(Y)$$
if $f\colon X\to Y$ is a submersion, but not necessarily proper.

In order to study further properties of the push-forward map, let us
prove some technical results.

\hfill


Assume now that $i\colon X\inj Y$ be a smooth {\itshape imbedding}.
Recall that we have the diagram (see (\ref{def-alp-beta}))
$$\PP_X\overset{\beta}{\leftarrow}\widetilde{X\times
_Y\PP_Y}\overset{\alp}{\to}\PP_Y.$$ Recall that by Proposition
\ref{P:push-forward-gener-imbed-currents} for $\phi\in \vmi_c(X)$
given by a pair of currents $(C_\phi,T_\phi)$ the push-forward
$i_*\phi$ is given by the pair $(0,\alp_*\beta^*T_\phi)$.

Recall also that, by the Definition of the pull-back under
immersions from Section \ref{Ss:pull-gener-immer}, if
$\psi\in\vmi(Y)$ is given by a pair $(C_\psi,T_\psi)$ and the
imbedding $i$ is transversal to $\psi$ then the pull-back $i^*\psi$
is given by the pair of currents $(i^*C_\psi,\beta_*\alp^*T_\psi)$.

Let us denote also
\begin{eqnarray*}
\Lam_\phi:= WF(C_\phi), \Lam_\psi:= WF(C_\psi),\\
\Gamma_\phi:=WF(T_\phi),  \Gamma_\psi:=WF(T_\psi).
\end{eqnarray*}

\begin{lemma}\label{LL:1}
Using the above notation let us assume that the imbedding $i\colon
X\inj Y$ is transversal to $\psi\in \vmi(Y)$ (in the sense of
Definition \ref{D:map-val-transversal}), the pairs of subsets
$(\Lam_\phi,\Gamma_\phi)$ and
$(i^*\Lam_\psi,\beta_*\alp^*\Gamma_\psi)$ satisfy
(\ref{pro-a})-(\ref{pro-e}), and also the pairs of sets
$(\emptyset,\alp_*\beta^*\Gamma_\phi)$ and $(\Lam_\psi,\Gamma_\psi)$
satisfy (\ref{pro-a})-(\ref{pro-e}).

Then $i_*\phi\cdot \psi$ and $\phi\cdot i^*\psi$ are well defined,
and
\begin{eqnarray}\label{E:adjoint-imbed}
\int_Y i_*\phi\cdot \psi=\int_X \phi\cdot i^*\psi.
\end{eqnarray}
\end{lemma}
{\bf Proof.} The well-definedness of the expressions is clear. This
implies that the two bilinear maps
$\vmi_{\Lam_\phi,\Gamma_\phi,c}(X)\times
\vmi_{\Lam_\psi,\Gamma_\psi}(Y)\to \CC$ given respectively by
$(\mu,\nu)\mapsto \int_Yi_*\mu\cdot \nu$ and $(\mu,\nu)\mapsto
\int_X \mu\cdot i^*\nu$ are jointly sequentially continuous. Thus in order to
show that they coincide it is enough to assume that $\psi$ is
smooth. But for smooth $\psi$ the equality (\ref{E:adjoint-imbed})
is satisfied by the definition of $i_*$ from Section
\ref{Ss:push-gener-immer}. \qed

\hfill

\begin{lemma}\label{LL:2}
Let $Y=X\times Z,\, f\colon X\times Z\to X$ be the natural
projection. Let $i\colon X\inj X\times Z$ be a closed imbedding of
the form $i(x)=(x,z_0)$ where $z_0\in Z$ is a fixed point. Let
$\phi\in \vmi_c(X)$ be a generalized valuation with compact support.
Let $\xi\in V^\infty(X)$ be a smooth valuation.

Then $\phi$, $\psi:=f^*\xi$, and $i\colon X\inj Y$ satisfy the
assumptions of Lemma \ref{LL:1}. Consequently the generalized
valuations $i_*\phi \cdot f^*\xi$ and $\phi\cdot i^*(f^*\xi)$ are
well defined, and
$$\int_Yi_*\phi \cdot f^*\xi=\int_X\phi\cdot i^*(f^*\xi).$$
\end{lemma}
{\bf Proof.} Let $\psi:=f^*\xi$ be given by the pair of currents
$(C_\psi,T_\psi)$. Since $\xi$ is smooth we have
\begin{eqnarray}\label{E:sept1}
\Lam_\psi=\emptyset,\\\label{E:sept2} \Gamma_\psi\subset
T^*_{Y\times_X\PP_X}(\PP_Y).
\end{eqnarray}
The condition (\ref{E:sept2}) implies readily that
$\beta_*\alp^*\Gamma_\psi=\emptyset$. Hence the first assumption of
Lemma \ref{LL:1} follows: the pairs $(\Lam_\phi,\Gamma_\phi)$ and
$(i^*\Lam_\psi,\beta_*\alp^*\Gamma_\psi)=(\emptyset,\emptyset)$
satisfy (\ref{pro-a})-(\ref{pro-e}).

\hfill

Let us check now the second assumption, namely that the pairs
$(\emptyset, \alp_*\beta^*\Gamma_\phi)$ and
$(\Lam_\psi,\Gamma_\psi)=(\emptyset, \Gamma_\psi)$ satisfy
(\ref{pro-a})-(\ref{pro-e}). The conditions
(\ref{pro-a})-(\ref{pro-c}) are satisfied trivially.

Let us check (\ref{pro-d}). Assume that $(y,[\xi_1],u_1,0)\in
\alp_*\beta^*\Gamma_\phi$, $(y,[\xi_2],u_2,0)\in \Gamma_\psi$. But
(\ref{E:sept2}) implies that $u_2=0$ which is a contradiction. Thus
(\ref{pro-d}) is proved.

Let us check (\ref{pro-e}). For simplicity of the notation we may
and will assume that $X,Z$ are vector spaces, $z_0=0_Z\in Z$. Assume
that
\begin{eqnarray*}
(y,[\xi])\in \PP_Y,\\
(u,\eta_1)\in(\alp_*\beta^*\Gamma_\phi)|_{(y,[\xi])},\\
(-u,\eta_2)\in \Gamma_\psi|_{(y,[-\xi])}.
\end{eqnarray*}
Clearly $y=(x,0_Z)$. Next (\ref{E:sept2}) implies that $u=0$ and
$\xi=(x^*,0_{Z^*})$ with $x^*\in X^*$. Denote
$(0,\eta_2'):=ds^*((0,\eta_2))\in
(T^*_{Y\times_X\PP_X}(\PP_Y))|_{(x,[\xi])}$. Thus we have to show
that $\eta_1\ne-\eta_2'$.

Otherwise we would have that
\begin{eqnarray}\label{E:sept5}
(0,\eta_1)\in \left(\alp_*\beta^*\Gamma_\phi\cap
T^*_{Y\times_X\PP_X}(\PP_Y)\right)\big|_{(y,[\xi])}.
\end{eqnarray}
Since $[\xi]=[x^*:0_{Z^*}]$, there exists a neighborhood $U\subset
\PP_Y$ of $(y,[\xi])$ such that $\alp\colon \alp^{-1}(U)\to U$ is a
diffeomorphism. Thus we will simply identify $\alp^{-1}(U)$ with $U$
via $\alp$ (recall that $\alp$ is the oriented blow up map along
$\PP_+(T^*_XY)\subset \PP_Y$). Then there exists a cotangent vector
$\lam\in (T^*\PP_X)|_{(x,[x^*])}$ such that
$(0,\eta_1)=d\beta^*(\lam)$. Let us denote by $F$ the tangent space
at $(y,[\xi])$ to the fiber $\beta^{-1}((x,[x^*]))$.  Then
$(0,\eta_1)$ vanishes on $F$. Also by (\ref{E:sept5}) $(0,\eta_1)$
vanishes on the tangent space of $Y\times_X\PP_X$. The latter space
is clearly transversal to $F$. Hence $(0,\eta_1)$ vanishes, which is
a contradiction. \qed

\hfill

\begin{corollary}\label{CC:3}
As in Lemma \ref{LL:2}, let $Y=X\times Z$, $f\colon X\times Z\to X$
be the natural projection. Let $i\colon X\inj X\times Z$ be a closed
imbedding of the form $[x\mapsto (x,z_0)]$ where $z_0\in Z$ is a
fixed point. Let $\phi\in \vmi_c(X)$ be a generalized valuation with
compact support. Let $\xi\in \vi(X)$ be a smooth valuation. Then
$$\int_Y i_*\phi\cdot f^*\xi=\int_X\phi\cdot \xi$$
where all the products are automatically well defined.
\end{corollary}
{\bf Proof.} By Lemma \ref{LL:2} the products are well defined, and
$$\int_Y i_*\phi\cdot f^*\xi=\int_X \phi\cdot i^*(f^*\xi).$$
But Proposition \ref{P:compos-immer-submer} implies immediately that
$i^*(f^*\xi)=\xi$. Corollary is proved. \qed


\hfill

Our next goal is to show that the pushforward under a submersion of
a constructible function (considered as a generalized valuation) is
equal to integration of this function with respect to the Euler
characteristic along the fibers (see Ch. 9 of
\cite{kashiwara-schapira} for the detailed discussion of this
notion). However we will do this explicitly only under rather
restrictive assumptions. This generality will be sufficient for some
integral geometric applications discussed below, but from the point
of view of the general theory it does not seem to be very
satisfactory. Our first result in this direction the following
lemma.

\begin{lemma}\label{LL:4}
Let $X,Z$ be smooth manifolds without boundary, let $Y=X\times Z$.
Let $f\colon X\times Z\to X$ be the natural projection. Let $P\in
\cp(Z)$. Then the wave fronts of the generalized valuation
$\Xi_{\cp}(\One_{X\times P})\in \vmi(Y)$ satisfy the condition of
transversality (\ref{D:push-forw-gener-val-submer}) (note that here
the roles of $X$ and $Y$ are interchanged) and
$$f_*(\Xi_{\cp}(\One_{X\times P}))=\chi(P)\cdot \chi.$$
\end{lemma}
{\bf Proof.} We have $N(X\times P)=Y\times_Z N(P)\subset
Y\times_Z\PP_Z\subset \PP_Y$. Clearly $N(X\times P)$ is disjoint
from $Y\times_X\PP_X$. Hence the push-forward
$f_*(\Xi_\cp(\One_{X\times P}))$ is well defined, i.e. the condition
of transversality (\ref{D:push-forw-gener-val-submer}) is satisfied.

Triangulating $P$ we may assume that $Z=\RR^l$ and $P$ is a linear
simplex. For $N\in \NN$ let us denote $P_N:=\frac{1}{N}\cdot P$. Also let
us denote for brevity
\begin{eqnarray*}
\Lam:=T^*Y\backslash \underline{0},\\
\Gamma:=(Y\times_Z\PP_Z)\times_{\PP_Y}(T^*\PP_Y\backslash
\underline{0})\subset T^*\PP_Y.
\end{eqnarray*}
It is easy to see that the pairs of sets $(\Lam,\Gamma)$ and
$(\emptyset, T^*_{Y\times_X\PP_X}\PP_Y\backslash \underline{0})$
satisfy the conditions (\ref{pro-a})-(\ref{pro-e}).

Below we will show that
\begin{eqnarray}\label{E:converge-may}
\Xi_\cp(\One_{X\times P_N})\to \Xi_\cp(\One_{X\times\{0_Z\}})
\mbox{ in } \vmi_{\Lam,\Gamma}(Y) \mbox{ as } N\to \infty.
\end{eqnarray}
First let us finish the proof of lemma. Evidently
$f_*(\Xi_\cp(\One_{X\times P_N}))=f_*(\Xi_\cp(\One_{X\times P}))$
for any $N$. Hence by (\ref{E:converge-may})
$f_*(\Xi_\cp(\One_{X\times P}))=f_*(\Xi_\cp(\One_{X\times
\{0_Z\}}))$. Let $i\colon X\inj X\times Z$ be the imbedding
$i(x)=(x,0_Z)$. Let $\phi:=\Xi_\cp(\One_X)\in V^{-\infty}(X)$.
Applying Corollary \ref{CC:3} to $\phi,f,$ and $i$ we deduce the
result (formally we may assume that $X$ is compact).

Thus it remains to prove (\ref{E:converge-may}). It is clear that
$[[X\times P_N]]\to [[X\times \{0_Z\}]]$ in
$C^{-\infty}_\Lam(Y)=C^{-\infty}(Y)$. Next $N(X\times
P_N)=X\times N(P_N)\subset \PP_Y$. By \cite{part3}, Lemma
2.1.3, we have $[[N(P_N)]]\to [[N(\{0_Z\})]]$ in $\cd_{\dim
Z-1}(\PP_Z)$. Hence $[[X\times N(P_N)]]\to [[X\times
N(\{0_Z\})]]=[[N(X\times \{0_Z\})]]$ in $\cd_{\dim_Y-1}(\PP_Y)$.

It remains to show a seemingly more precise statement that
$[[N(X\times P_N)]]\to [[N(X\times \{0_Z\})]]$ in
$C^{-\infty}_\Gamma(\PP_Y,\Omega^{\dim Y})$, which is obvious since
for any $N\in \NN$ one has an inclusion $N(X\times P_N)\subset
Y\times_Z\PP_Z$.
Lemma is proved. \qed

\hfill

\begin{proposition}\label{PP:5}
Let $X=\RR^m, Z=\RR^{n-m}, Y=X\times Z=\RR^n$. Let $f\colon Y\to X$
be the natural projection. Let $A\subset Y$ be a compact convex set
with smooth boundary and strictly positive Gauss curvature. Then the
push-forward $f_*(\Xi_\cp(\One_A))$ is well defined and is equal to
$\Xi_\cp(\One_{f(A)})$.
\end{proposition}
{\bf Proof.} Clearly we may and will assume that $n>m$. The Gauss
map $\pt A\to S^{n-1}$ is a diffeomorphism. This implies that $N(A)$
is transversal to $Y\times_X\PP_X$. Hence the wave fronts of
$\Xi_\cp(\One_A)$ satisfy the condition
(\ref{D:push-forw-gener-val-submer}) (with $X$ and $Y$
interchanged). Hence $f_*(\Xi_\cp(\One_A))$ is well defined.

For $\eps\in \{1,\frac{1}{2},\frac{1}{3},\frac{1}{4},\dots\}\cup \{0\}$
 let us denote by $f_\eps$ the linear map $\RR^n\to
\RR^n$ which is given by $f_\eps(y,z)=(y,\eps z)$. For $\eps >0$ let
us denote by  $g_\eps\colon \PP_{\RR^n}\to \PP_{\RR^n}$ the induced
map. Denote also $A_\eps:=f_\eps(A)$. Notice that $A_0=f(A)\times
\{0_Z\}$.

Let us fix a small open conic neighborhood $\cu\subset
T^*\PP_Y\backslash \underline{0}$ of the set
$T^*_{Y\times_X\PP_X}\PP_Y\backslash \underline{0}$ (it will be
clear from the discussion below how small $\cu$ should be). Denote
by $\Gamma$ the complement of $\cu$ in $T^*\PP_Y\backslash
\underline{0}$. Denote $\Lam:=T^*Y\backslash \underline{0}$. First
it is easy to see that the pair of sets $(\Lam,\Gamma)$ satisfies
the condition (\ref{D:push-forw-gener-val-submer}). We will show
that
\begin{eqnarray}\label{sept10}
\Xi_\cp(\One_{A_\eps})\to \Xi_\cp(\One_{A_0}) \mbox{ in }
\vmi_{\Lam,\Gamma}(Y) \mbox{ as } \eps\to +0.
\end{eqnarray}
This, sequential continuity of $f_*\colon V^{-\infty}_{\Lam,\Gamma}(Y)\to V^{\infty}$,
and Corollary \ref{CC:3} will imply our lemma easily.

\hfill

Thus it remains to prove (\ref{sept10}). It is clear that $A_\eps\to
A_0$ in the Hausdorff metric. Hence by \cite{part3}, Lemma 2.1.3,
$[[N(A_\eps)]]\to [[N(A_0)]]$ in $\cd_{n-1}(\PP_Y)$. It is also
clear that $[[A_\eps]]\to [[A_0]]=0$ in $\cd_n(Y)=C^{-\infty}(Y)$.
Hence $\Xi_\cp(\One_{A_\eps})\to \Xi_\cp(\One_{A_0})$ in $\vmi(Y)$.

Thus it remains to show that $[[N(A_\eps)]]\to [[N(A_0)]]\mbox{ in }
C^{-\infty}_{\Gamma}(\PP_Y,\Omega^n) \mbox{ as } \eps\to +0$.
Let us fix a conic subset $K\subset \cu$ such that $K/\RR_{>0}$ is
compact, and fix also a point $s_0\in \PP_Y$. We have to show that
there exists a coordinate system $\{(t_1,\dots, t_{2n-1})\}$ in a
neighborhood $\ct$ of $s_0$ such that for any smooth $(n-1)$-form
$\omega$ with compact support contained in $\ct$, and for any $N\in
\NN$ one has
\begin{eqnarray}\label{sept11}
\lim_{\eps\to +0}\sup_{\xi\in
K}\left(\big|\int_{N(A_\eps)}e^{i\sum_{\alp=1}^{2n-1}t_\alp\xi_\alp}\omega(t)-
\int_{N(A_0)}e^{i\sum_{\alp=1}^{2n-1}t_\alp\xi_\alp}\omega(t)\big|\cdot
(1+|\xi|^N)\right)= 0
\end{eqnarray}
where the integrations are with respect to $t$.

Since $N(A_\eps)$ converges to $N(A_0)$ in the Hausdorff metric, we
may and will assume that $s_0\in N(A_0)$. Next since we may restrict
our considerations to a neighborhood of $Y\times_X\PP_X$ in $\PP_Y$,
we may and will assume in addition that $s_0\in Y\times_X\PP_X$.
Thus we assume that $s_0\in (Y\times_X\PP_X)\cap N(A_0)$.

Then we can write $s_0=(y_0,[n_0])$ where $y_0\in \pt A_0,\, n_0\in
X^*\times \{0_{Z^*}\}\subset Y^*$. Let $\tilde y\in \pt A$ be the
only point on the boundary of $A$ such that the outer normal at
$\tilde y$ is equal to $n_0$.

Let us introduce coordinates $\{(x_1,\dots, x_m)\}$ on $X$ such that
$(0,\dots,0)$ corresponds to $0_X\in X$, and the axis $x_m$ is
parallel to $n_0$, but has opposite orientation. Let us also choose
coordinates $\{(z_1,\dots,z_{n-m})\}$ on $Z$ such that $(0,\dots,0)$
corresponds to $0_Z$.

Then in a small neighborhood $\cv$ of $\tilde y$ there exists a
smooth strictly convex function
$F(x_1,\dots,x_{m-1};z_1,\dots,z_{n-m})$ such that in this
neighborhood $\cv$ the set $A$ coincides with the set $\{x_m\geq
F(x_1,\dots,x_{m-1};z_1,\dots,z_{n-m})\}$. Let us denote for brevity
\begin{eqnarray*}
\bar x:=(x_1,\dots,x_{m-1}),\, z=(z_1,\dots,z_{n-m}),\\
\nabla_{\bar x}F:=\left(\frac{\pt F}{\pt x_1},\dots, \frac{\pt
F}{\pt x_{m-1}}\right),\\
\nabla_zF:=\left(\frac{\pt F}{\pt z_1},\dots, \frac{\pt F}{\pt
z_{n-m}}\right).
\end{eqnarray*}
Then, in $\cv$, the normal cycle $N(A)$ is given by
\begin{eqnarray*}
N(A)=\{(\bar x,F(\bar x,z),z;[\nabla_{\bar x}F(\bar
x,z):-1:\nabla_zF(\bar x,z)])\}.
\end{eqnarray*}
Then obviously
\begin{eqnarray*}
N(A_\eps)=\{(\bar x,F(\bar x,z),\eps z;[\nabla_{\bar x}F(\bar
x,z):-1:\eps^{-1}\nabla_zF(\bar x,z)])\}.
\end{eqnarray*}
Let us consider the change of coordinates $(\bar x,z)\mapsto (\bar
x,\nabla_zF)$. It is locally invertible. Indeed the Jacobian is
equal to $\left[\begin{array}{cc}
             Id_{m-1}&\ast\\
                  0&\left(\frac{\pt^2F}{\pt z_i\pt z_j}\right)
          \end{array}\right]$ and the matrix $\left(\frac{\pt^2F}{\pt z_i\pt
          z_j}\right)$ is invertible since $F$ is strictly convex.

Let us denote $w:=\nabla_z F$, and rewrite the parametrization of
$N(A_\eps)$ in coordinates $(\bar x,w)$. Denote
\begin{eqnarray*}
G(\bar x,w):=F(\bar x,z),\\
H(\bar x,w):=\nabla_{\bar x}F(\bar x,z),\\
K(\bar x,w):=z.
\end{eqnarray*}

Then $$N(A_\eps)=\{\left(\bar x,G(\bar x,w),\eps K(\bar x,w);[H(\bar
x,w):-1:\eps^{-1}w]\right)\}.$$ In our fixed small neighborhood
$\ct$ of $s_0\in \PP_Y$ we can reparameterize $N(A_\eps)$ by
choosing $\theta:=\eps^{-1}w$. Then
$$N(A_\eps)=\{(\bar x,G(\bar x,\eps \theta), \eps K(\bar
x,\eps\theta);[H(\bar x,\eps\theta):-1:\theta])\}.$$ To abbreviate
the notation, now we will write points given in homogeneous
coordinates as $[a:-1:b]$, in linear coordinates as $(a,b)$. Then in
appropriate coordinate chart on $\PP_Y$ the normal cycle $N(A_\eps)$
is locally parameterized by $(\bar x,\theta)$ as
$$N(A_\eps)=\{(\bar x,G(\bar x,\eps\theta),\eps K(\bar x,\eps
\theta);H(\bar x,\eps\theta),\theta)\}.$$ In the same coordinate
chart the normal cycle of the limiting set $A_0$ is parameterized by
$(\bar x,\theta)$ as
$$N(A_0)=\{(\bar x,G(\bar x,0),0;H(\bar x,0),\theta)\}.$$
Applying an appropriate diffeomorphism of the ambient space $\PP_Y$,
we may get the following parametrization of the normal cycles:
\begin{eqnarray*}
N(A_\eps)=\{(\bar x,\bar G(\bar x,\theta, \eps),\bar K(\bar
x,\theta,\eps);\bar H(\bar x,\theta,\eps),\theta)\},\\
N(A_0)=\{(\bar x,0,0;0,\theta)\},
\end{eqnarray*}
where $\bar G,\bar K,\bar H$ are some $C^\infty$-smooth functions of
$\bar x,\theta,$ and $\eps\in [0,1]$, such that
$$\bar G(\bar x,\theta,\eps=0)=\bar H(\bar x,\theta,\eps=0)=\bar K(\bar
x,\theta,\eps=0)=0.$$ To simplify the notation we will denote
$$(\bar G,\bar H,\bar K)=:M,\, (\bar x,\theta)=:v.$$
Notice that $M$ takes values in $\RR^{n}$, $v\in \RR^{n-1}$. Clearly
$M$ is $C^\infty$-smooth in $(v,\eps)$, and $M(v,0)=0$.

In this notation we have
\begin{eqnarray}\label{sept12}
N(A_0)=\{(v,0_{\RR^n})\},\\\label{sept13}
N(A_\eps)=\{(v,M(v,\eps))\}.
\end{eqnarray}
From this description it is clear that the conic neighborhood $\cu$
of $T^*_{Y\times_X\PP_X}\PP_Y\backslash\underline{0}$ can be chosen
so small that $T^*_{N(A_\eps)}\PP_Y\backslash\underline{0}\subset
\Gamma$ for all $\eps\in [0,1]$ (recall that $\Gamma$ is the
complement of $\cu$).

Let us fix an arbitrary conic closed subset $T\subset
(\RR^{n-1}\times \RR^n)^*\backslash \{0\}$ such that $T/\RR_{>0}$ is
compact and
\begin{eqnarray}\label{sept14}
T\cap (\{0_{(\RR^{n-1})^*}\}\times \RR^{n*})=\emptyset.
\end{eqnarray}

In order to finish the proof of the proposition, i.e. of
(\ref{sept11}), we have to show that for any smooth compactly
supported $(n-1)$-form $\ome$ on $\RR^{2n-1}$ and any $N\in \NN$ one
has
\begin{eqnarray}\label{sept15}
\lim_{\eps\to +0}\sup_{\xi\in
T}\left(\big|\int_{N(A_\eps)}e^{i<\xi,y>}\omega(y)-
\int_{N(A_0)}e^{i<\xi,y>}\omega(y)\big|\cdot(1+|\xi|^N)\right)=0.
\end{eqnarray}

Let us show that
\begin{eqnarray}\label{sept17}
\overline{\lim_{\eps\to +0}}\sup_{\xi\in
T}\left(\big|\int_{N(A_\eps)}e^{i<\xi,y>}\omega(y)\big|\cdot
(1+|\xi|^N)\right)<\infty.
\end{eqnarray}

We will write $\xi=(\eta,\zeta)$ with $\eta\in
(\RR^{n-1})^*,\zeta\in (\RR^n)^*$. Let us denote also by
$\omega_\eps$ the pull-back of $\ome$ under the map $v\mapsto
(v,M(v,\eps))$. Thus $\ome_\eps$ is an $(n-1)$-form on $\RR^{n-1}$.
In this notation (\ref{sept17}) is rewritten
\begin{eqnarray}\label{sept18}
\overline{\lim_{\eps\to +0}}sup_{\xi\in
T}\left(\big|\int_{\RR^{n-1}}e^{i(<\eta,v>+<\zeta,M(v,\eps)>)}\ome_\eps(v)\big|\cdot
(1+(|\eta|+|\zeta|)^N)\right)<\infty.
\end{eqnarray}

The assumption (\ref{sept14}) implies that there exists a constant
$C$ such that for any $\xi=(\eta,\zeta)\in T$ one has
$$|\zeta|\leq C|\eta|.$$
Subdividing $T$ into finitely many parts one may assume that there
exists $i_0\in \{1,\dots,n-1\}$ such that for any $(\eta,\zeta)\in
T$ with $\eta=(\eta_1,\dots,\eta_{n-1})$ one has
$$|\eta_j|\leq |\eta_{i_0}| \mbox{ for any } j=1,\dots,n-1.$$
For simplicity of the notation we will assume that $i_0=1$, i.e.
$|\eta_j|\leq |\eta_1|$.

For $\xi=(\eta,\zeta)\in T$ let us consider the smooth map
$$F_{\eta,\zeta,\eps}\colon \RR^{n-1}\to \RR^{n-1}$$ given by
$F_{\eta,\zeta,\eps}(v_1,
v_2,\dots,v_{n-1}):=(v_1+<\frac{\zeta}{\eta_1},M(v,\eps)>,v_2,\dots,v_{n-1})$.

Let us fix a large ball $B\subset \RR^{n-1}$ containing the supports
of all forms $\ome_\eps$, $\eps\in [0,1]$. Since for any
$(\eta,\zeta)\in T$ one has $|\zeta|\leq C'|\eta_1|$ for some
constant $C'$, and since $M(v,0)=0$, it follows that there exists
$\eps_0>0$ such that for any $\eps\in [0,\eps_0]$ and any
$(\eta,\zeta)\in T$ the map $F_{\eta,\zeta,\eps}$ is a
diffeomorphism of the ball $B$ onto its image. Moreover any
$C^k$-norm of it and its inverse can be estimated by a constant
independent of $(\eta,\zeta)\in T$, $\eps\in [0,\eps_0]$.

The expression in the exponent in (\ref{sept18}) is
\begin{eqnarray*}
<\eta,v>+<\zeta,M(v,\eps)>=\sum_{l=1}^{n-1}\eta_lv_l+<\zeta,M(v,\eps)>=\\
\eta_1\cdot(v_1+<\frac{\zeta}{\eta_1}M(v,\eps),v>)+\sum_{l=2}^{n-1}\eta_lv_l=<\eta,F_{\eta,\zeta,\eps}(v)>.
\end{eqnarray*}

Hence after the change of variables using $F_{\eta,\zeta,\eps}$, the
inequality (\ref{sept18}) becomes equivalent to
\begin{eqnarray}\label{sept19}
\overline{\lim_{\eps\to +0}}\left(\sup_{(\eta,\zeta)\in
T}\big|\int_{\RR^{n-1}}e^{i<\eta,v>}\cdot
(F_{\eta,\zeta,\eps}^{-1})^*\ome_\eps\big|\cdot
(1+|\eta|^N)\right)<\infty.
\end{eqnarray}

But for $(\eta,\zeta)\in T,\eps\in [0,\eps_0]$ the forms
$(F_{\eta,\zeta,\eps}^{-1})^*\ome_\eps$ have uniformly bounded
supports, and all their $C^k$-norms are uniformly bounded. This
readily implies (\ref{sept19}). Thus (\ref{sept17}) is proved.

Similarly, but simpler, one proves that there exists a constant
$C_N$ such that
\begin{eqnarray}\label{sept16}
\sup_{\xi\in T}\big|\int_{N(A_0)}e^{i<\xi,y>}\omega(y)\big|\cdot
(1+|\xi|^N)<C_N.
\end{eqnarray}

\hfill

It remains to show that (\ref{sept16}), (\ref{sept17}) imply
(\ref{sept15}). To prove it, let us fix $\kappa>0$. There exist
$\eps_1>0$ and $C>0$ such that
\begin{eqnarray*}
\sup_{\eps\in [0,\eps_1]}\sup_{\xi\in T}
\big|\int_{N(A_\eps)}e^{i<\xi,y>}\ome(y)\big|\cdot
(1+|\xi|^{N+1})\leq C.
\end{eqnarray*}
Then for $\xi\in T$ such that $|\xi|>\frac{2C}{\kappa}$ we have for
any $\eps\in (0,\eps_1)$
\begin{eqnarray*}
\big|\int_{N(A_\eps)}e^{i<\xi,y>}\ome(y)-\int_{N(A_0)}e^{i<\xi,y>}\ome(y)\big|\cdot(1+|\xi|^N)\leq
\frac{2C}{|\xi|}<\kappa.
\end{eqnarray*}
On the other hand for $|\xi|\leq \frac{2C}{\kappa}$ we have
\begin{eqnarray*}
\sup_{|\xi|\leq
\frac{2C}{\kappa}}\big|\int_{N(A_\eps)}e^{i<\xi,y>}\ome(y)-
\int_{N(A_0)}e^{i<\xi,y>}\ome(y)\big|=\\
\sup_{|\xi|\leq \frac{2C}{\kappa}}\big| \int_{\RR^{n-1}}
e^{i<\eta,v>}\left(e^{i<\zeta,M(v,\eps)>}\ome_\eps(v)-\ome(v)\right)\big|
\mbox{ as } \eps \to +0
\end{eqnarray*}
which follows from the facts that the expression under the integral
tends to 0 uniformly in $\xi$ when $|\xi|\leq \frac{2C}{\kappa}$ .
This implies (\ref{sept19}). Hence Proposition \ref{PP:5} is proved.
\qed

\section{Valuations and integral geometry.}\label{S:IG}

\subsection{General Radon transform on
valuations.}\label{Ss:radon-val} Let $X_1\overset{q_1}{\leftarrow}Z
\overset{q_2}{\to} X_2$ be a double fibration. Then we have the map
$i_Z:=q_1\times q_2$,
$$i_Z\colon Z\to X_1\times X_2$$
which is a closed imbedding by the definition of a double fibration.
The goal of this section is to define a map
$$\phi\mapsto q_{2*}(\gamma\cdot q_1^*\phi)$$
where $\gamma$ is a fixed smooth valuation on $Z$ ($\gamma$ is
called kernel), and $\phi$ is a valuation on $X_1$, either smooth or
generalized. We will define such a map under appropriate assumptions
on the classes of valuations considered (Theorem
\ref{T:radon-val-1}) or also on the double fibration (Theorem
\ref{T:radon-val-1.5}).

\hfill

We will denote $X:=X_1\times X_2$.

\begin{theorem}\label{T:radon-val-1}
Let $X_1\overset{q_1}{\leftarrow}Z \overset{q_2}{\to} X_2$ be a
double fibration such that $q_2$ is proper. Let $\gamma\in
V^{\infty}(Z)$. Then the map $\phi\mapsto q_{2*}(\gamma\cdot
q_1^*\phi)$ is a well defined linear continuous map
$$R_\gamma\colon V^\infty(X_1)\to V^{-\infty}(X_2).$$
\end{theorem}


First we will need a general lemma where we use the notation
$\pi_{X*}\Gamma$ introduced in Section \ref{Ss:wave fronts} after
Remark \ref{R:push-forward-dens}.
\begin{lemma}\label{L:oct010}
Let $X$ be a smooth manifold. Let $\Gamma\subset T^*\PP_X\backslash
\underline{0}$ be closed conic subset such that
$\pi_{X*}(\Gamma)=\emptyset$. Then the product
$$ \vmi_{\emptyset,\Gamma}(X)\times V^\infty(X)\to \vmi_{\emptyset,\Gamma}(X)$$
is a bilinear jointly sequentially continuous map.
\end{lemma}
\def\vlg{\vmi_{\emptyset,\Gamma}}
{\bf Proof.} In the proof we may and will assume that $X$ is a
vector space in order to simplify the notation. Let $\phi\in
\vlg(X),\, \psi\in V^\infty(X)$ be represented by pairs of currents
$(C_\phi,T_\phi)$, $(C_\psi,T_\psi)$ respectively. Thus
\begin{eqnarray*}
WF(C_\phi)=\emptyset,\, WF(T_\phi)\subset \Gamma,\\
WF(C_\psi)=\emptyset,\, WF(T_\psi)= \emptyset.
\end{eqnarray*}
By (\ref{E:may0001})-(\ref{E:may0002}) the product $\phi\cdot \psi$
is given by the pair $(C,T)$ where
\begin{eqnarray*}
C=C_\phi\cap C_\psi,\\
T=(-1)^n \bar p_*\bar\Phi^*(q_1^*T_\phi\cap
q_2^*T_\psi)+\pi_X^*C_\phi\cap T_\psi+T_\phi\cap \pi_X^*C_\psi
\end{eqnarray*}
where the maps $$\PP_X\overset{\bar
p}{\leftarrow}\bar\PP\overset{\bar\Phi}{\to}\PP_X\times_X\PP_X$$
defined in (\ref{E:bar-phi}) and (\ref{E:bar-p}) of Section
\ref{Ss:valuations-mfld}, and $q_i\colon\PP_X\times_X\PP_X\to
\PP_X$, $i=1,2,$ be the natural projections onto the first and
second factors respectively. Recall that $\bar\PP$ is the oriented
blow up of the space $\PP_{X\times X}\cap \pi^{-1}_{X\times
X}(\Delta)$ (where $\Delta \subset X\times X$ is the diagonal) along
the submanifold $\cm_0\bigsqcup \cm_1\bigsqcup \cm_2$ where
\begin{eqnarray*}
\cm_0:=\{(x,[\xi:-\xi])\},\\
\cm_1:=\{(x,[\xi:0])\},\\
\cm_2:=\{(x,[0:\xi])\}.
\end{eqnarray*}
We denote by $\cn_i,\, i=0,1,2,$ the preimages in $\bar \PP$ of the
sets $\cm_i$.

It is easy to see that
$$WF(q_1^*T_\phi\cap q_2^*T_\psi)\subset q_1^*(\Gamma).$$
Then in order to prove the lemma it suffices to show that for any
$\phi\in \vlg(X),\, \psi\in V^\infty(X)$ one has
\begin{eqnarray}\label{cnd-a}
WF(C_\phi\cap C_\psi)= \emptyset,\\\label{cnd-b}
WF(T_\phi\cap\pi_X^*C_\psi)\subset \Gamma,\\\label{cnd-c}
WF(\pi_X^*C_\phi\cap T_\psi)\subset \Gamma,\\\label{cnd-d} \bar
p_*\bar \Phi^*(q_1^*\Gamma)\subset \Gamma
\end{eqnarray}
where the notations of pull-back and push-forward of {\itshape sets}
were introduced in Section \ref{Ss:wave fronts} in Proposition
\ref{P:lifting} and after Remark \ref{R:push-forward-dens}
respectively.

The inclusions (\ref{cnd-a})-(\ref{cnd-c}) are obvious. It remains
to check (\ref{cnd-d}) which can be rewritten as
\begin{eqnarray}\label{cndd1}
\bar p_*(q_1\circ \bar \Phi)^*(\Gamma)\subset \Gamma.
\end{eqnarray}
Let us fix points
$$X_0=(x_0,[\xi_0])\in \PP_X, X_1=(x_1,[\xi_1])\in \PP_X,
\alp\in \bar \PP$$ such that
$$\bar p( \alp)=X_0,\, (q_1\circ \bar\Phi)(\alp)=X_1.$$
then necessarily $x_0=x_1$. We have to show that if $v\in
T^*_{X_0}\PP_X\backslash \{0\}$ and $w\in \Gamma|_{X_1}$ satisfy
\begin{eqnarray}\label{oct3.1.6}
d\bar p^*(v)=d(q_1\circ \bar \Phi)^*(w)
\end{eqnarray}
then $v\in \Gamma$ (here the differentials are taken at the point
$\alp\in \bar \PP$).

We will consider 4 cases:
\newline
\underline{Case 1:} $\alp\not\in \cn_0\bigsqcup\cn_1\bigsqcup
\cn_2$;
\newline
\underline{Case 2:} $\alp\in \cn_0$;
\newline
\underline{Case 3:} $\alp\in \cn_1$;
\newline
\underline{Case 4:} $\alp\in \cn_2$.

\hfill

Before we will start considering these cases let us introduce more
notation. Since we assume that $X$ is a vector space, we can write
$\PP_X=X\times \PP_+(X^*)$, and $\bar \PP=X\times Y$ where $Y$ is an
appropriate blow up space of $\PP_+(X^*\oplus X^*)$. The maps $\bar
p$ and $q_1\circ \bar \Phi$ preserve these decompositions, more
precisely $\bar p=(Id_X,\tilde p)$ where $\tilde p\colon Y\to
\PP_+(X^*)$, and $q_1\circ \bar\Phi=(Id_X,\tilde \Phi)$ where
$\tilde \Phi\colon Y\to \PP_+(X^*)$. It is easy to see that $\tilde
p$ and $\tilde \Phi$ are submersions as well as their restrictions
to the boundary of $Y$. Let us write the vector $v\in
T^*_{X_1}\PP_X$ in the form $v=(v',v'')$ where $v'\in T^*_{x_0}X$,
$v''\in T^*_{[\xi_0]}\PP_+(X^*)$. Let us show that if
(\ref{oct3.1.6}) is satisfied then
\begin{eqnarray}\label{october3.1.7}
v''\ne 0.
\end{eqnarray}
Indeed otherwise $d\bar p^*(v)=(v',0)$. Hence $d(q_1\circ \bar
\Phi)^*(w)=(v',0)$. But since $\tilde \Phi$ is a submersion $w$ has
the form $w=(v',0)$. This implies that $v'\in
\pi_{X*}(\Gamma)=\emptyset$. This is a contradiction.

\hfill

Now let us consider Case 1. Thus $\alp\not\in
\cn_0\bigsqcup\cn_1\bigsqcup \cn_2$. We may assume that $\alp\in
(\PP_X\times_X \PP_X)\backslash(\cm_0\bigsqcup\cm_1\bigsqcup
\cm_2)$. Then $\alp$ must have the form
$$\alp=(x_0,[\xi_1:\eta]) \mbox{ with } [\eta]\ne [-\xi_1]\mbox{ and }
[\xi_1+\eta]=[\xi_0].$$ Let us consider the subspace $S\subset
T_\alp\bar \PP$ defined as the tangent space at $\alp$ to the
submanifold $\{(x,[\xi_1:\zeta])|\,x\in X,\zeta\in X^*\}$ (here
$\xi_1$ is fixed). Let us denote also by $S'$ the tangent space to
$\{[\xi_1:\zeta]|\, \zeta\in X^*\}$ at $[\xi_1:\eta]$. In this
notation we have the following decompositions:
\begin{eqnarray}\label{decom1}
T_{X_0}\PP_X=X\oplus T_{[\xi_0]}(\PP_+(X^*)),\\\label{decom2}
T_{X_1}\PP_X=X\oplus T_{[\xi_1]}(\PP_+(X^*)),\\\label{decom3}
S=X\oplus S'.
\end{eqnarray}
It is clear that the restriction of the differential $d\bar
p|_S\colon S\to T_{X_0}\PP_X$ is an epimorphism, and moreover $d\bar
p|_S$ maps the decomposition (\ref{decom3}) to the decomposition
(\ref{decom1}).

Next the restriction of the differential $d(q_1\circ \bar
\Phi)|_S\colon S\to T_{X_1}\PP_X$ maps decomposition (\ref{decom3})
to (\ref{decom2}), and its kernel is equal to $S'$. Hence for $w\in
\Gamma|_{X_1}$ one has $d(q_1\circ \bar \Phi)^*(w)|_S\in X^*\times
\{0_{S^{'*}}\}$. This implies that if $d\bar p^*(v)=d(q_1\circ \bar
\Phi)^*(w)$ then
\begin{eqnarray}\label{oct3.1.10}
v=(v',0)\in X^*\times \{0\}\subset T_{X_0}^*\PP_X
\end{eqnarray}
(here we have used that $d\bar p|_S\colon S\to T_{X_0}\PP_X$ is an
epimorphism). This contradicts to (\ref{october3.1.7}), thus Case 1
is proved.


\hfill

Let us consider Case 2. Thus we assume that $\alp\in \cn_0$. It is
easy to see that there exists a subspace $S_0\subset T_\alp\bar \PP$
of the form $S_0=X\oplus S_0'$ with $S_0'\subset TY$, such that
$d\tilde p|_{S_0'}\colon S_0'\to T_{x_0}\PP_+(X^*)$ is an
isomorphism, and $d\tilde \Phi|_{S_0'}\colon S_0'\to
T_{x_1}(\PP_+(X^*))=T_{x_0}(\PP_+(X^*))$ vanishes. The second fact
implies that $d(q_1\circ \bar \Phi)^*|_{S'}(w)$ has the form
$(w',0_{S^{'*}})$ with $w'\in X^*$. This and the first fact imply
that $v=(w',0)$. This is again a contradiction to
(\ref{october3.1.7}), and Case 2 is proved.

\hfill

Let us consider Case 3. Namely assume that $\alp\in \cn_1$. This
implies that $[\xi_0]=[\xi_1]$. Hence $X_0=X_1$. It is easy to see
that there exists a linear subspace $S_1\subset T_\alp\bar \PP$ of
the form $S_1=X\oplus S_1'$ with $S_1'\subset TY$ such that $d\tilde
p|_{S_1'}\colon S_1'\to T_{[\xi_0]}\PP_+(X^*)$ and $d\tilde
\Phi|_{S_1'}\colon S_1'\to T_{[\xi_0]}\PP_+(X^*)$ are isomorphisms,
and $$d\tilde\Phi|_{S_1'}\circ(d\tilde p|_{S_1'})^{-1}\colon
T_{[\xi_0]}\PP_+(X^*)\to T_{[\xi_0]}\PP_+(X^*)$$ is the identity
map. This immediately  implies that if $d\bar p^*(v)=d(q_1\circ \bar
\Phi)^*(w)$ then $v=w$. Hence $v\in \Gamma$.

\hfill

Let us consider Case 4, i.e. $\alp\in \cn_2$. It is easy to see that
there exists a subspace $S_2\subset T_\alp\bar \PP$ of the form
$S_2=X\oplus S_2'$ with $S_2'\subset TY$ such that $d\tilde
p|_{S_2'}\colon S_2'\to T_{x_0}(\PP_+(X^*))$ vanishes, and $d\tilde
\Phi|_{S_2'}\colon S_2'\to T_{x_0}(\PP_+(X^*))$ is an isomorphism.
These two facts and the assumption $d\bar p^*(v)=d(q_1\circ \bar
\Phi)^*(w)$ imply that $w=(w',0)$. Hence $w'\in
\pi_{X*}(\Gamma)=\emptyset$ which is a contradiction.

Lemma is proved. \qed


\hfill

{\bf Proof of Theorem \ref{T:radon-val-1}.} For $i=1,2$ let us
denote for brevity
$$\Gamma_i:=T^*_{Z\times_{X_i}\PP_{X_i}}(Z)\backslash \underline{0}.$$
Obviously $\pi_{X_{i}\ast}(\Gamma_i)=\emptyset$. The map of
pull-back $$q_i^*\colon V^\infty(X_i)\to
\vmi_{\emptyset,\Gamma_i}(Z)$$ is a continuous linear operator (since $V^\infty(X_1)$ is a Fr\'echet space, sequential continuity implies
topological continuity). By
this and Lemma \ref{L:oct010} the map $[\phi\mapsto \gamma\cdot
q_1^*\phi]$ is a continuous linear operator $V^\infty(X_1)\to
\vmi_{\emptyset,\Gamma_1}(Z)$. Thus in order to finish the proof of
the theorem it remains to show that the push-forward map $q_{2*}\colon V^{-\infty}{\emptyset,\Gamma_1}(Z)\to V^{-\infty}(X_2)$
is well defined in the sense of Section \ref{Ss:pushforward-gen-submersions}; then it is automatically sequentially continuous.
  For this it
suffices to check that the pairs of sets $(\emptyset,\Gamma_1)$ and
$(\emptyset,\Gamma_2)$ satisfy the conditions
(\ref{pro-a})-(\ref{pro-e}). The conditions
(\ref{pro-a})-(\ref{pro-d}) are satisfied trivially in this case.

Let us check (\ref{pro-e}). It is easy to see that this condition is
equivalent the following general claim.

\begin{claim}\label{CLAIM:double-fibr}
If $X_1\overset{q_1}{\leftarrow}Z \overset{q_2}{\to} X_2$ is a
double fibration, then the sets $Z\times_{X_1}\PP_{X_1}$ and
$Z\times_{X_2}\PP_{X_2}$ intersect transversally in $\PP_Z$.
\end{claim}
Let us prove the claim. It is equivalent that the sets
$M_1:=Z\times_{X_1}(T^*X_1\backslash \underline{0})$ and
$M_2:=Z\times_{X_2}(T^*X_2\backslash \underline{0})$ intersect
transversally in $T^*Z$.

Let us fix a point $A\in M_1\cap M_2$. We have to show that
$T_AM_1+T_AM_2=T_A(T^*Z)$. Let us denote by $F$ the only fiber of
the natural projection $T^*Z\to Z$ which contains $A$. Since
$q_1,q_2$ are submersions we have
\begin{eqnarray}\label{octob3.1.1}
T_AF+T_AM_i=T_A(T^*Z),\, i=1,2.
\end{eqnarray}

On the other hand since $q_1\times q_2\colon Z\to X_1\times X_2$ is
an imbedding
\begin{eqnarray}\label{octob3.1.2}
T_AF\subset T_AM_1+T_AM_2.
\end{eqnarray}
From (\ref{octob3.1.1}) and (\ref{octob3.1.2}) we conclude that
$T_AM_1+T_AM_2=T_A(T^*Z)$. Hence $M_1,M_2$ are transversal and Claim
\ref{CLAIM:double-fibr} is proved.

Hence  the condition (\ref{pro-e}) is proved, and Theorem
\ref{T:radon-val-1} as well. \qed

\hfill


In order to formulate the second main result let us introduce more
notation. For $j=1,2$ let
$$p_j\colon Z\times_{X_j}\PP_{X_j}\to \PP_{X_j},\, i_j\colon
Z\times_{X_j}\PP_{X_j}\inj \PP_Z$$ be the obvious maps.

\begin{theorem}\label{T:radon-val-1.5}
Let $X_1\overset{q_1}{\leftarrow}Z \overset{q_2}{\to} X_2$ be a
double fibration. Assume in addition that
$T^*_{Z\times_{X_1}\PP_{X_1}}\PP_Z$ and
$i_{2*}p_2^*(T^*\PP_{X_2}\backslash\underline{0})$ are disjoint
subsets of $T^*\PP_Z$. Let $\gamma\in V^\infty(Z)$.


(1) If $q_1$ is proper then the map $V^{\infty}(X_2)\to
V^{-\infty}(X_1)$ defined in Theorem \ref{T:radon-val-1} by
$[\psi\mapsto q_{1*}(\gamma\cdot q_2^*\psi)]$ extends (uniquely) to
a continuous linear operator
$$V^{-\infty}(X_2)\to V^{-\infty}(X_1)$$
when both spaces are equipped with the weak topologies.

(2) If $q_2$ is proper then the map $[\phi\mapsto q_{2*}(\gamma\cdot
q_1^*\phi)]$ is a well defined continuous linear operator
$$V^\infty(X_1)\to V^{\infty}(X_2)$$
when both spaces are equipped with their usual Fr\'echet topologies.
\end{theorem}

\begin{remark}
Below in Proposition \ref{P:sufficient-cond} we will show that the
assumptions of Theorem \ref{T:radon-val-1.5} are satisfied provided
that the natural projection $T^*_Z(X_1\times X_2)\backslash\{0\}\to
T^*X_2\backslash\{0\}$ is a submersion.

\end{remark}


\hfill

{\bf Proof of Theorem \ref{T:radon-val-1.5}.} Let us start with part
(2). By Theorem \ref{T:radon-val-1} we have a separately continuous
bilinear map $V^\infty(X_1)\times V^\infty_c(X_2)\to \CC$ given by
$$(\phi,\psi)\mapsto \int_Z(q_1^*\phi\cdot \gamma)\cdot q_2^*\psi.$$
First let us show that the last map extends (uniquely, of course) to
a separately continuous bilinear map
$$V^\infty(X_1)\times \vmi_c(X_2)\to \CC.$$
Let us denote for brevity
$$\Lam:=q_2^*(T^*X_2\backslash\underline{0}),\,
\Gamma:=i_{2*}p_2^*(T^*\PP_{X_2}\backslash \underline{0}).$$ Also as
above in this section we denote
$\Gamma_j:=(T^*_{Z\times_{X_j}\PP_{X_j}}\PP_Z)\backslash\underline{0}$.
In this notation for any $\psi\in \vmi(X_2)$ we have $q_2^*\psi\in
\vmi_{\Lam,\Gamma}(Z)$. Also it is clear that
$\pi_{Z*}(\Gamma_j)=0$. Hence by Lemma \ref{L:oct010} the map
$V^\infty(X_1)\to V^{-\infty}_{\emptyset,\Gamma_1}(Z)$ given by
$\phi\mapsto \gamma\cdot q_1^*\phi$ is continuous.

Thus it suffices to show that the pairs of sets $(\Lam,\Gamma)$ and
$(\emptyset, \Gamma_1)$ satisfy the conditions
(\ref{pro-a})-(\ref{pro-e}). Explicitly in this case they mean:
\begin{eqnarray}\label{pro-na}
\Lam\cap s(\emptyset)=\emptyset,\\\label{pro-nb} \Gamma\cap
s(\pi_Z^*\emptyset)=\emptyset,\\\label{pro-nc} \Gamma_1\cap
s(\pi_Z^*\Lam)=\emptyset,\\\label{pro-nd} \mbox{if }
(z,[\xi_1],u_1,0)\in \Gamma_1 \mbox{ and }(z,[\xi],u,0)\in \Gamma
\mbox{ then } u_1\ne -u,\\\nonumber \mbox{if } (z,[\xi])\in \PP_Z,\,
(u,\eta_1)\in \Gamma_1|_{(z,[\xi])}, \mbox{ and } (-u,\eta)\in
\Gamma|_{(z,[-\xi])}, \mbox{ then }\\\label{pro-ne}
d\theta^*(0,\eta_1,\eta)\ne (0,l,-l)\in
T^*_{(z,[\xi],[\xi])}(\PP_Z\times_Z\PP_Z)
\end{eqnarray}
where $\theta\colon \PP_Z\times_Z\PP_Z\to \PP_Z\times_Z\PP_Z$ is
defined by $\theta(z,[\zeta_1],[\zeta_2])=(z,[\zeta_1],[-\zeta_2])$.

The conditions (\ref{pro-na}), (\ref{pro-nb}) are satisfied
trivially. The conditions (\ref{pro-nc}), (\ref{pro-nd}) are
satisfied since $\Gamma_1$ contains no elements of the form
$(z,[\xi_1],u_1,0)$. The condition (\ref{pro-ne}) easily follows
from our assumption $\Gamma_1\cap \Gamma=\emptyset$.

Thus we got a jointly sequentially continuous bilinear map
$$B\colon V^\infty(X_1)\times V^{-\infty}_c(X_2)\to \CC.$$
Then for any $\phi\in V^{\infty}(X_1)$ the map $V_c^{-\infty}(X_2)\to \CC$ given by $\psi\mapsto B(\phi,\psi)$
is sequentially continuous. But the topology on $V^{-\infty}_c(X_2)=(V^\infty(X_2))^*$ coincides with the weak topology.
Since $V^\infty(X_2)$ separable Fr\'echet space, any sequentially continuous (in the weak topology) linear functional on
its dual $(V^\infty(X_2))^*$ is continuous in the weak topology (by Ch. IV, \S 6, Corollary 3 in \cite{schaefer}).
This implies that $B$ is separately continuous.

Then by duality $B$ induces a continuous linear map $V^\infty(X_1)\to
V^\infty(X_2)$ where the target space is equipped with the {\itshape
weak} topology. It remains to prove the continuity of this map when
the target space is equipped with the usual Fr\'echet topology.
Let us denote temporarily this map by $R$. Since $V^\infty(X_1)$ and
$V^\infty(X_2)$ are Fr\'echet spaces, by the closed graph theorem
(see e.g. \cite{robertson-robertson}, Ch. VI, \S 3, Thm. 8) it
suffices to show that $R$ has a closed graph. Assume that a sequence
$\phi_N\to \phi$ in $V^\infty(X_1)$, and $R(\phi_N)\to \psi$ in
$V^\infty(X_2)$. Since $R$ is continuous is the weak topology, it
follows that $\psi=R(\phi)$. Hence the graph of $R$ is closed in the
Fr\'echet topology. Part (2) is proved.

\hfill

Let us prove now part (1). By Theorem \ref{T:radon-val-1} we have a
separately continuous bilinear map
$$V^\infty_c(X_1)\times V^\infty(X_2)\to \CC$$
given by $(\phi,\psi)\mapsto \int_Zq_1^*\phi\cdot (\gamma\cdot
q_2^*\psi)$. Let us show that this map extends (uniquely) to a
separately continuous bilinear map
$$V^\infty_c(X_1)\times V^{-\infty}(X_2)\to \CC.$$ First observe that if
$\phi,\psi$ are smooth then $q_1^*\phi\in
\vmi_{\emptyset,\Gamma_1}(Z)$, $q_2^*\psi\in
\vmi_{\emptyset,\Gamma_2}(Z)$, and by Lemma \ref{L:oct010}
$$q_1^*\phi\cdot\gamma \in \vmi_{\emptyset,\Gamma_1}(Z),\, \gamma\cdot
q_2^*\psi\in \vmi_{\emptyset,\Gamma_2}(Z).$$ Hence for smooth
$\phi,\psi$ the products $(q_1^*\phi\cdot\gamma)\cdot q_2^*\psi$ and
$q_1^*\phi\cdot (\gamma\cdot q_2^*\psi)$ are defined and are equal
to each other by the associativity proven in \cite{alesker-bernig},
Theorem 8.3. In particular
\begin{eqnarray}\label{oct00}
\int_Zq_1^*\phi\cdot (\gamma\cdot
q_2^*\psi)=\int_Z(q_1^*\phi\cdot\gamma)\cdot q_2^*\psi.
\end{eqnarray}
Let us show that the right hand side of (\ref{oct00})
extends to a jointly sequentially continuous bilinear map
$V^\infty_c(X_1)\times V^{-\infty}(X_2)\to \CC$. To prove it, it
suffices to show that the pairs of sets $(\Lam,\Gamma)$ and
$(\emptyset,\Gamma_2)$ satisfy the conditions
(\ref{pro-a})-(\ref{pro-e}). But this has been proven in the proof
of part (2) of the theorem. Next, as in the proof of part (2), this bilinear functional
separately continuous since $V^{-\infty}(X_2)=(V_c^\infty(X_2))^*$ and any sequentially continuous
(in the weak topology) linear functional on it is continuous. \qed

\hfill

Now we would like to give a sufficient condition on a double
fibration implying the assumption of Theorem \ref{T:radon-val-1.5}.
Recall that $Z\subset X_1\times X_2$ is a closed submanifold.
Consider the conormal bundle $T^*_Z(X_1\times X_2)\subset
T^*X_1\times T^*X_2$. Since $q_i$'s are submersions, the projections
$T^*_Z(X_1\times X_2)\to T^*X_i$, $i=1,2,$ map $T^*_Z(X_1\times
X_2)\backslash \underline{0}$ to $ T^*X_i\backslash\underline{0}$.
We are going to prove
\begin{proposition}\label{P:sufficient-cond}
Assume that the double fibration $X_1\overset{q_1}{\leftarrow}Z
\overset{q_2}{\to} X_2$ satisfies that $$T^*_Z(X_1\times
X_2)\backslash \underline{0}\to T^*X_2\backslash\underline{0}$$ is a
submersion. Then the assumption of Theorem \ref{T:radon-val-1.5} is
satisfied, namely
$$\left(T^*_{Z\times_{X_1}\PP_{X_1}}\PP_Z\right)\cap
i_{2*}p_2^*(T^*\PP_{X_2}\backslash\underline{0})=\emptyset.$$
\end{proposition}
{\bf Proof.} The map $T^*_Z(X_1\times X_2)\backslash
\underline{0}\to T^*X_2\backslash\underline{0}$ is a submersion if
and only if the induced map
\begin{eqnarray}\label{map-proj}
\PP_+(T^*_Z(X_1\times X_2))\to \PP_{X_2}
\end{eqnarray}
is a submersion.

By Claim \ref{CLAIM:double-fibr} the submanifolds
$Z\times_{X_1}\PP_{X_1}$ and $Z\times_{X_2}\PP_{X_2}$ intersect
transversally in $\PP_Z$. Let us denote their intersection by $\cc$.
We claim that the restriction of the obvious projection $p_2\colon
Z\times_{X_2}\PP_{X_2}\to \PP_{X_2}$ to $\cc$ is a submersion. For
consider the map
$$J\colon \PP_+(T^*_Z(X_1\times X_2))\to Z\times_{X_2}\PP_{X_2}$$
given by $J(z,[\xi_1:\xi_2]):=(z,[\xi_2])$. It is easy to see that
the image of $J$ is equal to $\cc$, and $J$ is a diffeomorphism
$$J\colon\PP_+(T^*_Z(X_1\times X_2))\tilde\to \cc.$$
Moreover the composition
$$\PP_+(T^*_Z(X_1\times
X_2))\overset{J}{\to}Z\times_{X_2}\PP_{X_2}\overset{p_2}{\to}\PP_{X_2}$$
coincides with the map (\ref{map-proj}) which, by the assumption, is
a submersion. This implies that the restriction of $p_2$ to $\cc$
$$p_2|_{\cc}\colon \cc\to \PP_{X_2}$$
is a submersion as we wanted.

Let us fix now a point $A\in \cc$. Let us denote
\begin{eqnarray*}
E:=T_A\PP_Z,\\
L_i:=T_A(Z\times_{X_i}\PP_{X_i}),\, i=1,2.
\end{eqnarray*}
Since $Z\times_{X_1}\PP_{X_1}$ and $Z\times_{X_2}\PP_{X_2}$ are
transversal, one has
\begin{eqnarray}\label{oct:pluss}
L_1+L_2=E.
\end{eqnarray}

Obviously the fiber of $T^*_{Z\times_{X_1}\PP_{X_1}}\PP_Z$ over $A$
is
$$\left(T^*_{Z\times_{X_1}\PP_{X_1}}\PP_Z\right)|_A=L_1^\perp\subset
E^*.$$ It is easy to see that the fiber of
$i_{2*}p_2^*(T^*\PP_{X_2})$ over $A$ is equal to
$$\left( Ker(dp_2\colon L_2\to T_{p_2(A)}\PP_{X_2})\right)^\perp.$$
Thus to prove the proposition we have to check that
\begin{eqnarray}\label{oct3.1.22}
L_1+ Ker(dp_2\colon L_2\to T_{p_2(A)}\PP_{X_2})=E.
\end{eqnarray}
Since, as we have proven, $p_2|_\cc\colon \cc\to \PP_{X_2}$ is a
submersion, and obviously $T_A\cc=L_1\cap L_2$, the map
$$dp_2\colon L_1\cap L_2\to T_{p_2(A)}\PP_{X_2}$$ is onto. To prove
(\ref{oct3.1.22}) let us fix an arbitrary element $x\in E$. By
(\ref{oct:pluss} ) there exists a presentation $x=l_1+l_2$ with
$l_i\in L_i$. By the surjectivity of $dp_2\big|_{L_1\cap L_2}$,
there exists $\tau\in L_1\cap L_2$ such that
$$dp_2(l_2)=dp_2(\tau).$$
Set $k:=l_2-\tau\in Ker(dp_2|_{L_2})$. We have
$$x=(l_1+\tau)+k\in L_1+Ker(dp_2|_{L_2}).$$
Thus (\ref{oct3.1.22}) is proved. Hence the proposition is proved.
\qed

Theorem \ref{T:radon-val-1.5} and Proposition
\ref{P:sufficient-cond} imply immediately the following corollary
which we would like to state separately since it will be used in
this form in the applications.

\begin{corollary}\label{COR:radon-val-2}
Let $X_1\overset{q_1}{\leftarrow}Z \overset{q_2}{\to} X_2$ be a
double fibration. Assume that the projection $T^*_Z(X_1\times
X_2)\backslash \underline{0}\to T^*X_2\backslash\underline{0}$ is a
submersion.

(1) If $q_1$ is proper then the map $V^{\infty}(X_2)\to
V^{-\infty}(X_1)$ defined by Theorem \ref{T:radon-val-1} as
$[\psi\mapsto q_{1*}(\gamma\cdot q_2^*\psi)]$ extends (uniquely) to
a topologically continuous linear operator
$$V^{-\infty}(X_2)\to V^{-\infty}(X_1)$$
when both spaces are equipped with the weak topologies.

(2) If $q_2$ is proper then the map $[\phi\mapsto q_{2*}(\gamma\cdot
q_1^*\phi)]$ is a well defined topologically continuous linear operator
$$V^\infty(X_1)\to V^{\infty}(X_2)$$
when both spaces are equipped with their usual Fr\'echet topologies.
\end{corollary}


\begin{example}\label{Example:double-fibrations}
Let us give some examples of double fibrations satisfying the
assumptions of Corollary \ref{COR:radon-val-2}.

(1) Let $X=\RR^n$. Let $Y=\ca Gr_{k,n}$ be the Grassmannian of
affine $k$-dimensional subspaces in $\RR^n$, $0<k<n$. Let $Z$ be the
incidence variety, i.e.
$$Z:=\{(x,H)\in\RR^n\times\ca Gr_{k,n}|\, x\in H\}.$$
Then $X=\RR^n\overset{q_1}{\leftarrow}Z\overset{q_2}{\to}\ca
Gr_{k,n}=Y$ with the obvious projections $q_1,q_2$ is a double
fibration, and the map
\begin{eqnarray*}
T^*_Z(X\times Y)\backslash\underline{0}\to
T^*X\backslash\underline{0}
\end{eqnarray*}
is a submersion.

(2) The next example will be particularly important in Section
\ref{Ss:constructible}. Let $X=\RR\PP^n$ be the real projective
space, i.e. the space of real lines in $\RR^{n+1}$ passing through
0. Let $Y=\RR\PP^{n\vee}$ be the dual projective space, i.e. the
space of real hyperplanes in $\RR^{n+1}$ passing through 0. Let
$Z:=\{(l,H)\in \RR\PP^n\times\RR\PP^{n+1}|\, l\subset H\}$ be again
the incidence variety. Then
$X=\RR\PP^n\overset{q_1}{\leftarrow}Z\overset{q_2}{\to}\RR\PP^{n\vee}=Y$
is a double fibration, and both projections
\begin{eqnarray*}
T^*_Z(X\times Y)\backslash\underline{0}\to
T^*X\backslash\underline{0},\,\, T^*_Z(X\times
Y)\backslash\underline{0}\to T^*Y\backslash\underline{0}
\end{eqnarray*}
are submersions. We will prove in Appendix that this double
fibration satisfies the assumptions of Corollary
\ref{COR:radon-val-2}.

(3) The example (2) can be generalized as follows. Let $\KK$ be
either $\RR$, $\CC$, or $\HH$ (where $\HH$ denotes the
non-commutative field of quaternions). Fix natural numbers $k,m,n$
such that $0<k\ne m<n$. Denote by $\grk_{k,n}$ the Grassmannian of
$\KK$-linear subspaces in $\KK^n$ of $\KK$-dimension $k$. Let
$X=\grk_{k,n}$, $Y=\grk_{m,n}$. Let $Z$ be the incidence variety
again. Then $X\overset{q_1}{\leftarrow}Z\overset{q_2}{\to}Y$ is a
double fibration. If $k=1$ or $n-1$ then the projection
\begin{eqnarray*}
T^*_Z(X\times Y)\backslash\underline{0}\to
T^*X\backslash\underline{0}.
\end{eqnarray*}
is a submersion. We refer to Appendix for the proof of this fact.

(4) Let $X$ (resp. $Y$) be the Grassmannian of affine
$\KK$-subspaces in $\KK^n$ of $\KK$-dimension $k$ (resp. $m$).
Assume that $k=1$ or $n-1$, and $k\ne m, 1<m<n$. Let $Z$ be the
incidence variety. In this case all the properties of the
corresponding double fibration remain the same as in the previous
examples.
\end{example}

\subsection{Relations to the Gelfand style integral geometry.}\label{Ss:gelfand}
In this section we describe the explicit relation between the Radon
transform on valuations introduced in Section \ref{Ss:radon-val} and
the classical integral geometry in the spirit of Gelfand's school
\cite{gelfand-graev-vilenkin}.

In a rather general setting, the Gelfand style integral geometry can
be described as follows. For brevity we will denote by
$\cm^\infty(X)$ the space of smooth densities (i.e. smooth measures)
on a manifold $X$. Let
$$X_1\overset{q_1}{\leftarrow}Z\overset{q_2}{\to}X_2$$
be a double fibration. Assume that $q_2$ is proper. Let us fix a
smooth density $m$ on $Z$. Consider the following operator
$$R\colon C^\infty(X_1)\to \cm^\infty(X_2)$$
given by
\begin{eqnarray}\label{D:radon-funct}
R(f):=q_{2*}(m\cdot q_1^*f)
\end{eqnarray}
where $q_1^*$ is the
usual pull-back of functions, and $q_{2*}$ is the usual push-forward
on densities. Since $q_2$ is a submersion the image of $R$ lies in
smooth densities, and $R$ is a continuous linear operator.

\hfill

All classically known Radon transforms on functions have such a
form, e.g. the classical Radon transform over hyperplanes in $\RR^n$
or in projective space $\RR\PP^n$, and the Radon transform for a
pair of Grassmann manifolds \cite{gelfand-graev-rosu}. In classical
situations the manifold $Z$ is often certain incidence variety acted
by a Lie group in a transitive way, and $m$ is a Haar measure.

\hfill

In the above high generality we know only one non-trivial result
\cite{guillemin}. Let us state a consequence of it. Let us assume
that all the manifolds $X_1,X_2,Z$ are compact. Assume moreover that
the projection
$$T^*_Z(X_1\times X_2)\backslash\underline{0}\to T^*X_2
\backslash\underline{0}$$ is an injective immersion. Assume finally
that the density $m$ on $Z$ is strictly positive everywhere. Then
the Radon transform $R$ has a finite dimensional kernel.

\hfill

Now let us show that the general Radon transform on functions
(\ref{D:radon-funct}) can be considered as a special case of the
general Radon transform on valuations from Section
\ref{Ss:radon-val}.

First let us recall that the space of smooth densities
$\cm^\infty(Z)$ is a subspace of smooth valuations $V^\infty(Z)$.
Thus we may consider the density $m$ as a valuation. Using Theorem
\ref{T:radon-val-1} consider the operator
\begin{eqnarray*}
R'\colon V^\infty(X_1)\to V^{-\infty}(X_2)
\end{eqnarray*}
given by
\begin{eqnarray}\label{radon-new}
R'(\phi):=q_{2*}(m\cdot q_{1}^*\phi)
\end{eqnarray}
where now $q_1^*$ and $q_{2*}$ denote the pull-back and push-forward
on valuations.

Recall that $V^\infty(X_1)$ has a closed subspace
$W_1(V^\infty(X_1))$ such that the quotient space
$V^\infty(X_1)/W_1(V^\infty(X_1))$ is canonically isomorphic to
$C^\infty(X_1)$. Moreover by Proposition \ref{P:pull-back-filtr-gen}
$$q_1^*(W_1(V^\infty(X_1)))\subset W_1(V^{-\infty}(Z)).$$
Next by \cite{part4}, equality (7.3.1), $\cm^\infty(Z)\cdot
W_1(V^{-\infty}(Z))=0$. Hence the operator $R'$ admits a unique
factorization
$$\btriangle<1`1`1;800>[V^\infty(X_1)`C^\infty(X_1)`V^{-\infty}(X_2);`R'`R''].$$
It is easy to see that the operator $R''$ takes values in smooth
measures $\cm^\infty(Y)\subset V^{-\infty}(Y)$ and coincides with
the classical Radon transform $R$ on functions
(\ref{D:radon-funct}).

\subsection{Relations to the integral geometry of constructible
functions.}\label{Ss:constructible} Let us discuss first the
integral geometry of constructible functions in the general setting.
The notion of constructible function may have slightly different
meanings in various contexts; but the general ideas are the same,
though the technical details may be different. Here we will use the
space $\cf(X)$ of constructible functions on a real analytic
manifold $X$ described in Section \ref{Ss:detailed-overview}.







Let $f\colon X\to Y$ be a real analytic map of real analytic
manifolds. We have the obvious pull-back map
$$f^*\colon \cf(Y)\to \cf(X)$$
preserving the class of constructible functions. If the map $f$ is
proper then we have the push-forward map
$$f_*\colon \cf(X)\to \cf(Y)$$
which is uniquely characterized by the following property: for a
subanalytic subset $P\subset X$ one has
$$f_*(\One_{P})(y)=\chi(P\cap f^{-1}(y))$$
for any $y\in Y$. One can show that $f_*$ maps constructible
functions to constructible ones (see \cite{kashiwara-schapira}, \S
9.7). This push-forward map $f_*$ on constructible functions is also
called integration with respect the Euler characteristic along the
fibers of $f$. In particular if $Y$ is a single point $f_*\phi$ is
called the integral of $\phi$ with respect to the Euler
characteristic and is denoted by $\int_X \phi$.

\hfill

Let us return now to integral geometry. Let
$X_1\overset{q_1}{\leftarrow}Z\overset{q_2}{\to}X_2$ be a double
fibration when $X_1,X_2,Z$ are real analytic manifolds, and
$q_1,q_2$ are real analytic maps. Assume that $q_2$ is proper. Then
we have the operator on constructible functions $R\colon \cf(X_1)\to
\cf(X_2)$ given by
$$R(\phi):=q_{2*}q_1^*\phi$$
which is another version of the Radon transform.

As far as we know this type of transforms (e.g. the inversion
problem) in a somewhat different context of {\itshape complex
algebraic} constructible functions on the complex projective space
$\CC\PP^n$ was studied by Viro \cite{viro}. For {\itshape real}
projective space $\RR\PP^n$ and for functions constructible in a
more restrictive sense the problem was studies by Khovanskii and
Pukhlikov \cite{khovanskii-pukhlikov}. Let us discuss their
inversion formula in a slightly greater generality of constructible
functions in the sense of Section \ref{Ss:detailed-overview} as it
was obtained in \cite{schapira2}. Let $\RR\PP^n$ denote the real
projective space, i.e. the space of real lines in $\RR^{n+1}$
passing through zero. Let $\RR\PP^{n\vee}$ be the dual projective
space, i.e. the space of linear hyperplanes in $\RR^{n+1}$. Let $Z$
be the incidence variety, i.e.
$$Z:=\{(l,H)\in\RR\PP^n\times\RR\PP^{n\vee}|\, l\subset H\}.$$ We
have the natural double fibration
\begin{eqnarray}\label{dob-fibr-proj}
\RR\PP^n\overset{q_1}{\leftarrow}Z\overset{q_2}{\to}\RR\PP^{n\vee}.
\end{eqnarray}

The following inversion formula was proved by Khovanskii-Pukhlikov
\cite{khovanskii-pukhlikov} for more restrictive class of
constructible functions; in the present form it was proved by
Schapira \cite{schapira2}.
\begin{theorem}\label{T:kh-p}
Let $R:=q_{2*}q_1^*\colon \cf(\RR\PP^n)\to \cf(\RR\PP^{n\vee})$ be
the Radon transform. Consider the (dual) Radon transform
$$R^t:=q_{1*}q_2^*\colon \cf(\RR\PP^{n\vee})\to \cf(\RR\PP^n).$$
Then for any constructible function $\phi\in \cf(\RR\PP^n)$ one has
\begin{eqnarray}\label{radon-inversion}
(-1)^{n-1}R^tR(\phi)=\phi+\frac{1}{2}((-1)^{n-1}-1)\cdot
\int_{\RR\PP^n}\phi.
\end{eqnarray}
\end{theorem}
\begin{remark}\label{R:radon-inject}
In particular if $n$ is odd then the Radon transform $R$ is
injective, and (\ref{radon-inversion}) is the inversion formula. If
$n$ is even then the kernel of $R$ consists precisely of the
constant functions.
\end{remark}

\hfill

Now we will discuss a valuation theoretic generalization of Theorem
\ref{T:kh-p} which is the main result of this section. By Corollary
\ref{COR:radon-val-2}(1) and Example \ref{Example:double-fibrations}
we have the continuous linear operator $\tilde R$ on generalized
valuations
$$\tilde R:=q_{2*}q_1^*\colon V^{-\infty}(\RR\PP^n)\to V^{-\infty}(\RR\PP^{n\vee}).$$
(Here we just take the kernel $\gamma$ to be equal to the Euler
characteristic $\chi$.) Note that by Corollary
\ref{COR:radon-val-2}(2) $\tilde R$ maps smooth valuations to smooth
ones: $\tilde R(V^{\infty}(\RR\PP^n))\subset
V^\infty(\RR\PP^{n\vee})$. Similarly we define
$$\tilde R^t:=q_{1*}q_2^*\colon V^{-\infty}(\RR\PP^{n\vee})\to V^{-\infty}(\RR\PP^n)$$
which is also continuous and maps smooth valuations to smooth ones.
We have
\begin{theorem}\label{T:radon-val-constr}
For any generalized valuations $\phi\in V^{-\infty}(\RR\PP^n)$ one
has
\begin{eqnarray}\label{inversion-radon-val-constr}
(-1)^{n-1}\tilde R^t\tilde R(\phi)=\phi+\frac{1}{2}((-1)^{n-1}-1)
\left(\int_{\RR\PP^n}\phi\right)\cdot \chi.
\end{eqnarray}
In particular if $n$ is odd then the Radon transform $\tilde R$ is
injective, and (\ref{inversion-radon-val-constr}) is the inversion
formula. If $n$ is even then the kernel of $\tilde R$ consists
precisely of the multiples of the Euler characteristic $\chi$.
\end{theorem}

\hfill

In order to prove Theorem \ref{T:radon-val-constr} we will need some
preparations. Below we will denote for brevity $\RR\PP^n$ by
$\PP^n$, and $\RR\PP^{n\vee}$ by $\PP^{n\vee}$. Let
$$r\colon \RR^{n+1}\backslash 0\to \PP^n$$
be the natural map.
\begin{definition}\label{def1}
A subset $K\subset \PP^n$ is called {\itshape convex} if

(i) there exists a convex cone $\tilde K\subset \RR^{n+1}$ such that
$r(\tilde K\backslash 0)=K$;

(ii) there exists a hyperplane $H\subset \PP^n$ such that $H\cap
K=\emptyset$.
\end{definition}

\begin{remark}\label{rem2}
Equivalently, a subset $K\subset \PP^n$ is convex if and only if
there exists a hyperplane $H\subset \PP^n$ which does not intersect
$K$, and $K$ is a convex subset in the usual (affine) sense of
affine space $\PP^n\backslash H \simeq \RR^n$. This condition does
not depend on a choice of $H$ due to equivalence to Definition
\ref{def1}.
\end{remark}

\begin{definition}\label{def3}
A closed convex subset $K\subset \PP^n$ with non-empty interior is
called {\itshape strictly convex} if

(i) $\pt K$ is smooth;

(ii) $\pt K$ has a strictly positive Gauss curvature when $K$ is
considered as a subset of an affine space $\PP^n\backslash H$, where
$H$ is a hyperplane not intersecting $K$.
\end{definition}
It is easy to see that the last definition is independent of a
choice of hyperplane $H$.

\hfill

\def\dup{\PP^{n\vee}}

The dual projective space $\dup$ can naturally be identified with
the projectivization of the the dual space $(\RR^{n+1})^*$. It can
be identified with the set of hyperplanes in $\PP^n$. Denote by
$$r^\vee\colon (\RR^{n+1})^*\backslash \{0\}\to \dup$$
the natural map.

For a convex compact subset $K\subset\PP^n$ let us define the dual
set $K^\vee\subset \dup$ as follows. Let $\tilde K\subset \RR^{n+1}$
be a closed convex cone satisfying $r(\tilde K\backslash 0)=K$.
$\tilde K$ is defined uniquely up to multiplication by $-1$.
Consider the dual cone in $(\RR^{n+1})^*$
$$\tilde K^o:=\{\xi\in (\RR^{n+1})^*|\,\, \xi(x)\geq 0 \,\,\forall x\in
\tilde K\}.$$ Define
$$K^\vee:=r^\vee(\tilde K^o).$$
If $K$ has non-empty interior then $K^\vee$ is also convex in the
sense of Definition \ref{def1}. In this case $(K^\vee)^\vee=K$. It
is easy to see that the dual of a strictly convex set is strictly
convex.

The following another description of $K^\vee$ will be useful. Assume
that $K$ has non-empty interior. Then $K^\vee$ is equal to the set
of hyperplanes which do not intersect the interior of $K$.

\hfill

Before the proof of Theorem \ref{T:radon-val-constr} we will prove
the following technical proposition.
\begin{proposition}\label{P:radon-convex}
Let $K\subset \PP^n$ be a closed convex set. Then
$$R(\Xi_\cp(\One_K))=\Xi_\cp(\One_{\dup\backslash int(K^\vee)}),$$
where $int (K^\vee)$ denotes the interior of $K^\vee$.
\end{proposition}
\begin{remark}
By the abuse of notation we denote by $\Xi_\cp$ the map defined not
only on the space $\cf(X)$ of constructible functions (in the
subanalytic sense), but also on indicator functions of convex
compact sets and also of compact submanifolds with corners (for all
of them the notion of normal cycle is well defined too). The
non-trivial results from the previous sections will be applied for
indicator functions of convex sets with smooth strictly convex
boundary. These modifications of the map $\Xi_\cp$ are purely
formal, and there is no need to introduce separate notation for
them.
\end{remark}

{\bf Proof of Proposition \ref{P:radon-convex}.} {\underline{Step
1.} Let us show that it is enough to prove the statement for $\pt K$
being smooth strictly convex.

For a closed convex set $K$ there exists a sequence $\{K_j\}$ of
strictly convex sets which converges to $K$ in the Hausdorff metric.
Then it is easy to see that $K_j^\vee\to K^\vee$ in the Hausdorff
metric. But then
\begin{eqnarray*}
\Xi_\cp(\One_{K_j})\to \Xi_\cp(\One_{K}) \mbox{ in } \vmi(\PP^n),\\
\Xi_\cp(\One_{K_j^\vee})\to \Xi_\cp(\One_{K^\vee}) \mbox{ in }
\vmi(\dup).
\end{eqnarray*}

But for a compact convex set $B\subset \dup$ with non-empty
interior, $B\in \cp(\PP^{n\vee})$, one readily has
$$\Xi_\cp (\One_{\dup \backslash int(B)})=\chi-\sigma\left(\Xi_\cp(\One_B)\right)$$
where the operator $\sigma\colon \vmi(\dup)\to \vmi(\dup)$ is a
continuous linear operator introduced in \cite{part4} called the
Euler-Verdier involution. Hence
$$\Xi_\cp(\One_{\dup\backslash int(K_j^\vee)})\to \Xi_\cp(\One_{\dup\backslash
int(K^\vee)})\mbox{ in } \vmi(\dup).$$ Thus since $R\colon
\vmi(\PP^n)\to \vmi(\dup)$ is continuous in the weak topology we may
assume that $\pt K$ is smooth and strictly convex.

\hfill

\underline{Step 2.} $R(\Xi_\cp(\One_K))$ vanishes on $int(K^\vee)$.
Indeed since $int(K^\vee)$ consists of hyperplanes in $\PP^n$ not
intersecting $K$, then $q_1^{-1}(K) \cap
q_2^{-1}(int(K^\vee))=\emptyset$.

\hfill

\underline{Step 3.} Let us show that the restriction of
$R(\Xi_\cp(\One_K))$ to the open set $\dup\backslash K^\vee$ is
equal to $\chi$.

First let us show that
\begin{eqnarray}\label{sept4.5}
q_2\colon q_1^{-1}(\pt K)\cap q_2^{-1}(\dup \backslash K^\vee)\to
\dup\backslash K^\vee
\end{eqnarray}
is a submersion. Let us fix $x_0\in \pt K,\, H_0\in \dup\backslash
K^\vee$ such that $x_0\in H_0$. Let us denote by $A_1$ the tangent
space of $\PP^n$ at $x_0$, $A_2$ the tangent space of $\dup$ at
$H_0$, $B$ the tangent space of $\pt K$ at $x_0$, $C$ the tangent
space of $Z$ at $(x_0,H_0)$, and $D$ the tangent space of
$q_1^{-1}(\pt K)$ at $(x_0,H_0)$. Then we have
\begin{eqnarray*}
B\subset A_1\overset{dq_1}{\leftarrow}C\overset{dq_2}{\to}A_2,\\
D\subset C.
\end{eqnarray*}
We want to show
\begin{eqnarray}\label{sept4.6}
dq_2(D)=A_2.
\end{eqnarray}
But clearly $D=(dq_1)^{-1}(B)$. Denote also $F:=Ker(dq_2)$. Then
(\ref{sept4.6}) is equivalent to
$$(dq_1^{-1})(B)+F=C.$$
This is equivalent to
\begin{eqnarray}\label{sept4.7}
B+(dq_1)(F)=A_1.
\end{eqnarray}
But $F$ is the tangent space at $H_0$ of the fiber $q_2^{-1}(H_0)$.
Hence $(dq_1)(F)$ is the tangent space of $q_1(q_2^{-1}H_0)$ at
$x_0$. But obviously $$q_1(q_2^{-1}H_0)=\{x\in \PP^n|\, x\in
H_0\}=H_0.$$ The condition (\ref{sept4.7}) means that the hyperplane
$H_0$ is transversal to $\pt K$ at $x_0$. But it is easy to see that
any hyperplane $H_0\in \dup\backslash K^\vee$ is transversal to $\pt
K$ since $K$ is convex. Thus we have shown that the map
(\ref{sept4.5}) is a submersion.

\hfill

This implies that for any point $H\in \dup\backslash K^\vee$ there
exists a neighborhood $\cu$ of $H$, a neighborhood $\cv$ of
$q_1^{-1}(K)\cap q_2^{-1}(\cu)$ such that $\cv\subset
q_2^{-1}(\cu)$, and a diffeomorphism
$$\phi\colon \cv\tilde\to \cu\times \RR^{n-1}$$
such that $\phi(q_1^{-1}(K)\cap q_2^{-1}(\cu))=\cu\times B$ where
$B\subset \RR^{n-1}$ is the closed unit ball, and such that the
following diagram is commutative:
\def\pu{p_\cu}
$$\Vtriangle[\cv`\cu\times \RR^{n-1}`\cu;\phi`q_2`p_\cu]$$
where $\pu\colon \cu\times \RR^{n-1}\to \cu$ is the natural
projection. Then Proposition \ref{P:pull-submer-gener-constr} and
Lemma \ref{LL:4} imply that the restriction of
$q_{2*}q_1^*(\Xi_\cp(\One_K))=q_{2*}(\Xi_\cp(\One_{q_1^{-1}(K)}))$
to $\cu$ is equal to $p_{\cu *}(\Xi_\cp(\One_{\cu\times
B}))=\chi(B)\cdot \chi=\chi$. Hence the restriction of
$R(\Xi_\cp(\One_K))$ to $\dup\backslash K^\vee$ is equal to $\chi$.

\hfill

\underline{Step 4.} Now let us fix a point $H^0\in \pt K^\vee$, and
show that the restriction of $R(\Xi_\cp(\One_K))$ to a small
neighborhood of $H^0$ is equal to $\Xi_\cp(\One_{\dup\backslash
int(K^\vee)})$.

Let $\tilde K\subset \RR^{n+1}$ be a closed convex cone over $K$. We
can choose a coordinate system $\{(x_1,\dots,x_{n+1})\}$ in
$\RR^{n+1}$ such that the hyperplane $\{(x_1,\dots,x_n,0)\}$
intersects $\tilde K$ only at 0, and such that the point
$(0,\dots,0,1)\in int(\tilde K)$.

Let $\{y_1,\dots,y_{n+1}\}$ be the dual coordinates in
$(\RR^{n+1})^*$. The the dual cone $\tilde K^o$ is given by
$$\tilde K^o=\{(y_1,\dots,y_{n+1})|\, \sum_{i=1}^{n+1}x_iy_i\geq 0
\, \mbox{ for all } (x_1,\dots,x_{n+1})\in \tilde K\}.$$ It is clear
that if $(y_1,\dots,y_{n+1})\in \tilde K^o\backslash \{0\}$ then
$y_{n+1}>0$. In particular the hyperplane $\{(y_1,\dots,y_n,0)\}$
intersects $\tilde K^o$ only at 0.

Thus we have shown that $K$ is contained in the affine chart of
$\PP^n$ given by $\{[x_1:\dots:x_n:1]\}$ in homogeneous coordinates,
and $K^\vee$ is contained in $\{[y_1:\dots:y_n:1]\}$. We will
identify both these charts with $\RR^n$ using the coordinates
$(x_1,\dots,x_n)$ and $(y_1,\dots,y_n)$ respectively. In these
coordinates
\begin{eqnarray*}
(0,\dots,0)\in K,\\
K^\vee=\{(y_1,\dots,y_n)|\, \sum_{i+1}^n x_iy_i\geq -1 \mbox{ for
all } (x_1,\dots,x_n)\in K\}.
\end{eqnarray*}

Notice that $-K^\vee=\{(y_1,\dots,y_n)|\, \sum_{i+1}^n x_iy_i\leq 1
\mbox{ for all } (x_1,\dots,x_n)\in K\}$ is the dual body of $K$ in
the usual linear sense.

We will reduce our statement to Proposition \ref{PP:5} applying
certain appropriate diffeomorphism. More precisely we will show that
there exist a neighborhood $\cu\subset \dup$ of an arbitrary fixed
point $H^0\in \pt K^\vee$, a neighborhood $\cv$ of $q_1^{-1}(K)\cap
q_2^{-1}(\cu)$ such that $\cv \subset q_2^{-1}(\cu)$, and open
imbeddings
\begin{eqnarray*}
\phi\colon \cv\to \RR^n\times \RR^{n-1},\\
\psi\colon \cu\to \RR^n
\end{eqnarray*}
such that the following two conditions are satisfied:

(1) the diagram is commutative:
$$\square[\cv`\RR^n\times\RR^{n-1}`\cu`\RR^n;\phi`q_2`p`\psi]$$
where $p\colon \RR^n\times \RR^{n-1}\to \RR^n$ is the projection
onto the first factor;

(2) $\phi(q_1^{-1}(K)\cap q_2^{-1}(\cu))=:A$ is a closed convex
subset of $\RR^n\times \RR^{n-1}$ with smooth strictly convex
boundary, and the restriction $p|_A$ is proper.

Clearly if we will show the above, Proposition \ref{PP:5} finishes
the proof of our proposition.

\hfill

For $\eps >0$ let us consider the map
\begin{eqnarray}\label{sept4.8}
\phi_1\colon \pt K^\vee\times (-\eps,\eps)\to \RR^n
\end{eqnarray}
given by $\phi_1(s,t):=(1+t)^{-1}s$. Clearly $\phi_1$ is a
diffeomorphism onto an open neighborhood of $\pt K^\vee$ in $\RR^n$
for $\eps$ small enough.

For any point $s\in \pt K^\vee$ there exists a unique point $x\in
\pt K$ such that $<s,x>=-1$. This defines a map
$$\Phi\colon \pt K^\vee\to \pt K$$
which is a diffeomorphism since $\pt K$ is strictly convex. For
$y\in \RR^n\backslash \{0\}$ let us denote the hyperplane
$$H_y:=\{(x_1,\dots,x_n)\in \RR^n|\, \sum_{i=1}^nx_iy_i=-1\}.$$


If $H=\phi_1(s,t)$, $s\in \pt K^\vee,\, t\in (-\eps,\eps)$, then
$x\in H_s$ if and only if
\begin{eqnarray}\label{oct1}
<s,x>=-(1+t).
\end{eqnarray}

For $s\in \pt K^\vee$ belonging to a small neighborhood of the fixed
point $H^0\in \pt K^\vee$ let us choose a basis
$\xi_1(s),\dots,\xi_{n-1}(s)$ of $s^\perp$ which depends smoothly on
$s$. Then any $x\in H_{\phi_1(s,t)}$ has a decomposition
\begin{eqnarray}\label{oct2}
x=\sum_{j=1}^{n-1}\lam_j\xi_j(s)+(1+t)\Phi(s)
\end{eqnarray}
with $\lam_j\in \RR$ (the coefficient of $\Phi(s)$ is equal to $1+t$
due to (\ref{oct1})). Hence on the incidence variety $Z$ near the
pair $(\Phi(H^0),H^0)$ we have a coordinate system given by
$(s;\lam_1,\dots,\lam_{n-1};t)$.

Let us describe the set $q_1^{-1}(K)$ in these coordinates. For any
$s\in \pt K^\vee$ there exists a unique function
$$F_s\colon s^\perp\to \RR$$
such that near $\Phi(s)$ we have
\begin{eqnarray}\label{oct3}
K=\{h-\tau\Phi(s)|\, h\in s^\perp,\, \tau\geq F_s(h)\}.
\end{eqnarray}

This function $F_s$ is $C^\infty$-smooth and strictly convex.
Moreover the function
$$F(s,h):=F_s(h)$$
is $C^\infty$-smooth on the set $\{(s,h)|\, s\in \pt K^\vee,\, h\in
s^\perp\}$ near $s=H^0,\, h=0$. Also $0\in s^\perp$ is a critical
point of $F_s$, and $F_s(0)=-1$.

Assume we are given a pair $(x,y)\in Z\subset \PP^n\times
\PP^{n\vee}$ corresponding to $(s;\lam_1,\dots,\lam_{n-1};t)$. Then
$x\in K$ if and only if
\begin{eqnarray*}
-t\geq 1+F_s(\sum_j\lam_j\xi_j(s)).
\end{eqnarray*}
Indeed by (\ref{oct2}) $x=\sum\lam_j\xi_j(s)+(1+t)\Phi(s)$. Then by
(\ref{oct3}), $x\in K$ if and only if $-(1+t)\geq
F_s\left(\sum\lam_j\xi_j(s)\right)$.

\hfill

Now let us denote for brevity
$$G_s(\lam_1,\dots,\lam_{n-1}):=1+F_s\left(\sum_{j=1}^{n-1}\lam_j\xi_j(s)\right).$$
Thus we have shown that in the above coordinates on $Z$ near
$(0,H^0)$ the set $q_1^{-1}(K)$ is given by
$$q_1^{-1}(K)=\{(s;\lam_1,\dots,\lam_{n-1};t)\in \pt K^\vee\times
\RR^{n-1}\times \RR|\, -t\geq G_s(\lam_1,\dots,\lam_{n-1})\}.$$

It is clear that the map $(s,\lam_1,\dots,\lam_{n-1})\mapsto
G_s(\lam_1,\dots,\lam_{n-1})$ is $C^\infty$-smooth, and when $s$ is
fixed the function $G_s$ is strictly convex. In particular for fixed
$s$, $G_s$ is a Morse function. Notice also that $G_s$ has a minimum
at $0\in \RR^{n-1}$, and $G_s(0)=0$.

By the Morse lemma with parameters (see e.g. \cite{bruce-giblin}, p.
97) there exists a $C^\infty$ map
$$D\colon \pt K^\vee\times \RR^{n-1}\to \pt K^\vee\times \RR^{n-1}$$
of the form $D(s;\mu_1,\dots,\mu_{n-1})=(s;\lam_1,\dots,\lam_{n-1})$
such that for $s$ close to $H^0$ the map $\RR^{n-1}\to \RR^{n-1}$
given by $(\mu_1,\dots, \mu_{n-1})\mapsto
D(s;\mu_1,\dots,\mu_{n-1})=:D_s(\mu_1,\dots,\mu_{n-1})$ is a
diffeomorphism of a neighborhood of $0\in \RR^{n-1}$, and
$$(G_s\circ D_s)(\mu_1,\dots,
\mu_{n-1})=\mu_1^2+\dots+\mu_{n-1}^2.$$

\hfill

In the new coordinates defined by
$$(s;\mu_1,\dots,
\mu_{n-1};\tau):=(s;D_s^{-1}(\lam_1,\dots,\lam_{n-1});-t)$$ the set
$q_1^{-1}(K)$ is locally given by
$$q_1^{-1}(K)=\{(s;\mu_1,\dots,\mu_{n-1};\tau)|\, s\in \pt K^\vee,\,
\tau\geq \mu_1^2+\dots+\mu_{n-1}^2\}.$$ Moreover in these
coordinates the map $q_2\colon Z\to \dup$ is given by
$$q_2(s;\mu_1,\dots,
\mu_{n-1};\tau)=\phi_1(s,-\tau).$$

\hfill

Let us fix a coordinate system $\{(s_1,\dots,s_{n-1})\}$ on $\pt
K^\vee$ near $H^0$ such that $(0,\dots,0)$ corresponds to $H^0$.
Finally let us introduce the last our coordinate system
$\{(s_1,\dots,s_{n-1};\mu_1,\dots,\mu_{n-1};r)\}$ on $Z$ near
$(\Phi(H^0),H^0)$ when $r$ is defined by
$$r:=\tau+\sum_{j=1}^{n-1}s_j^2.$$
Then $\{(s_1,\dots,s_{n-1};r)\}$ is a coordinate system on
$\PP^{n\vee}$ near $H^0$. Then $q_1^{-1}(K)$ is given by
\begin{eqnarray}\label{oct4}
q_1^{-1}(K)=\{(s_1,\dots,s_{n-1};\mu_1,\dots,\mu_{n-1};r)|\, r\geq
\sum_{j=1}^{n-1}s_j^2+\sum_{j=1}^{n-1}\mu_j^2\},
\end{eqnarray}
and the map $q_2$ is given by
\begin{eqnarray}\label{oct5}
q_2(s_1,\dots,s_{n-1};\mu_1,\dots,\mu_{n-1};r)=(s_1,\dots,s_{n-1};r).
\end{eqnarray}
Thus the formula (\ref{oct4}) defines a Euclidean ball in
$\RR^{2n-1}$, and the map $q_2$ given by (\ref{oct5}) is just a
projection. This coordinate system gives the maps $\phi,\psi$ as
promised in the beginning of Step 4. Hence Step 4 is completed.

Steps 1-4 imply Proposition \ref{P:radon-convex}. \qed

\hfill

{\bf Proof of Theorem \ref{T:radon-val-constr}.} The linear operator
$$\tilde R^t\tilde R\colon \vmi(\PP^n)\to \vmi(\PP^n)$$
is continuous in the weak topology by Corollary \ref{COR:radon-val-2}
and Example \ref{Example:double-fibrations}(2).
It is not hard to see (as in
\cite{part4}, Section 8.1) that linear combinations of valuations of
the form $\Xi_\cp(\One_K)$, with $K\subset \PP^n$ being a compact
convex polyhedral subset, are dense in $V^{-\infty}(\PP^n)$. Hence
it is sufficient to prove the inversion formula
(\ref{inversion-radon-val-constr}) for such valuations. By
Proposition \ref{P:radon-convex} the operators $\tilde R,\tilde R^t$
coincide on such valuations with the Radon transform on
constructible functions with respect to the Euler characteristic.
Note also that again by Proposition \ref{P:radon-convex} the
operator $\tilde R$ maps valuations of this form on $\PP^n$ to
valuations of this form on $\dup$: indeed if $K\subset \PP^n$ is a
compact convex polyhedral subset then the dual $K^\vee$ is also
compact convex and polyhedral. But by the Khovanskii-Pukhlikov
theorem \cite{khovanskii-pukhlikov} the formula
(\ref{inversion-radon-val-constr}) is satisfied for the integration
of constructible functions with respect to the Euler characteristic.
\qed

\subsection{Relations to the Crofton style integral
geometry.}\label{Ss:chern} In this section we explain the relation
of operations on valuations to yet another classical kind of
integral geometry. However the exposition in this section is less
formal and rigorous than in the rest of the article.

Let us remind the classical Crofton intersection formula and the
general kinematic formulas of Chern \cite{chern-general-kinemat}.
(For the Crofton formula we refer to Schneider's book
\cite{schneider-book}, \S 4.5; for the kinematic formulas see also
Santalo's book \cite{santalo} and references therein.)

\hfill

First we will need to remind the notion of intrinsic volume.
Consider the standard Euclidean space $\RR^n$. Let $M\subset \RR^n$
be a compact submanifold with corners. For $\eps>0$ let us denote by
$M_\eps$ the $\eps$-neighborhood of $M$ (or $\eps$-tube in some
terminology). For sufficiently small $\eps>0$ the volume of $M_\eps$
is a polynomial of degree $n$ with respect to $\eps$:
$$vol(M_\eps)=\sum_{i=0}^n\eps^{n-i}\kappa_{n-i} V_i(M)$$
where $\kappa_j$ is the $j$-dimensional volume of the
$j$-dimensional unit ball. The coefficient $V_i(M)$ is called the
$i$-th intrinsic volume of $M$. It is  not hard to see that $V_i$ is
a smooth valuation. $V_0$ is the Euler characteristic, and $V_n$ is
the $n$-dimensional volume. But the most surprising property of
$V_i$ is due to H. Weyl \cite{weyl}; it says that $V_i(M)$ is
independent of an {\itshape isometric} imbedding of $M$ into
$\RR^n$. It fact $V_i$ can be written in terms of the curvature
tensor of $M$ (see Weyl's paper \cite{weyl}, see also
\cite{hormander-convexity}).

\hfill

Let us denote by $\overline{O(n)}$ the group of all affine
isometries of $\RR^n$. Let $m$ denote the Haar measure on
$\overline{O(n)}$ (we do not specify the normalization of it since
we will not write down here the constants anyway). Let us denote by
$\ca Gr_{k,n}$ the Grassmannian of affine $k$-dimensional subspaces
in $\RR^n$. The group $\overline{O(n)}$ acts transitively on it. Let
$\mu$ denote a Haar measure of it with respect to the action of the
group $\overline{O(n)}$. Then the Crofton intersection formula says
that
\begin{eqnarray}\label{crofton}
\int_{\ca Gr_{k,n}}V_i(M\cap E) d\mu(E)=\alp_{i,k,n} V_{n+i-k}
\end{eqnarray}
where $\alp_{i,k,n}$ depends on $i,k,n$ only.

\hfill

The general kinematic formula of Chern \cite{chern-general-kinemat}
says that for two compact submanifolds with boundary $M,N$ one has
\begin{eqnarray}\label{kinemat}
\int_{\overline{O(n)}}V_i(M\cap g(N))dm(g)=\sum_j c_{i,j,n}
V_j(M)V_{n+i-j}(N)
\end{eqnarray}
where the constants $c_{i,j,n}$ depend only on $i,j,n$ and can be
written down explicitly.

\hfill

Let us explain how to rewrite the left hand side in the formulas of
Crofton and Chern using the operations of pull-back, push-forward,
and the product of valuations.

Let us start with the Crofton formula. Let us denote by
$$Z:=\{(x,E)\in\RR^n\times \ca Gr_{k,n}|\, x\in E\}$$
the incidence variety. We have the double fibration
$\RR^n\overset{q_1}{\leftarrow}Z\overset{q_2}{\to}\ca Gr_{k,n}$.
Clearly $q_1$ is proper.

One can see the following formula for the left hand side of
(\ref{crofton})
\begin{eqnarray}\label{crofton-valuat}
\int_{\ca Gr_{k,n}}V_i(M\cap E) d\mu(E)=\left(V_i\cdot
q_{1*}q_2^*\mu\right)(M)
\end{eqnarray}
where the pull-back $q_2^*$, the push-forward $q_{1*}$, and the
product are taken in the space of valuations. Let us also notice
that the formula (\ref{crofton-valuat}) holds in a greater
generality: not only for a Haar measure $\mu$ but for {\itshape any}
smooth measure on the affine Grassmannian.

\hfill

Let us return to the general kinematic formula (\ref{kinemat}). The
left hand side of it can be rewritten as follows. Let us consider
the double fibration
$$\RR^n\overset{p}{\leftarrow}\overline{O(n)}\times
\RR^n\overset{a}{\to}\RR^n$$ where $p$ is the projection $p(g,x)=x$,
and $a$ is the action map $a(g,x)=g(x)$.

\begin{eqnarray}\label{kinemat-valuat}
\int_{\overline{O(n)}}V_i(M\cap
g(N))dm(g)=\left(p_*\left(a^*[\Xi_\cp(\One_N)]\cdot [m\boxtimes
V_i]\right)\right)(M)
\end{eqnarray}
where $\boxtimes$ in the right hand side is the exterior product on
valuations introduced in Section \ref{S:exter-prod}. Observe that
the operator $[\phi\mapsto p_*\left(a^*\phi\cdot (m\boxtimes
V_i)\right)]$ takes values in $\overline{O(n)}$-invariant
valuations.

\hfill

Finally we would like to rewrite in this language some integral
geometric expressions in complex integral geometry. They are not
fully understood and still under investigation. The case of the
affine complex space $\CC^n$ is understood better (see
\cite{alesker-jdg-03}, \cite{bernig-fu-hermit}), and we consider
here the less studied case of the complex projective space
$\CC\PP^n$. To introduce the complex version of the Crofton formula
(actually only the left hand side of it, as the right hand side is
still to be understood), let us denote by $\grc_{k+1,n+1}$ the
Grassmannian of complex linear subspaces in $\CC^{n+1}$ of complex
dimension $k+1$. Every element of it can be identified with the
$k$-dimensional complex projective subspace in $\CC\PP^n$. This
Grassmannian is acted by the unitary group $U(n+1)$ in a transitive
way. We will denote by $\mu$ a Haar measure on the Grassmannian. Let
$M\subset \CC\PP^n$ be a smooth submanifold with boundary (or even
with corners). In this notation, one would like to have a closed
formula for
\begin{eqnarray}\label{complex-crofton}
\int_{E\in\grc_{k+1,n+1}}V_i(M\cap E)d\mu(E).
\end{eqnarray}
Let us consider the double fibration
\begin{eqnarray*}
\CC\PP^n\overset{q_1}{\leftarrow}Z\overset{q_2}{\to}\grc_{k+1,n+1}
\end{eqnarray*}
where $Z:=\{(l,E)\in \CC\PP^n\times \grc_{k+1,n+1}|\, l\subset E\}$
is the incidence variety again. Then the integral
(\ref{complex-crofton}) can be rewritten as
\begin{eqnarray}\label{complex-crofton-2}
(V_i\cdot q_{1*}q_2^*(\mu))(M).
\end{eqnarray}
This expression is an $U(n+1)$-invariant valuation on $\CC\PP^n$.
The similarity with the expression (\ref{crofton-valuat}) in the
Euclidean case is evident. Notice that in the recent article
\cite{abardia-gallego-solanes} a closed formula for
(\ref{complex-crofton}) in the case $i=0$ has been obtained.

Analogously the left hand sides of Crofton and kinematic type
formulas can be written also on the complex hyperbolic space where
they are also not well understood.

\def\gr{{}\!^ {\KK} Gr}
\section{Appendix.}\label{S:appendix}\setcounter{subsection}{1} In this appendix we prove that
some pairs of the Grassmannians satisfy the assumptions of Corollary
\ref{COR:radon-val-2} (see Proposition \ref{P:submersion} below).

Let $\KK$ denote either $\RR$, $\CC$, or $\HH$ (where $\HH$ is the
non-commutative field of quaternions). Let $V$ be a $\KK$-vector
space of $\KK$-dimension $n$. Let $\gr_k$ denote the Grassmannian of
$\KK$-linear subspaces in $V$ of $\KK$-dimension $k$. Fix now $k,m$
with $0<k<m<n$. Let
$$Z\subset \gr_k\times \gr_m$$
be the incidence variety, i.e. $Z:=\{(l,H)|\, l\subset H\}$. It is
well known that there exist canonical isomorphisms
\begin{eqnarray*}
T_l(\gr_k)=Hom_\KK(l,V/l)=l^*\otimes_\KK V/l,\\
T_H(\gr_m)=Hom_\KK(H,V/H)=H^*\otimes_\KK V/H.
\end{eqnarray*}
Notice that if $\KK=\HH$, we consider $V,l,H$ as right $\KK$-vector
spaces. Then $$l^*:=Hom_\KK(l,\KK)\tilde\to Hom_\RR(l,\RR)$$ is a
left $\KK$-vector space, and similarly for $H^*$.

For a point $(l,H)\in Z$ let us denote by
\begin{eqnarray*}
i\colon l\inj H,\\
p\colon V/l\twoheadrightarrow V/H
\end{eqnarray*}
the natural maps.

\begin{proposition}\label{P:tangent}
For $(l,H)\in Z$ the tangent space $T_{(l,H)}Z$ is naturally
isomorphic to the space of pairs
$$(\xi,\eta)\in T_l(\gr_k)\oplus T_H(\gr_m)=Hom_\KK(l,V/l)\oplus
Hom_\KK(H,V/H)$$ which make the following diagram commutative
$$\square[l`V/l`H`V/H;\xi`i`p`\eta].$$
\end{proposition}
{\bf Proof.} Let us fix splittings
\begin{eqnarray*}
H=l\oplus A,\\
V=H\oplus B=l\oplus A\oplus B.
\end{eqnarray*}
Consider the open neighborhood of $l$ (resp. $H$) consisting of
subspaces which are graphs of $\KK$-linear maps
$$l\to A\oplus B,\, (\mbox{ resp. } H\to B).$$
Given such $\phi\colon l\to A\oplus B,\, \psi\colon H\to B$, let us
denote by $\tilde\phi\in \gr_k$, $\tilde\psi\in \gr_m$ the
corresponding subspaces. The map $\phi$ has the form
$\phi=(\phi_A,\phi_B)$ with
$$\phi_A\colon l\to A,\, \phi_B\colon l\to B.$$
Then $\tilde\phi\subset \tilde\psi$ if and only if for any $x\in l$
the vector $(x,\phi_A(x),\phi_B(x))$ belongs to the graph of $\psi$.
Namely there exists $y=(z,a)\in H=l\oplus A$ such that
$$(x,\phi_A(x),\phi_B(x))=(z,a,\psi(z,a)).$$
Equivalently $z=x,a=\phi_A(x)$, and
\begin{eqnarray}\label{E:small}
\phi_B(x)=\psi(x,\phi_A(x))=\psi(x,0_A)+\psi(0_l,\phi_A(x)).
\end{eqnarray}
In order to describe the tangent space $T_{(l,H)}Z$ we have to
assume that $\phi,\psi$ are infinitesimally small. Then in
(\ref{E:small}) we have to neglect the summand $\psi(0_l,\phi_A(x))$
which has the second order of smallness. Thus infinitesimally small
$\phi,\psi$ satisfy
$$\phi_B(x)=\psi(x,0_A) \mbox{ for any } x\in l.$$
Equivalently
$$p\circ \phi=\psi\circ i.$$
This proves the proposition. \qed

\hfill

By Proposition \ref{P:tangent} we have the short exact sequence
\begin{eqnarray}\label{E:e-seq}
0\to T_{(l,H)}Z\to (l^*\otimes_\KK V/l)\oplus (H^*\otimes_\KK
V/H)\to l^*\otimes_\KK V/H\to 0
\end{eqnarray}
where the third arrow is $(id_{l^*}\otimes p)\oplus(i^*\otimes
-id_{V/H})$. Let us describe now the fiber over $(l,H)$ of the
{\itshape conormal} bundle of $Z\subset \gr_k\times \gr_m$. Clearly
this fiber is equal to $(T_{(l,H)}Z)^\perp$. Observe also that
$(V/l)^*=l^\perp$, $(V/H)^*=H^\perp$. Then taking the dual sequence
of (\ref{E:e-seq}) we immediately deduce the following corollary.
\begin{corollary}\label{Cor:fiber}
The fiber $(T_{(l,H)}Z)^\perp$ of the conormal bundle of $Z\subset
\gr_k\times \gr_m$ over a point $(l,H)\in Z$ is canonically
isomorphic to the image of the injective map
$$l\otimes_\KK H^\perp\to (l\otimes_\KK l^\perp)\oplus(H\otimes_\KK
H^\perp)$$ where the map is $(id_l\otimes p^*)\oplus (i\otimes
-id_{H^\perp})$.
\end{corollary}

\hfill

Now we will prove the main result of this appendix.
\begin{proposition}\label{P:submersion} Let $k=1$ or $n-1$. Let
$m\ne k$. Then the natural projection
$$f\colon T^*_Z(\gr_k\times \gr_m)\backslash \underline{0}\to
T^*(\gr_k)\backslash \underline{0}$$ is a submersion.
\end{proposition}
{\bf Proof.} By taking the orthogonal complement if necessary, we
may and assume that $k=1$. Recall that
$$T^*(\gr_k)=\{(l,u)|\, l\in \gr_k,\, u\in l\otimes_\KK l^\perp\}.$$
By Corollary \ref{Cor:fiber} we can identify
\begin{eqnarray}\label{E:Ap1}
T^*_Z(\gr_k\times\gr_m)\simeq \{(l,H,v)|\, l\subset H,\, v\in
l\otimes_\KK H^\perp\}. \end{eqnarray} Then the projection $f$ is
given by
\begin{eqnarray}\label{E:Ap2}
f(l,H,v)=(l,(id_l\otimes p^*)(v))
\end{eqnarray}
where, we remind, $id_l\otimes p^*\colon l\otimes_\KK H^\perp\to
l\otimes_\KK l^\perp$ is the natural imbedding.

Let us denote for brevity
\begin{eqnarray*}
S:=T^*_Z(\gr_k\times \gr_m)\backslash\underline{0},\\
T:=T^*(\gr_k)\backslash\underline{0}.
\end{eqnarray*}
We have the commutative diagram
$$\qtriangle[S`T`\gr_k;f`g_1`g_2]$$
where $g_1,g_2$ are the obvious maps. Moreover all maps in this
diagram commute with the natural action of the group $GL(n,\KK)$.
Since the last group acts transitively on $\gr_k$, it follows that
$f\colon S\to T$ is a submersion if and only if for any
(equivalently, for some) $l\in \gr_k$ the restriction of $f$
$$f|_{g_1^{-1}(l)}\colon g_1^{-1}(l)\to g_2^{-1}(l)$$
is a submersion. But for a fixed $l\in \gr_k$ we have
\begin{eqnarray*}
g_1^{-1}(l)=\{(H,v)|\, H\supset l,\, v\in (l\otimes_\KK
H^\perp)\backslash\{0\}\},\\
g_2^{-1}(l)=(l\otimes_\KK l^\perp)\backslash\{0\}.
\end{eqnarray*}
Then the restriction of $f$ to $g_1^{-1}(l)$ is given by
$$f(H,v)=(id_l\otimes p^*)(v).$$

Let us denote by $G\subset GL(n,\KK)$ the stabilizer of $l$. Since
$k=1$, $G$ acts transitively on $g_2^{-1}(l)$. Moreover the map $f$
commutes with this action. This readily implies that
$f|_{g_1^{-1}(l)}\colon g_1^{-1}(l)\to g_2^{-1}(l)$ is a submersion.
Hence $f$ is a submersion. \qed

\end{document}

%% file: diagram.tex
\makeatletter

\def\diagram{\m@th\leftwidth=\z@ \rightwidth=\z@ \topheight=\z@
\botheight=\z@ \setbox\@picbox\hbox\bgroup}

\def\enddiagram{\egroup\wd\@picbox\rightwidth\unitlength
\ht\@picbox\topheight\unitlength \dp\@picbox\botheight\unitlength
\hskip\leftwidth\unitlength\box\@picbox}

\def\bfig{\begin{diagram}}
\def\efig{\end{diagram}}
\newcount\wideness \newcount\leftwidth \newcount\rightwidth
\newcount\highness \newcount\topheight \newcount\botheight

\def\ratchet#1#2{\ifnum#1<#2 \global #1=#2 \fi}

\def\putbox(#1,#2)#3{%
\horsize{\wideness}{#3} \divide\wideness by 2
{\advance\wideness by #1 \ratchet{\rightwidth}{\wideness}}
{\advance\wideness by -#1 \ratchet{\leftwidth}{\wideness}}
\vertsize{\highness}{#3} \divide\highness by 2
{\advance\highness by #2 \ratchet{\topheight}{\highness}}
{\advance\highness by -#2 \ratchet{\botheight}{\highness}}
\put(#1,#2){\makebox(0,0){$#3$}}}

\def\putlbox(#1,#2)#3{%
\horsize{\wideness}{#3}
{\advance\wideness by #1 \ratchet{\rightwidth}{\wideness}}
{\ratchet{\leftwidth}{-#1}}
\vertsize{\highness}{#3} \divide\highness by 2
{\advance\highness by #2 \ratchet{\topheight}{\highness}}
{\advance\highness by -#2 \ratchet{\botheight}{\highness}}
\put(#1,#2){\makebox(0,0)[l]{$#3$}}}

\def\putrbox(#1,#2)#3{%
\horsize{\wideness}{#3}
{\ratchet{\rightwidth}{#1}}
{\advance\wideness by -#1 \ratchet{\leftwidth}{\wideness}}
\vertsize{\highness}{#3} \divide\highness by 2
{\advance\highness by #2 \ratchet{\topheight}{\highness}}
{\advance\highness by -#2 \ratchet{\botheight}{\highness}}
\put(#1,#2){\makebox(0,0)[r]{$#3$}}}

\def\adjust[#1]{} 

\newcount \coefa
\newcount \coefb
\newcount \coefc
\newcount\tempcounta
\newcount\tempcountb
\newcount\tempcountc
\newcount\tempcountd
\newcount\xext
\newcount\yext
\newcount\xoff
\newcount\yoff
\newcount\gap%
\newcount\arrowtypea
\newcount\arrowtypeb
\newcount\arrowtypec
\newcount\arrowtyped
\newcount\arrowtypee
\newcount\height
\newcount\width
\newcount\xpos
\newcount\ypos
\newcount\run
\newcount\rise
\newcount\arrowlength
\newcount\halflength
\newcount\arrowtype
\newdimen\tempdimen
\newdimen\xlen
\newdimen\ylen
\newsavebox{\tempboxa}%
\newsavebox{\tempboxb}%
\newsavebox{\tempboxc}%

\newdimen\w@dth

\def\setw@dth#1#2{\setbox\z@\hbox{\m@th$#1$}\w@dth=\wd\z@
\setbox\@ne\hbox{\m@th$#2$}\ifnum\w@dth<\wd\@ne \w@dth=\wd\@ne \fi
\advance\w@dth by 1.2em}

\def\t@^#1_#2{\allowbreak\def\n@one{#1}\def\n@two{#2}\mathrel
{\setw@dth{#1}{#2}
\mathop{\hbox to \w@dth{\rightarrowfill}}\limits
\ifx\n@one\empty\else ^{\box\z@}\fi
\ifx\n@two\empty\else _{\box\@ne}\fi}}
\def\t@@^#1{\@ifnextchar_{\t@^{#1}}{\t@^{#1}_{}}}
\def\to{\@ifnextchar^{\t@@}{\t@@^{}}}

\def\t@left^#1_#2{\def\n@one{#1}\def\n@two{#2}\mathrel{\setw@dth{#1}{#2}
\mathop{\hbox to \w@dth{\leftarrowfill}}\limits
\ifx\n@one\empty\else ^{\box\z@}\fi
\ifx\n@two\empty\else _{\box\@ne}\fi}}
\def\t@@left^#1{\@ifnextchar_{\t@left^{#1}}{\t@left^{#1}_{}}}
\def\toleft{\@ifnextchar^{\t@@left}{\t@@left^{}}}

\def\two@^#1_#2{\allowbreak
\def\n@one{#1}\def\n@two{#2}\mathrel{\setw@dth{#1}{#2}
\mathop{\vcenter{\lineskip\z@\baselineskip\z@
                 \hbox to \w@dth{\rightarrowfill}%
                 \hbox to \w@dth{\rightarrowfill}}%
       }\limits
\ifx\n@one\empty\else ^{\box\z@}\fi
\ifx\n@two\empty\else _{\box\@ne}\fi}}
\def\tw@@^#1{\@ifnextchar _{\two@^{#1}}{\two@^{#1}_{}}}
\def\two{\@ifnextchar ^{\tw@@}{\tw@@^{}}}

\def\tofr@^#1_#2{\def\n@one{#1}\def\n@two{#2}\mathrel{\setw@dth{#1}{#2}
\mathop{\vcenter{\hbox to \w@dth{\rightarrowfill}\kern-1.7ex
                 \hbox to \w@dth{\leftarrowfill}}%
       }\limits
\ifx\n@one\empty\else ^{\box\z@}\fi
\ifx\n@two\empty\else _{\box\@ne}\fi}}
\def\t@fr@^#1{\@ifnextchar_ {\tofr@^{#1}}{\tofr@^{#1}_{}}}
\def\tofro{\@ifnextchar^ {\t@fr@}{\t@fr@^{}}}

\def\mon{\mathop{\m@th\hbox to
      14.6\P@{\lasyb\char'51\hskip-2.1\P@$\arrext$\hss
$\mathord\rightarrow$}}\limits} 
\def\leftmono{\mathrel{\m@th\hbox to
14.6\P@{$\mathord\leftarrow$\hss$\arrext$\hskip-2.1\P@\lasyb\char'50%
}}\limits} 
\mathchardef\arrext="0200       

\def\settypes(#1,#2,#3){\arrowtypea#1 \arrowtypeb#2 \arrowtypec#3}
\def\settoheight#1#2{\setbox\@tempboxa\hbox{#2}#1\ht\@tempboxa\relax}%
\def\settodepth#1#2{\setbox\@tempboxa\hbox{#2}#1\dp\@tempboxa\relax}%
\def\settokens`#1`#2`#3`#4`{%
     \def\tokena{#1}\def\tokenb{#2}\def\tokenc{#3}\def\tokend{#4}}
\def\setsqparms[#1`#2`#3`#4;#5`#6]{%
\arrowtypea #1
\arrowtypeb #2
\arrowtypec #3
\arrowtyped #4
\width #5
\height #6
}
\def\setpos(#1,#2){\xpos=#1 \ypos#2}

\def\settriparms[#1`#2`#3;#4]{\settripairparms[#1`#2`#3`1`1;#4]}%

\def\settripairparms[#1`#2`#3`#4`#5;#6]{%
\arrowtypea #1
\arrowtypeb #2
\arrowtypec #3
\arrowtyped #4
\arrowtypee #5
\width #6
\height #6
}

\def\resetparms{\settripairparms[1`1`1`1`1;500]\width 500}

\resetparms

\def\mvector(#1,#2)#3{
\put(0,0){\vector(#1,#2){#3}}%
\put(0,0){\vector(#1,#2){26}}%
}
\def\evector(#1,#2)#3{{
\arrowlength #3
\put(0,0){\vector(#1,#2){\arrowlength}}%
\advance \arrowlength by-30
\put(0,0){\vector(#1,#2){\arrowlength}}%
}}

\def\horsize#1#2{%
\settowidth{\tempdimen}{$#2$}%
#1=\tempdimen
\divide #1 by\unitlength
}

\def\vertsize#1#2{%
\settoheight{\tempdimen}{$#2$}%
#1=\tempdimen
\settodepth{\tempdimen}{$#2$}%
\advance #1 by\tempdimen
\divide #1 by\unitlength
}

\def\putvector(#1,#2)(#3,#4)#5#6{{%
\ifnum3<\arrowtype
\putdashvector(#1,#2)(#3,#4)#5\arrowtype
\else
\ifnum\arrowtype<-3
\putdashvector(#1,#2)(#3,#4)#5\arrowtype
\else
\xpos=#1
\ypos=#2
\run=#3
\rise=#4
\arrowlength=#5
\ifnum \arrowtype<0
    \ifnum \run=0
        \advance \ypos by-\arrowlength
    \else
        \tempcounta \arrowlength
        \multiply \tempcounta by\rise
        \divide \tempcounta by\run
        \ifnum\run>0
            \advance \xpos by\arrowlength
            \advance \ypos by\tempcounta
        \else
            \advance \xpos by-\arrowlength
            \advance \ypos by-\tempcounta
        \fi
    \fi
    \multiply \arrowtype by-1
    \multiply \rise by-1
    \multiply \run by-1
\fi
\ifcase \arrowtype
\or \put(\xpos,\ypos){\vector(\run,\rise){\arrowlength}}%
\or \put(\xpos,\ypos){\mvector(\run,\rise)\arrowlength}%
\or \put(\xpos,\ypos){\evector(\run,\rise){\arrowlength}}%
\fi\fi\fi
}}

\def\putsplitvector(#1,#2)#3#4{
\xpos #1
\ypos #2
\arrowtype #4
\halflength #3
\arrowlength #3
\gap 140
\advance \halflength by-\gap
\divide \halflength by2
\ifnum\arrowtype>0
   \ifcase \arrowtype
   \or \put(\xpos,\ypos){\line(0,-1){\halflength}}%
       \advance\ypos by-\halflength
       \advance\ypos by-\gap
       \put(\xpos,\ypos){\vector(0,-1){\halflength}}%
   \or \put(\xpos,\ypos){\line(0,-1)\halflength}%
       \put(\xpos,\ypos){\vector(0,-1)3}%
       \advance\ypos by-\halflength
       \advance\ypos by-\gap
       \put(\xpos,\ypos){\vector(0,-1){\halflength}}%
   \or \put(\xpos,\ypos){\line(0,-1)\halflength}%
       \advance\ypos by-\halflength
       \advance\ypos by-\gap
       \put(\xpos,\ypos){\evector(0,-1){\halflength}}%
   \fi
\else \arrowtype=-\arrowtype
   \ifcase\arrowtype
   \or \advance \ypos by-\arrowlength
       \put(\xpos,\ypos){\line(0,1){\halflength}}%
       \advance\ypos by\halflength
       \advance\ypos by\gap
       \put(\xpos,\ypos){\vector(0,1){\halflength}}%
   \or \advance \ypos by-\arrowlength
       \put(\xpos,\ypos){\line(0,1)\halflength}%
       \put(\xpos,\ypos){\vector(0,1)3}%
       \advance\ypos by\halflength
       \advance\ypos by\gap
       \put(\xpos,\ypos){\vector(0,1){\halflength}}%
   \or \advance \ypos by-\arrowlength
       \put(\xpos,\ypos){\line(0,1)\halflength}%
       \advance\ypos by\halflength
       \advance\ypos by\gap
       \put(\xpos,\ypos){\evector(0,1){\halflength}}%
   \fi
\fi
}

\def\putmorphism(#1)(#2,#3)[#4`#5`#6]#7#8#9{{%
\run #2
\rise #3
\ifnum\rise=0
  \puthmorphism(#1)[#4`#5`#6]{#7}{#8}#9%
\else\ifnum\run=0
  \putvmorphism(#1)[#4`#5`#6]{#7}{#8}#9%
\else
\setpos(#1)%
\arrowlength #7
\arrowtype #8
\ifnum\run=0
\else\ifnum\rise=0
\else
\ifnum\run>0
    \coefa=1
\else
   \coefa=-1
\fi
\ifnum\arrowtype>0
   \coefb=0
   \coefc=-1
\else
   \coefb=\coefa
   \coefc=1
   \arrowtype=-\arrowtype
\fi
\width=2
\multiply \width by\run
\divide \width by\rise
\ifnum \width<0  \width=-\width\fi
\advance\width by60
\if l#9 \width=-\width\fi
\putbox(\xpos,\ypos){#4}
{\multiply \coefa by\arrowlength
\advance\xpos by\coefa
\multiply \coefa by\rise
\divide \coefa by\run
\advance \ypos by\coefa
\putbox(\xpos,\ypos){#5} }%
{\multiply \coefa by\arrowlength
\divide \coefa by2
\advance \xpos by\coefa
\advance \xpos by\width
\multiply \coefa by\rise
\divide \coefa by\run
\advance \ypos by\coefa
\if l#9%
   \putrbox(\xpos,\ypos){#6}%
\else\if r#9%
   \putlbox(\xpos,\ypos){#6}%
\fi\fi }%
{\multiply \rise by-\coefc
\multiply \run by-\coefc
\multiply \coefb by\arrowlength
\advance \xpos by\coefb
\multiply \coefb by\rise
\divide \coefb by\run
\advance \ypos by\coefb
\multiply \coefc by70
\advance \ypos by\coefc
\multiply \coefc by\run
\divide \coefc by\rise
\advance \xpos by\coefc
\multiply \coefa by140
\multiply \coefa by\run
\divide \coefa by\rise
\advance \arrowlength by\coefa
\ifcase\arrowtype
\or \put(\xpos,\ypos){\vector(\run,\rise){\arrowlength}}%
\or \put(\xpos,\ypos){\mvector(\run,\rise){\arrowlength}}%
\or \put(\xpos,\ypos){\evector(\run,\rise){\arrowlength}}%
\fi}\fi\fi\fi\fi}}

\newcount\numbdashes \newcount\lengthdash \newcount\increment

\def\howmanydashes{
\numbdashes=\arrowlength \lengthdash=40
\divide\numbdashes by \lengthdash
\lengthdash=\arrowlength
\divide\lengthdash by \numbdashes
\increment=\lengthdash
\multiply\lengthdash by 3
\divide\lengthdash by 5
}

\def\putdashvector(#1)(#2,#3)#4#5{%
\ifnum#3=0 \putdashhvector(#1){#4}#5
\else
\ifnum#2=0
\putdashvvector(#1){#4}#5\fi\fi}

\def\putdashhvector(#1,#2)#3#4{{%
\arrowlength=#3 \howmanydashes
\multiput(#1,#2)(\increment,0){\numbdashes}%
{\vrule height .4pt width \lengthdash\unitlength}
\arrowtype=#4 \xpos=#1
\ifnum\arrowtype<0 \advance\arrowtype by 7 \fi
\ifcase\arrowtype
\or \advance\xpos by 10
    \put(\xpos,#2){\vector(-1,0){\lengthdash}}
    \advance\xpos by 40
    \put(\xpos,#2){\vector(-1,0){\lengthdash}}
\or \advance \xpos by 10
    \put(\xpos,#2){\vector(-1,0){\lengthdash}}
    \advance\xpos by  \arrowlength
    \advance\xpos by  -50
    \put(\xpos,#2){\vector(-1,0){\lengthdash}}
\or \advance\xpos by 10
    \put(\xpos,#2){\vector(-1,0){\lengthdash}}
\or \advance\xpos by \arrowlength
    \advance\xpos by -\lengthdash
    \put(\xpos,#2){\vector(1,0){\lengthdash}}
\or {\advance\xpos by 10
    \put(\xpos,#2){\vector(1,0){\lengthdash}}}
    \advance\xpos by \arrowlength
    \advance\xpos by -\lengthdash
    \put(\xpos,#2){\vector(1,0){\lengthdash}}
\or \advance\xpos by \arrowlength
    \advance\xpos by -\lengthdash
    \put(\xpos,#2){\vector(1,0){\lengthdash}}
    \advance\xpos by -40
    \put(\xpos,#2){\vector(1,0){\lengthdash}}
   \fi
}}

\def\putdashvvector(#1,#2)#3#4{{%
\arrowlength=#3 \howmanydashes
\ypos=#2 \advance\ypos by -\arrowlength
\multiput(#1,#2)(0,\increment){\numbdashes}%
    {\vrule width .4pt height \lengthdash\unitlength}
\arrowtype=#4 \ypos=#2
\ifnum\arrowtype<0 \advance\arrowtype by 7 \fi
\ifcase\arrowtype
\or \advance\ypos by \arrowlength \advance\ypos by -40
    \put(#1,\ypos){\vector(0,1){\lengthdash}}
    \advance\ypos by -40
    \put(#1,\ypos){\vector(0,1){\lengthdash}}
\or \advance\ypos by 10
    \put(#1,\ypos){\vector(0,1){\lengthdash}}
    \advance\ypos by \arrowlength \advance\ypos by -40
    \put(#1,\ypos){\vector(0,1){\lengthdash}}
\or \advance\ypos by \arrowlength \advance\ypos by -40
    \put(#1,\ypos){\vector(0,1){\lengthdash}}
\or \advance\ypos by 10
    \put(#1,\ypos){\vector(0,-1){\lengthdash}}
\or \advance\ypos by 10
    \put(#1,\ypos){\vector(0,-1){\lengthdash}}
    \advance\ypos by \arrowlength \advance\ypos by -40
    \put(#1,\ypos){\vector(0,-1){\lengthdash}}
\or \advance\ypos by 10
    \put(#1,\ypos){\vector(0,-1){\lengthdash}}
    \advance\ypos by 40
    \put(#1,\ypos){\vector(0,-1){\lengthdash}}
\fi
}}

\def\puthmorphism(#1,#2)[#3`#4`#5]#6#7#8{{%
\xpos #1
\ypos #2
\width #6
\arrowlength #6
\arrowtype=#7
\putbox(\xpos,\ypos){#3\vphantom{#4}}%
{\advance \xpos by\arrowlength
\putbox(\xpos,\ypos){\vphantom{#3}#4}}%
\horsize{\tempcounta}{#3}%
\horsize{\tempcountb}{#4}%
\divide \tempcounta by2
\divide \tempcountb by2
\advance \tempcounta by30
\advance \tempcountb by30
\advance \xpos by\tempcounta
\advance \arrowlength by-\tempcounta
\advance \arrowlength by-\tempcountb
\putvector(\xpos,\ypos)(1,0)\arrowlength\arrowtype
\divide \arrowlength by2
\advance \xpos by\arrowlength
\vertsize{\tempcounta}{#5}%
\divide\tempcounta by2
\advance \tempcounta by20
\if a#8 %
   \advance \ypos by\tempcounta
   \putbox(\xpos,\ypos){#5}%
\else
   \advance \ypos by-\tempcounta
   \putbox(\xpos,\ypos){#5}%
\fi}}

\def\putvmorphism(#1,#2)[#3`#4`#5]#6#7#8{{%
\xpos #1
\ypos #2
\arrowlength #6
\arrowtype #7
\settowidth{\xlen}{$#5$}%
\putbox(\xpos,\ypos){#3}%
{\advance \ypos by-\arrowlength
\putbox(\xpos,\ypos){#4}}%
{\advance\arrowlength by-140
\advance \ypos by-70
\ifdim\xlen>0pt
   \if m#8%
      \putsplitvector(\xpos,\ypos)\arrowlength\arrowtype
   \else
   \putvector(\xpos,\ypos)(0,-1)\arrowlength\arrowtype
   \fi
\else
   \putvector(\xpos,\ypos)(0,-1)\arrowlength\arrowtype
\fi}%
\ifdim\xlen>0pt
   \divide \arrowlength by2
   \advance\ypos by-\arrowlength
   \if l#8%
      \advance \xpos by-40
      \putrbox(\xpos,\ypos){#5}%
   \else\if r#8%
      \advance \xpos by40
      \putlbox(\xpos,\ypos){#5}%
   \else
      \putbox(\xpos,\ypos){#5}%
   \fi\fi
\fi
}}

\def\putsquarep<#1>(#2)[#3;#4`#5`#6`#7]{{%
\setsqparms[#1]%
\setpos(#2)%
\settokens`#3`%
\puthmorphism(\xpos,\ypos)[\tokenc`\tokend`{#7}]{\width}{\arrowtyped}b%
\advance\ypos by \height
\puthmorphism(\xpos,\ypos)[\tokena`\tokenb`{#4}]{\width}{\arrowtypea}a%
\putvmorphism(\xpos,\ypos)[``{#5}]{\height}{\arrowtypeb}l%
\advance\xpos by \width
\putvmorphism(\xpos,\ypos)[``{#6}]{\height}{\arrowtypec}r%
}}

\def\putsquare{\@ifnextchar <{\putsquarep}{\putsquarep%
   <\arrowtypea`\arrowtypeb`\arrowtypec`\arrowtyped;\width`\height>}}
\def\square{\@ifnextchar< {\squarep}{\squarep
   <\arrowtypea`\arrowtypeb`\arrowtypec`\arrowtyped;\width`\height>}}
\def\squarep<#1>[#2`#3`#4`#5;#6`#7`#8`#9]{{
\setsqparms[#1]
\diagram
\putsquarep<\arrowtypea`\arrowtypeb`\arrowtypec`
\arrowtyped;\width`\height>
(0,0)[#2`#3`#4`{#5};#6`#7`#8`{#9}]
\enddiagram
}}                                                 
\def\putptrianglep<#1>(#2,#3)[#4`#5`#6;#7`#8`#9]{{%
\settriparms[#1]%
\xpos=#2 \ypos=#3
\advance\ypos by \height
\puthmorphism(\xpos,\ypos)[#4`#5`{#7}]{\height}{\arrowtypea}a%
\putvmorphism(\xpos,\ypos)[`#6`{#8}]{\height}{\arrowtypeb}l%
\advance\xpos by\height
\putmorphism(\xpos,\ypos)(-1,-1)[``{#9}]{\height}{\arrowtypec}r%
}}

\def\putptriangle{\@ifnextchar <{\putptrianglep}{\putptrianglep
   <\arrowtypea`\arrowtypeb`\arrowtypec;\height>}}
\def\ptriangle{\@ifnextchar <{\ptrianglep}{\ptrianglep
   <\arrowtypea`\arrowtypeb`\arrowtypec;\height>}}
\def\ptrianglep<#1>[#2`#3`#4;#5`#6`#7]{{
\settriparms[#1]
\diagram
\putptrianglep<\arrowtypea`\arrowtypeb`
\arrowtypec;\height>
(0,0)[#2`#3`#4;#5`#6`{#7}]
\enddiagram
}}                                            

\def\putqtrianglep<#1>(#2,#3)[#4`#5`#6;#7`#8`#9]{{%
\settriparms[#1]%
\xpos=#2 \ypos=#3
\advance\ypos by\height
\puthmorphism(\xpos,\ypos)[#4`#5`{#7}]{\height}{\arrowtypea}a%
\putmorphism(\xpos,\ypos)(1,-1)[``{#8}]{\height}{\arrowtypeb}l%
\advance\xpos by\height
\putvmorphism(\xpos,\ypos)[`#6`{#9}]{\height}{\arrowtypec}r%
}}

\def\putqtriangle{\@ifnextchar <{\putqtrianglep}{\putqtrianglep
   <\arrowtypea`\arrowtypeb`\arrowtypec;\height>}}
\def\qtriangle{\@ifnextchar <{\qtrianglep}{\qtrianglep
   <\arrowtypea`\arrowtypeb`\arrowtypec;\height>}}
\def\qtrianglep<#1>[#2`#3`#4;#5`#6`#7]{{
\settriparms[#1]
\width=\height                                
\diagram
\putqtrianglep<\arrowtypea`\arrowtypeb`
\arrowtypec;\height>
(0,0)[#2`#3`#4;#5`#6`{#7}]
\enddiagram
}}

\def\putdtrianglep<#1>(#2,#3)[#4`#5`#6;#7`#8`#9]{{%
\settriparms[#1]%
\xpos=#2 \ypos=#3
\puthmorphism(\xpos,\ypos)[#5`#6`{#9}]{\height}{\arrowtypec}b%
\advance\xpos by \height \advance\ypos by\height
\putmorphism(\xpos,\ypos)(-1,-1)[``{#7}]{\height}{\arrowtypea}l%
\putvmorphism(\xpos,\ypos)[#4``{#8}]{\height}{\arrowtypeb}r%
}}

\def\putdtriangle{\@ifnextchar <{\putdtrianglep}{\putdtrianglep
   <\arrowtypea`\arrowtypeb`\arrowtypec;\height>}}
\def\dtriangle{\@ifnextchar <{\dtrianglep}{\dtrianglep
   <\arrowtypea`\arrowtypeb`\arrowtypec;\height>}}
\def\dtrianglep<#1>[#2`#3`#4;#5`#6`#7]{{
\settriparms[#1]
\width=\height                                
\diagram
\putdtrianglep<\arrowtypea`\arrowtypeb`
\arrowtypec;\height>
(0,0)[#2`#3`#4;#5`#6`{#7}]
\enddiagram
}}

\def\putbtrianglep<#1>(#2,#3)[#4`#5`#6;#7`#8`#9]{{%
\settriparms[#1]%
\xpos=#2 \ypos=#3
\puthmorphism(\xpos,\ypos)[#5`#6`{#9}]{\height}{\arrowtypec}b%
\advance\ypos by\height
\putmorphism(\xpos,\ypos)(1,-1)[``{#8}]{\height}{\arrowtypeb}r%
\putvmorphism(\xpos,\ypos)[#4``{#7}]{\height}{\arrowtypea}l%
}}

\def\putbtriangle{\@ifnextchar <{\putbtrianglep}{\putbtrianglep
   <\arrowtypea`\arrowtypeb`\arrowtypec;\height>}}
\def\btriangle{\@ifnextchar <{\btrianglep}{\btrianglep
   <\arrowtypea`\arrowtypeb`\arrowtypec;\height>}}
\def\btrianglep<#1>[#2`#3`#4;#5`#6`#7]{{
\settriparms[#1]
\width=\height                               
\diagram
\putbtrianglep<\arrowtypea`\arrowtypeb`
\arrowtypec;\height>
(0,0)[#2`#3`#4;#5`#6`{#7}]
\enddiagram
}}

\def\putAtrianglep<#1>(#2,#3)[#4`#5`#6;#7`#8`#9]{{%
\settriparms[#1]%
\xpos=#2 \ypos=#3
{\multiply \height by2
\puthmorphism(\xpos,\ypos)[#5`#6`{#9}]{\height}{\arrowtypec}b}%
\advance\xpos by\height \advance\ypos by\height
\putmorphism(\xpos,\ypos)(-1,-1)[#4``{#7}]{\height}{\arrowtypea}l%
\putmorphism(\xpos,\ypos)(1,-1)[``{#8}]{\height}{\arrowtypeb}r%
}}

\def\putAtriangle{\@ifnextchar <{\putAtrianglep}{\putAtrianglep
   <\arrowtypea`\arrowtypeb`\arrowtypec;\height>}}
\def\Atriangle{\@ifnextchar <{\Atrianglep}{\Atrianglep
   <\arrowtypea`\arrowtypeb`\arrowtypec;\height>}}
\def\Atrianglep<#1>[#2`#3`#4;#5`#6`#7]{{
\settriparms[#1]
\width=\height                                     
\diagram
\putAtrianglep<\arrowtypea`\arrowtypeb`
\arrowtypec;\height>
(0,0)[#2`#3`#4;#5`#6`{#7}]
\enddiagram
}}

\def\putAtrianglepairp<#1>(#2)[#3;#4`#5`#6`#7`#8]{{%
\settripairparms[#1]%
\setpos(#2)%
\settokens`#3`%
\puthmorphism(\xpos,\ypos)[\tokenb`\tokenc`{#7}]{\height}{\arrowtyped}b%
\advance\xpos by\height
\puthmorphism(\xpos,\ypos)[\phantom{\tokenc}`\tokend`{#8}]%
{\height}{\arrowtypee}b%
\advance\ypos by\height
\putmorphism(\xpos,\ypos)(-1,-1)[\tokena``{#4}]{\height}{\arrowtypea}l%
\putvmorphism(\xpos,\ypos)[``{#5}]{\height}{\arrowtypeb}m%
\putmorphism(\xpos,\ypos)(1,-1)[``{#6}]{\height}{\arrowtypec}r%
}}

\def\putAtrianglepair{\@ifnextchar <{\putAtrianglepairp}{\putAtrianglepairp%
   <\arrowtypea`\arrowtypeb`\arrowtypec`\arrowtyped`\arrowtypee;\height>}}
\def\Atrianglepair{\@ifnextchar <{\Atrianglepairp}{\Atrianglepairp%
   <\arrowtypea`\arrowtypeb`\arrowtypec`\arrowtyped`\arrowtypee;\height>}}

\def\Atrianglepairp<#1>[#2;#3`#4`#5`#6`#7]{{
\settripairparms[#1]
\settokens`#2`
\width=\height                                
\diagram
\putAtrianglepairp                            
<\arrowtypea`\arrowtypeb`\arrowtypec`
\arrowtyped`\arrowtypee;\height>
(0,0)[{#2};#3`#4`#5`#6`{#7}]
\enddiagram
}}

\def\putVtrianglep<#1>(#2,#3)[#4`#5`#6;#7`#8`#9]{{%
\settriparms[#1]%
\xpos=#2 \ypos=#3
\advance\ypos by\height
{\multiply\height by2
\puthmorphism(\xpos,\ypos)[#4`#5`{#7}]{\height}{\arrowtypea}a}%
\putmorphism(\xpos,\ypos)(1,-1)[`#6`{#8}]{\height}{\arrowtypeb}l%
\advance\xpos by\height
\advance\xpos by\height
\putmorphism(\xpos,\ypos)(-1,-1)[``{#9}]{\height}{\arrowtypec}r%
}}

\def\putVtriangle{\@ifnextchar <{\putVtrianglep}{\putVtrianglep
   <\arrowtypea`\arrowtypeb`\arrowtypec;\height>}}
\def\Vtriangle{\@ifnextchar <{\Vtrianglep}{\Vtrianglep
   <\arrowtypea`\arrowtypeb`\arrowtypec;\height>}}
\def\Vtrianglep<#1>[#2`#3`#4;#5`#6`#7]{{
\settriparms[#1]
\width=\height                                 
\diagram
\putVtrianglep<\arrowtypea`\arrowtypeb`
\arrowtypec;\height>
(0,0)[#2`#3`#4;#5`#6`{#7}]
\enddiagram
}}

\def\putVtrianglepairp<#1>(#2)[#3;#4`#5`#6`#7`#8]{{
\settripairparms[#1]%
\setpos(#2)%
\settokens`#3`%
\advance\ypos by\height
\putmorphism(\xpos,\ypos)(1,-1)[`\tokend`{#6}]{\height}{\arrowtypec}l%
\puthmorphism(\xpos,\ypos)[\tokena`\tokenb`{#4}]{\height}{\arrowtypea}a%
\advance\xpos by\height
\puthmorphism(\xpos,\ypos)[\phantom{\tokenb}`\tokenc`{#5}]%
{\height}{\arrowtypeb}a%
\putvmorphism(\xpos,\ypos)[``{#7}]{\height}{\arrowtyped}m%
\advance\xpos by\height
\putmorphism(\xpos,\ypos)(-1,-1)[``{#8}]{\height}{\arrowtypee}r%
}}

\def\putVtrianglepair{\@ifnextchar <{\putVtrianglepairp}{\putVtrianglepairp%
    <\arrowtypea`\arrowtypeb`\arrowtypec`\arrowtyped`\arrowtypee;\height>}}
\def\Vtrianglepair{\@ifnextchar <{\Vtrianglepairp}{\Vtrianglepairp%
    <\arrowtypea`\arrowtypeb`\arrowtypec`\arrowtyped`\arrowtypee;\height>}}
\def\Vtrianglepairp<#1>[#2;#3`#4`#5`#6`#7]{{
\settripairparms[#1]
\settokens`#2`
\diagram
\putVtrianglepairp                             
<\arrowtypea`\arrowtypeb`\arrowtypec`
\arrowtyped`\arrowtypee;\height>
(0,0)[{#2};#3`#4`#5`#6`{#7}]
\enddiagram
}}

\def\putCtrianglep<#1>(#2,#3)[#4`#5`#6;#7`#8`#9]{{%
\settriparms[#1]%
\xpos=#2 \ypos=#3
\advance\ypos by\height
\putmorphism(\xpos,\ypos)(1,-1)[``{#9}]{\height}{\arrowtypec}l%
\advance\xpos by\height
\advance\ypos by\height
\putmorphism(\xpos,\ypos)(-1,-1)[#4`#5`{#7}]{\height}{\arrowtypea}l%
{\multiply\height by 2
\putvmorphism(\xpos,\ypos)[`#6`{#8}]{\height}{\arrowtypeb}r}%
}}

\def\putCtriangle{\@ifnextchar <{\putCtrianglep}{\putCtrianglep
    <\arrowtypea`\arrowtypeb`\arrowtypec;\height>}}
\def\Ctriangle{\@ifnextchar <{\Ctrianglep}{\Ctrianglep
    <\arrowtypea`\arrowtypeb`\arrowtypec;\height>}}
\def\Ctrianglep<#1>[#2`#3`#4;#5`#6`#7]{{
\settriparms[#1]
\width=\height                               
\diagram
\putCtrianglep<\arrowtypea`\arrowtypeb`
\arrowtypec;\height>
(0,0)[#2`#3`#4;#5`#6`{#7}]
\enddiagram
}}                                           
\def\putDtrianglep<#1>(#2,#3)[#4`#5`#6;#7`#8`#9]{{%
\settriparms[#1]%
\xpos=#2 \ypos=#3
\advance\xpos by\height \advance\ypos by\height
\putmorphism(\xpos,\ypos)(-1,-1)[``{#9}]{\height}{\arrowtypec}r%
\advance\xpos by-\height \advance\ypos by\height
\putmorphism(\xpos,\ypos)(1,-1)[`#5`{#8}]{\height}{\arrowtypeb}r%
{\multiply\height by 2
\putvmorphism(\xpos,\ypos)[#4`#6`{#7}]{\height}{\arrowtypea}l}%
}}

\def\putDtriangle{\@ifnextchar <{\putDtrianglep}{\putDtrianglep
    <\arrowtypea`\arrowtypeb`\arrowtypec;\height>}}
\def\Dtriangle{\@ifnextchar <{\Dtrianglep}{\Dtrianglep
   <\arrowtypea`\arrowtypeb`\arrowtypec;\height>}}
\def\Dtrianglep<#1>[#2`#3`#4;#5`#6`#7]{{
\settriparms[#1]
\width=\height                              
\diagram
\putDtrianglep<\arrowtypea`\arrowtypeb`
\arrowtypec;\height>
(0,0)[#2`#3`#4;#5`#6`{#7}]
\enddiagram
}}                                          
\def\setrecparms[#1`#2]{\width=#1 \height=#2}%

\def\recursep<#1`#2>[#3;#4`#5`#6`#7`#8]{{\m@th
\width=#1 \height=#2
\settokens`#3`
\settowidth{\tempdimen}{$\tokena$}
\ifdim\tempdimen=0pt
  \savebox{\tempboxa}{\hbox{$\tokenb$}}%
  \savebox{\tempboxb}{\hbox{$\tokend$}}%
  \savebox{\tempboxc}{\hbox{$#6$}}%
\else
  \savebox{\tempboxa}{\hbox{$\hbox{$\tokena$}\times\hbox{$\tokenb$}$}}%
  \savebox{\tempboxb}{\hbox{$\hbox{$\tokena$}\times\hbox{$\tokend$}$}}%
  \savebox{\tempboxc}{\hbox{$\hbox{$\tokena$}\times\hbox{$#6$}$}}%
\fi
\ypos=\height
\divide\ypos by 2
\xpos=\ypos
\advance\xpos by \width
\bfig
\putCtrianglep<-1`1`1;\ypos>(0,0)[`\tokenc`;#5`#6`{#7}]%
\puthmorphism(\ypos,0)[\tokend`\usebox{\tempboxb}`{#8}]{\width}{-1}b%
\puthmorphism(\ypos,\height)[\tokenb`\usebox{\tempboxa}`{#4}]{\width}{-1}a%
\advance\ypos by \width
\putvmorphism(\ypos,\height)[``\usebox{\tempboxc}]{\height}1r%
\efig
}}

\def\recurse{\@ifnextchar <{\recursep}{\recursep<\width`\height>}}

\def\puttwohmorphisms(#1,#2)[#3`#4;#5`#6]#7#8#9{{%
\puthmorphism(#1,#2)[#3`#4`]{#7}0a
\ypos=#2
\advance\ypos by 20
\puthmorphism(#1,\ypos)[\phantom{#3}`\phantom{#4}`#5]{#7}{#8}a
\advance\ypos by -40
\puthmorphism(#1,\ypos)[\phantom{#3}`\phantom{#4}`#6]{#7}{#9}b
}}

\def\puttwovmorphisms(#1,#2)[#3`#4;#5`#6]#7#8#9{{%
\putvmorphism(#1,#2)[#3`#4`]{#7}0a
\xpos=#1
\advance\xpos by -20
\putvmorphism(\xpos,#2)[\phantom{#3}`\phantom{#4}`#5]{#7}{#8}l
\advance\xpos by 40
\putvmorphism(\xpos,#2)[\phantom{#3}`\phantom{#4}`#6]{#7}{#9}r
}}

\def\puthcoequalizer(#1)[#2`#3`#4;#5`#6`#7]#8#9{{%
\setpos(#1)%
\puttwohmorphisms(\xpos,\ypos)[#2`#3;#5`#6]{#8}11%
\advance\xpos by #8
\puthmorphism(\xpos,\ypos)[\phantom{#3}`#4`#7]{#8}1{#9}
}}

\def\putvcoequalizer(#1)[#2`#3`#4;#5`#6`#7]#8#9{{%
\setpos(#1)%
\puttwovmorphisms(\xpos,\ypos)[#2`#3;#5`#6]{#8}11%
\advance\ypos by -#8
\putvmorphism(\xpos,\ypos)[\phantom{#3}`#4`#7]{#8}1{#9}
}}

\def\putthreehmorphisms(#1)[#2`#3;#4`#5`#6]#7(#8)#9{{%
\setpos(#1) \settypes(#8)
\if a#9 %
     \vertsize{\tempcounta}{#5}%
     \vertsize{\tempcountb}{#6}%
     \ifnum \tempcounta<\tempcountb \tempcounta=\tempcountb \fi
\else
     \vertsize{\tempcounta}{#4}%
     \vertsize{\tempcountb}{#5}%
     \ifnum \tempcounta<\tempcountb \tempcounta=\tempcountb \fi
\fi
\advance \tempcounta by 60
\puthmorphism(\xpos,\ypos)[#2`#3`#5]{#7}{\arrowtypeb}{#9}
\advance\ypos by \tempcounta
\puthmorphism(\xpos,\ypos)[\phantom{#2}`\phantom{#3}`#4]{#7}{\arrowtypea}{#9}
\advance\ypos by -\tempcounta \advance\ypos by -\tempcounta
\puthmorphism(\xpos,\ypos)[\phantom{#2}`\phantom{#3}`#6]{#7}{\arrowtypec}{#9}
}}

\def\setarrowtoks[#1`#2`#3`#4`#5`#6]{%
\def\toka{#1}
\def\tokb{#2}
\def\tokc{#3}
\def\tokd{#4}
\def\toke{#5}
\def\tokf{#6}
}
\def\hex{\@ifnextchar <{\hexp}{\hexp<1000`400>}}
\def\hexp<#1`#2>[#3`#4`#5`#6`#7`#8;#9]{%
\setarrowtoks[#9]
\yext=#2 \advance \yext by #2
\xext=#1 \advance\xext by \yext
\bfig
\putCtriangle<-1`0`1;#2>(0,0)[`#5`;\tokb``\tokd]
\xext=#1 \yext=#2 \advance \yext by #2
\putsquare<1`0`0`1;\xext`\yext>(#2,0)[#3`#4`#7`#8;\toka```\tokf]
\advance \xext by #2
\putDtriangle<0`1`-1;#2>(\xext,0)[`#6`;`\tokc`\toke]
\efig
}
\makeatother